\documentclass[11pt]{article}

\usepackage{color}
\usepackage{latexsym}
\usepackage{amssymb}
\usepackage{amsmath,amsfonts,theorem,euscript,array,enumerate,mathrsfs}
\usepackage{graphicx}
\newtheorem{Theorem}{Theorem}[part]
\newtheorem{Definition}{Definition}[part]
\newtheorem{Proposition}{Proposition}[part]

\newtheorem{Lemma}{Lemma}[part]
\newtheorem{Corollary}{Corollary}[part]
\newtheorem{Remark}{Remark}[part]

\def\esssup_#1{\underset{#1}{\mathrm{ess\,sup\, }}}
\def\essinf_#1{\underset{#1}{\mathrm{ess\,inf\, }}}

\def \trans{^{\scriptscriptstyle{\intercal}}}

\def \trans{^{\scriptscriptstyle{\intercal }}}

\def \A{\mathbb{A}}

\def \N{\mathbb{N}}
\def \R{\mathbb{R}}

\def \E{\mathbb{E}}
\def \F{\mathbb{F}}

\def \H{\mathbb{H}}
\def \L{\mathbb{L}}
\def \P{\mathbb{P}}
\def \Q{\mathbb{Q}}
\def \D{\mathbb{D}}
\def \S{\mathbb{S}}

\def \X{\mathbb{X}}

\def \Bc{{\cal B}}

\def \Dc{{\cal D}}

\def \Fc{{\cal F}}

\def \Mc{{\cal M}}

\def \Sc{{\cal S}}

\def \Uc{{\cal U}}

\def \eps{\varepsilon}

\def \ep{\hbox{ }\hfill$\Box$}

\allowdisplaybreaks

\begin{document}

\title{A regularization approach to functional It\^o calculus and strong-viscosity solutions to path-dependent PDEs}

\author{
Andrea COSSO\thanks{Laboratoire de Probabilit\'es et Mod\`eles Al\'eatoires, CNRS, UMR 7599, Universit\'e Paris Diderot, France. E-mail: \sf andrea.cosso@polimi.it}
\qquad\quad
Francesco RUSSO\thanks{ENSTA ParisTech, Unit\'e de Math\'ematiques appliqu\'ees, 828, boulevard des Mar\'echaux, F-91120 Palaiseau, France. E-mail: \sf francesco.russo@ensta-paristech.fr}
}

\maketitle


\begin{abstract}
First, we revisit functional It\^o/path-dependent calculus started by B. Dupire, R. Cont and D.-A. Fourni\'e, using the formulation of calculus via regularization. Relations with the corresponding Banach space valued calculus 
introduced by C. Di Girolami and the second named author 
are explored. The second part of the paper is devoted to the study of the Kolmogorov type equation associated with the so called window Brownian motion, called path-dependent heat equation, for which well-posedness at the level of classical solutions is established. Then, a notion of strong approximating solution, called strong-viscosity solution, is introduced which is supposed to be a substitution tool to the viscosity solution. For that kind of solution, we also prove existence and uniqueness. The notion of strong-viscosity solution motivates the last part of the paper which is devoted to explore this new concept of solution for general semilinear PDEs in the finite dimensional case. We prove an equivalence result between the classical viscosity solution and the new one. The definition of strong-viscosity solution for semilinear PDEs is inspired by the notion of {\it good} solution, and it is based again on an approximating procedure. 
\end{abstract}

\vspace{5mm}

\noindent {\bf Key words:} Horizontal and vertical derivative; functional It\^o/path-dependent calculus; 
 strong-viscosity solutions; calculus via regularization.

\vspace{5mm}

\noindent {\bf 2010 Math Subject Classification:}  35D35; 35D40; 35K10;
 60H05;  60H10; 60H30.

\newpage

\section{Introduction}

\setcounter{equation}{0} \setcounter{Assumption}{0}
\setcounter{Theorem}{0} \setcounter{Proposition}{0}
\setcounter{Corollary}{0} \setcounter{Lemma}{0}
\setcounter{Definition}{0} \setcounter{Remark}{0}

Recently, 
 a new branch of  stochastic calculus has appeared, known as \emph{functional It\^o calculus}, which results to be an extension of classical It\^o calculus to 
functionals depending on the all path of a stochastic process and not only on its current value, see Dupire \cite{dupire}, Cont and Fourni\'e 
\cite{contfournie10,contfournie,contfournie13}.
On the other hand,  C. Di Girolami, the second named author
 and more recently
G. Fabbri, have introduced in a series of papers (\cite{DGR, DGRnote, 
digirrusso12, DGR2, digirfabbrirusso13}), a stochastic calculus
via regularizations for processes taking values in
a separable Banach space $B$, which includes (when $B = C([-T,0]$), as 
applications a path-dependent type calculus having similar objectives.

In the first part of the present paper, we revisit functional It\^o calculus by means of stochastic calculus via regularization. We recall that 
developed functional It\^o calculus and derived a functional It\^o's formula 
using discretization techniques of F\"ollmer \cite{follmer} type, instead of regularization techniques. Let us illustrate another difference with respect to \cite{contfournie10}. One of the main issues of functional It\^o calculus is the definition of the functional (or pathwise) derivatives, i.e., the horizontal derivative (calling in only the past values of the trajectory) and the vertical derivative (calling in only the present value of the trajectory). 
In \cite{contfournie10}, it is essential to consider functionals defined on the space of c\`adl\`ag trajectories, since the definition of functional derivatives necessitates of discontinuous paths. Therefore, if a functional is defined only on the space of continuous trajectories (because, e.g., it depends on the paths of a continuous process as Brownian motion), we have to extend it anyway to the space of c\`adl\`ag trajectories, even though, in general, there is no 
 unique way to extend it. In contrast to this approach, we introduce an intermediate space between the space of continuous trajectories $C([-T,0])$ and the space of c\`adl\`ag trajectories $\D([-T,0])$, denoted $\mathscr C([-T,0])$, which allows us to define functional derivatives. $\mathscr C([-T,0])$ is the space of bounded trajectories on $[-T,0]$, continuous on $[-T,0[$ and with possibly a jump at $0$. We endow $\mathscr C([-T,0])$ with a topology such that $C([-T,0])$ is dense in $\mathscr C([-T,0])$ with respect to this topology. Therefore, any functional $\Uc\colon[0,T]\times C([-T,0])\rightarrow\R$, continuous with respect to the topology of $\mathscr C([-T,0])$, admits a unique extension to $\mathscr C([-T,0])$, denoted $u\colon[0,T]\times\mathscr C([-T,0])\rightarrow\R$. We present some significant functionals for which a continuous extension exists. Then, we develop the functional It\^o calculus for $u\colon[0,T]\times\mathscr C([-T,0])\rightarrow\R$.

Notice that we use a slightly different notation with respect to \cite{contfournie10}. In particular, in place of a map $\Uc\colon [0,T]\times C([-T,0])\rightarrow\R$, in \cite{contfournie10} a family of maps $F=(F_t)_{t\in[0,T]}$, with $F_t\colon C([0,t])\rightarrow\R$, is considered. However, we can always move from one formulation to the other. Indeed, given $F=(F_t)_{t\in[0,T]}$, where each $F_t\colon C([0,t])\rightarrow\R$, we can define $\Uc\colon [0,T]\times C([-T,0])\rightarrow\R$ as follows:
\[
\Uc(t,\eta) \ := \ F_t(\eta(\cdot+T)|_{[0,t]}), \qquad (t,\eta)\in[0,T]\times C([-T,0]).
\]
Vice-versa, let $\Uc\colon [0,T]\times C([-T,0])\rightarrow\R$ and define $F=(F_t)_{t\in[0,T]}$ as
\begin{equation}
\label{E:F=u}
F_t(\tilde\eta) \ := \ \Uc(t,\eta), \qquad (t,\tilde\eta)\in[0,T]\times C([0,t]),
\end{equation}
where $\eta$ is the element of $C([-T,0])$ obtained from $\tilde\eta$ firstly translating $\tilde\eta$ on the interval $[-t,0]$, then extending it in a constant way up to $-T$, namely $\eta(x) := \tilde\eta(x+t)1_{[-t,0]}(x) + \tilde\eta(-t)1_{[-T,-t)}(x)$, for any $x\in[-T,0]$. Observe that, in principle, the map $\Uc$ contains more information than $F$, since in \eqref{E:F=u} we do not take into account the values of $\Uc$ at $(t,\eta)\in[0,T]\times C([-T,0])$ with $\eta$ not constant on the interval $[-T,-t]$. Despite this, the equivalence between the two notations is guaranteed by the fact that, as it will be clear later, when we consider the composition of $\Uc$ with a stochastic process, this extra information plays no role. Our formulation has two advantages. Firstly, we can work with a single map instead of a family of maps. In addition, the time variable and the path have two distinct roles in our setting, as for the time variable and the space variable in the classical It\^o calculus. This, in particular, allows us to define the horizontal derivative independently of the time derivative, so that, the horizontal derivative defined in \cite{contfournie10} corresponds to the sum of our horizontal derivative and of the time derivative. We mention that an alternative approach to functional derivatives was introduced in \cite{buckdahn_ma_zhang13}.

We end the first part of the paper showing how our functional It\^o's formula is strictly related to the It\^o's formula derived in the framework of Banach space valued stochastic calculus via regularization, for the case of window processes. This latter and brand new branch of stochastic calculus and stochastic analysis has been recently conceived, deeply studied, and developed in many directions in  \cite{DGR2, digirrusso12, DGRnote},
\cite{digirfabbrirusso13} and for more details \cite{DGR}. 
For the particular case of window processes, we also refer to Theorem 6.3 and Section 7.2 in \cite{digirfabbrirusso13}. In the present paper, we prove formulae which allow to express functional derivatives in terms of differential operators arising in the Banach space valued stochastic calculus via regularization, with the aim of identifying the building blocks of our functional It\^o's formula with the terms appearing in the It\^o's formula for window processes.

Dupire \cite{dupire} introduced also the concept of \emph{path-dependent partial differential equation}, to which the second part of the present paper is devoted. Di Girolami and the second named author, in Chapter 9
of \cite{DGR}, considered
a similar equation in the framework of Banach space valued calculus, for which we refer also to \cite{flandoli_zanco13}. 
We focus on path-dependent nonlinear Kolmogorov equations driven by the path-dependent heat operator, for which we provide a definition of classical solution. We prove a uniqueness result for this kind of solution, by means of probabilistic methods based on the theory of backward stochastic differential equations (BSDEs). More precisely, we show that, if a classical solution exists, then it can be expressed through the  solution of
 a certain backward stochastic differential equation. Therefore, from the uniqueness of the BSDE it follows that there exists at most one classical solution. Then, we prove an existence result for classical solutions. However, this notion of solution turns out to be unsuitable to deal with all significant examples. As a matter of fact, if we consider the path-dependent PDE arising in the hedging problem of lookback contingent claims, we can not expect too much regularity of the solution (this example is studied in detail in subsection \ref{SubS:HedgingExample}). Therefore, we are led to consider a weaker notion of solution to the path-dependent nonlinear Kolmogorov equation. In particular, we are interested in a viscosity-type solution, namely a solution which is not required to be differentiable, but only locally uniformly continuous.

The issue of providing a suitable definition of viscosity solutions for path-dependent PDEs has attracted a great interest. We recall that Ekren, Keller, Touzi, and Zhang~\cite{ektz} and Ekren, Touzi, and Zhang~\cite{etzI,etzII} recently provided a definition of viscosity solution to path-dependent PDEs, replacing the classical minimum/maximum property, which appears in the standard definition of viscosity solution, with an optimal stopping problem under nonlinear expectation~\cite{etzOptStop}. We also recall that other definitions of viscosity solutions for path-dependent PDEs were given by Peng~\cite{peng12} and Tang and Zhang~\cite{tangzhang13}. In contrast with the above cited papers, our definition of solution is not inspired by the standard definition of viscosity solution given in terms of test functions or jets.
In fact, our weak solution, called \emph{strong-viscosity solution}, to the path-dependent nonlinear Kolmogorov equation is defined, in a few words, as the pointwise limit of classical solutions to perturbed equations. We notice that our definition is more similar in spirit to the concept of 
\emph{good solution}, which turned out to be equivalent to the definition of $L^p$-viscosity solution for certain fully nonlinear partial differential equations, see, e.g., \cite{cerutti_escauriaza_fabes93}, \cite{crandall_kocan_soravia_swiech94}, \cite{jensen96}, and \cite{jensen_kocan_swiech02}. It has also some similarities with the \emph{vanishing viscosity} method, which represents one of the primitive ideas leading to the conception of the modern definition of viscosity solution. Our definition is likewise inspired by the notion of \emph{strong solution} (which justifies the name of our solution), as defined for example in \cite{cerrai01}, \cite{gozzi_russo06a}, and \cite{gozzi_russo06b}, even though strong solutions are required to be more regular than simply locally uniformly continuous (this regularity is usually required to prove uniqueness of strong
 solutions, which for example in \cite{gozzi_russo06a} and \cite{gozzi_russo06b} is based on a Fukushima-Dirichlet decomposition).
 Instead, our definition of strong-viscosity solution to the path-dependent nonlinear Kolmogorov equation is only required to be locally uniformly continuous and with polynomial growth. The term \emph{viscosity} in the name of our solution is also justifies by the fact that in the finite dimensional case we have an equivalence result between the notion of strong-viscosity solution and that of viscosity solution. We prove a uniqueness theorem for strong-viscosity solutions using the theory of backward stochastic differential equations and we provide an existence result.

We conclude the second part of the paper analyzing more in detail the notion of strong-viscosity solution for semilinear equations, focusing on the more understandable finite dimensional case. In particular, we notice that, if from one hand we consider only equations, as well as perturbed equations, whose classical solutions admit a nonlinear Feynman-Kac representation formula in terms of BSDEs, on the other hand our definition of strong-viscosity solution has the advantage that the comparison theorem follows directly from the comparison theorem for BSDEs. In other words, the comparison theorem for strong-viscosity solutions can be proved using probabilistic methods, in contrast to real analysis' tools which characterize comparison theorems for viscosity solutions and revealed to be arduous to extend to the infinite dimensional setting, see, e.g., \cite{fabbrigozziswiech14}. We present two definitions of strong-viscosity solution, one of them is more in the spirit of the standard definition of viscosity 
solution, since it is required to be both a strong-viscosity subsolution and a strong-viscosity supersolution. A strong-viscosity supersolution (resp. subsolution) is defined, in few words, to be the pointwise limit of classical supersolutions (resp. subsolutions) to perturbed semilinear equations. We prove, using the theory of backward stochastic differential equations, that a comparison theorem for strong-viscosity sub and supersolutions holds, therefore obtaining a uniqueness result for our strong-viscosity solution. More precisely, we prove that every strong-viscosity supersolution (resp. subsolution) can be represented by a supersolution (resp. subsolution) of a BSDE. Indeed, every strong-viscosity supersolution is the limit of a sequence of classical supersolutions, which admit a representation in terms of supersolutions of BSDEs. Then, using a limit theorem for BSDEs (partly inspired by the monotonic limit theorem  of Peng~\cite{peng00}), we derive a limit BSDE supersolution, which turns out to be a probabilistic representation for our strong-viscosity supersolution. Therefore, as already mentioned, the comparison theorem for strong-viscosity sub and supersolutions is a consequence of the comparison theorem for BSDEs. We conclude investigating the equivalence between the notion of strong-viscosity solution and the standard definition of viscosity solution. 

The rest of the paper is organized as follows. In section \ref{S:Functional} we develop functional It\^o calculus via regularization: after a brief introduction on finite dimensional stochastic calculus via regularization in subsection \ref{SubS:Background}, we introduce and study the space $\mathscr C([-T,0])$ in subsection \ref{SubS:Preliminaries}; then, we define the pathwise derivatives and we prove the functional It\^o's formula in subsection \ref{SubS:PathwiseDerivatives}; in subsection \ref{SubS:Comparison}, instead, we discuss the relation between functional It\^o calculus via regularization and Banach space valued stochastic calculus via regularization for window processes. In section \ref{S:StrongViscositySolutions}, on the other hand, we study path-dependent PDEs. More precisely, in subsection \ref{SubS:ClassicalSolutions} we discuss classical solutions to the path-dependent nonlinear Kolmogorov equation; in subsection \ref{SubS:HedgingExample} we present a significant hedging example to motivate the introduction of a weaker notion of solution; in subsection \ref{SubS:StringViscositySolutions} we provide the definition of strong-viscosity solution to the path-dependent nonlinear Kolmogorov equation; finally, in subsection \ref{SubS:FiniteDimensionalCase} we explore the notion of strong-viscosity solution for more general PDEs (not only PDEs driven by the heat operator), in the finite dimensional case.

\section{Functional It\^o calculus: a regularization approach}
\label{S:Functional}

\setcounter{equation}{0} \setcounter{Assumption}{0}
\setcounter{Theorem}{0} \setcounter{Proposition}{0}
\setcounter{Corollary}{0} \setcounter{Lemma}{0}
\setcounter{Definition}{0} \setcounter{Remark}{0}

\subsection{Background: finite dimensional calculus via regularization}
\label{SubS:Background}




The theory of stochastic calculus via regularization has been developed in several papers, starting from \cite{russovallois91, russovallois93}. We recall below only the results used in the present paper, and we refer to \cite{russovallois07} for a survey on the subject. 
 We emphasize that integrands are allowed to be anticipating. Moreover, the integration theory and calculus appears to be close to a pure pathwise approach even though there is still a probability space behind.

\vspace{3mm}

Fix a probability space $(\Omega,\Fc,\P)$ and $T\in]0,\infty[$. Let $\F=(\Fc_t)_{t\in[0,T]}$ denote a filtration satisfying the usual conditions. Let $X=(X_t)_{t\in[0,T]}$ (resp. $Y=(Y_t)_{t\in[0,T]}$) be a real continuous (resp. $\P$-a.s. integrable) process. Every real continuous process $X = (X_t)_{t\in[0,T]}$ is naturally extended to all $t\in\R$ setting $X_t = X_0$, $t\leq0$, and $X_t=X_T$, $t\geq T$. We also define a $C([-T,0])$-valued process $\X=(\X_t)_{t\in\R}$, called the {\bf window process} associated with $X$, defined by
\[
\X_t := \{X_{t+x},\,x\in[-T,0]\}, \qquad t\in\R.
\]

\begin{Definition}
Suppose that, for every $t\in[0,T]$,  the following limit
\begin{equation}
\label{DPFI}
\int_{0}^{t}Y_sd^-X_s \ := \ \lim_{\eps\rightarrow0^+}\int_{0}^{t} Y_s\frac{X_{s+\eps}-X_s}{\eps}ds,
\end{equation}
exists in probability. If the obtained  random function admits a continuous modification, that process is denoted by  $\int_0^\cdot Yd^-X$ and called \textbf{forward integral of $Y$ with respect to $X$}.
\end{Definition}

\begin{Definition}
A family of processes $(H_t^{(\eps)})_{t\in[0,T]}$ is said to converge to $(H_t)_{t\in[0,T]}$ in the \textbf{ucp sense}, if $\sup_{0\leq t\leq T}|H_t^{(\eps)}-H_t|$ goes to $0$ in probability, as $\eps\rightarrow0^+$.
\end{Definition}

\begin{Proposition}
Suppose that the limit \eqref{DPFI} exists in the ucp sense. Then, the forward integral $\int_0^\cdot Yd^-X$ of $Y$ with respect to $X$ exists.
\end{Proposition}

Let us introduce the concept of covariation, which is a crucial notion in stochastic calculus via regularization. Let us suppose that
$X, Y$ are continuous processes.

\begin{Definition}
\label{D33}
The \textbf{covariation of $X$ and $Y$} is defined by
\[						
\left[X,Y\right]_{t} \ = \ \left[Y,X\right]_{t}  \ = \ \lim_{\eps\rightarrow 0^{+}}
\frac{1}{\eps} \int_{0}^{t} (X_{s+\eps}-X_{s})(Y_{s+\eps}-Y_{s})ds, \qquad t \in [0,T],
\]
if the limit exists in probability for every  $t \in [0,T]$, provided that the limiting random function  admits a continuous version $($this is the case if the limit holds in the ucp sense$)$. If $X=Y,$ $X$ is said to be a \textbf{finite quadratic variation process} and we set $[X]:=[X,X]$. 
\end{Definition}

The forward integral and the covariation generalize the classical It\^o integral and covariation for semimartingales. In particular, we have the following result, for a proof we refer to, e.g., \cite{russovallois07}.

\begin{Proposition}
\label{Pproperties}
The following properties hold:
\begin{enumerate}
\item[\textup{(i)}] Let $S^1,S^2$ be continuous $\F$-semimartingales. Then, $[S^1,S^2]$ is the classical bracket $[S^1,S^2]=\langle M^1,M^2\rangle$, where $M^1$ $($resp. $M^2$$)$ is the local martingale part of $S^1$ $($resp. $S^2$$)$.
\item[\textup{(ii)}] Let $V$ be a continuous bounded variation process and $Y$ be a c\`{a}dl\`{a}g process $($or vice-versa$)$$;$ then $[V] =[Y,V]= 0$. Moreover $\int_0^\cdot Y d^-V=\int_0^\cdot Y dV $,
is the {\bf Lebesgue-Stieltjes integral}.
\item[\textup{(iii)}] If $W$ is a Brownian motion and $Y$ is an $\F$-progressively measurable process such that $\int_0^T Y^2_s ds < \infty$, $\P$-a.s., then $\int_0^\cdot Yd^- W$ exists and equals the It\^o integral  $\int_0^\cdot YdW$.
\end{enumerate}
\end{Proposition}

We end this crash introduction to finite dimensional stochastic calculus via regularization presenting one of its cornerstones: It\^o's formula. It is a  well-known result in the theory of semimartingales, but it also  extends to the framework of  finite quadratic variation processes. For a proof we refer to Theorem 2.1 of \cite{russovallois95}.

\begin{Theorem}
\label{ITOFQV}
Let $F:[0,T]\times \R\longrightarrow \R$ be of class $ C^{1,2}\left( [0,T]\times \R \right)$ and $X=(X_t)_{t\in[0,T]}$ be
a real continuous finite quadratic variation process. Then, the following \textbf{It\^o's formula} holds, $\P$-a.s.,
\begin{align}
F(t,X_{t}) \ &= \ F(0,X_{0}) + \int_{0}^{t}\partial_t  F(s,X_s)ds + \int_{0}^{t} \partial_{x} F(s,X_s)d^-X_s \notag \\
&\quad \ + \frac{1}{2}\int_{0}^{t} \partial^2_{x\, x} F(s,X_s)d[X]_s, \qquad\qquad\qquad 0\leq t\leq T. \label{E:ITOFQV}
\end{align}
\end{Theorem}

\subsubsection{The deterministic calculus via regularization}

\label{S211}

A useful particular case of finite dimensional stochastic calculus via regularization arises when $\Omega$ is a singleton, i.e., when the calculus becomes deterministic. In addition, in this deterministic framework we will make use of the \emph{definite integral} on an interval $[a,b]$, where $a<b$ are two real numbers. Typically, we will consider $a=-T$ or $a=-t$ and $b=0$.

\vspace{3mm}

We start with two conventions. By default, every bounded variation function $f\colon[a,b]\rightarrow\R$ will be considered as c\`{a}dl\`{a}g. Moreover, if $f\colon[a,b]\rightarrow\R$ is a c\`{a}dl\`{a}g function, we extend it naturally to another c\`{a}dl\`{a}g function $\tilde f$ on the real line as follows:
\[
\tilde f(x) \ = \
\begin{cases}
f(b), \qquad\qquad &x>b, \\
f(x),  &a\leq x\leq b, \\
0,  &x<a.
\end{cases}
\]

\begin{Definition}
\label{D:DeterministicIntegral}
Let $f\colon[a,b]\rightarrow\R$ be a c\`{a}dl\`{a}g function and $g\colon[a,b]\rightarrow\R$ be in $L^1([a,b])$.\\
\textup{(i)} Suppose that the following limit
\[
\int_{[a,b]}g(s)d^-f(s) \ := \ \lim_{\eps\rightarrow0^+}\int_a^b g(s)\frac{f(s+\eps)-f(s)}{\eps}ds,
\]
exists and it is finite. Then, the obtained quantity is denoted by $\int_{[a,b]} gd^-f$ and called \textbf{$($deterministic, definite$)$ forward integral of $g$ with respect to $f$ $($on $[a,b]$$)$}.\\
\textup{(ii)} Suppose that the following limit
\[
\int_{[a,b]}g(s)d^+f(s) \ := \ \lim_{\eps\rightarrow0^+}\int_a^b g(s)\frac{f(s)-f(s-\eps)}{\eps}ds,
\]
exists and it is finite. Then, the obtained quantity is denoted by $\int_{[a,b]} gd^+f$ and called \textbf{$($deterministic, definite$)$ backward integral of $g$ with respect to $f$ $($on $[a,b]$$)$}.
\end{Definition}

\noindent Let us now introduce the deterministic covariation.

\begin{Definition}
\label{D:DeterministicQuadrVar}
Let $f,g\colon[a,b]\rightarrow\R$ be continuous functions and suppose that $0\in[a,b]$. The \textbf{$($deterministic$)$ covariation of $f$ and $g$ $($on $[a,b]$$)$} is defined by
\[						
\left[f,g\right](x) \ = \ \left[g,f\right](x)  \ = \ \lim_{\eps\rightarrow 0^{+}}
\frac{1}{\eps} \int_{0}^x (f(s+\eps)-f(s))(g(s+\eps)-g(s))ds, \qquad x\in[a,b],
\]
if the limit exists and it is finite for every  $x\in[a,b]$. If $f=g$, we set $[f]:=[f,f]$ and it is called \textbf{quadratic variation of $f$ $($on $[a,b]$$)$}.
\end{Definition}

We notice that in Definition \ref{D:DeterministicQuadrVar} the quadratic variation $[f]$ is continuous on $[a,b]$, since $f$ is a continuous function. We conclude this subsection with an integration by parts formula for the deterministic forward and backward integrals.

\begin{Proposition}
\label{P:BVI}
Let $f\colon[a,b]\rightarrow\R$ be a c\`{a}dl\`{a}g function and $g\colon[a,b]\rightarrow\R$ be a bounded variation function. Then, the following \textbf{integration by parts formulae} hold:
\begin{equation}
\label{E:IbyP-}
\int_{[a,b]}  g(s) d^-f(s) \ = \ g(b^-) f(b)  - \int_{]a,b]} f(s) dg(s)
\end{equation}
and, if $f$ is continuous,
\begin{equation}
\label{E:IbyP}
\int_{[a,b]}  g(s) d^+f(s) \ = \ g(b) f(b)  - \int_{]a,b]} f(s) dg(s).
\end{equation}
\end{Proposition}
\textbf{Proof.}
Let us prove the integration by parts formula \eqref{E:IbyP} relative to the backward integral. By its definition, we are led to consider the following expression, for $\eps>0$,
\begin{align*}
\int_a^b g(s)\frac{f(s)-f(s-\eps)}{\eps}ds \ &= \ \frac{1}{\eps}\int_{b-\eps}^b g(s+\eps)f(s) ds - \frac{1}{\eps}\int_{a-\eps}^a g(s+\eps)f(s) ds \\
&\quad \ - \int_a^b \frac{g(s+\eps)-g(s)}{\eps}f(s) ds.
\end{align*}
Since $g(s)=g(b)$, for $s\geq b$, we get
\[
\frac{1}{\eps}\int_{b-\eps}^b g(s+\eps)f(s) ds \ = \ g(b)\frac{1}{\eps}\int_{b-\eps}^b f(s) ds \ \overset{\eps\rightarrow0^+}{\longrightarrow} \ g(b)f(b).
\]
Moreover, since $f(s)=0$, for $s<a$, we obtain
\[
\frac{1}{\eps}\int_{a-\eps}^a g(s+\eps)f(s) ds \ = \ 0.
\]
Finally, by Fubini's theorem we have
\begin{align*}
\int_a^b \frac{g(s+\eps)-g(s)}{\eps}f(s) ds \ &= \ \int_a^b \frac{1}{\eps}\bigg(\int_{]s,s+\eps]}dg(r)\bigg)f(s) ds \\
&= \ \int_{]a,b+\eps]} \frac{1}{\eps}\bigg(\int_{a\vee(r-\eps)}^{b\wedge r}f(s) ds\bigg)dg(r).
\end{align*}
Recalling that $g(s)=g(b)$, for $s\geq b$, and $f(s)=0$, for $s<a$, we obtain
\begin{align*}
\int_a^b \frac{g(s+\eps)-g(s)}{\eps}f(s) ds \ &= \ \int_{]a,b]} \frac{1}{\eps}\bigg(\int_{a\vee(r-\eps)}^{b\wedge r}f(s) ds\bigg)dg(r) \\
&= \ \int_{]a,b]} \frac{1}{\eps}\bigg(\int_{r-\eps}^r f(s) ds\bigg)dg(r).
\end{align*}
Since $f$ is continuous at $r\in]a,b]$, it follows that $\int_{r-\eps}^r f(s) ds/\eps\rightarrow f(r)$ as $\eps\rightarrow0^+$. Therefore, by Lebesgue's dominated convergence theorem, we find
\[
\int_{]a,b]} \frac{1}{\eps}\bigg(\int_{r-\eps}^r f(s) ds\bigg)dg(r) \ \overset{\eps\rightarrow0^+}{\longrightarrow} \ \int_{]a,b]} f(r)dg(r),
\]
which implies the validity of \eqref{E:IbyP}. The integration by parts formula relative to the forward integral is proved analogously.
\ep

\subsection{The spaces $\mathscr C([-T,0])$ and $\mathscr C([-T,0[)$}
\label{SubS:Preliminaries}

Let $C([-T,0])$ denote the set of real continuous functions on $[-T,0]$, endowed with supremum norm $\|\eta\|_\infty = \sup_{x\in[-T,0]}|\eta(x)|$, for any $\eta\in C([-T,0])$.

\begin{Remark} 
{\rm
We shall develop functional It\^o calculus via regularization firstly for time-independent functionals $\Uc\colon C([-T,0])\rightarrow\R$, since we aim at emphasizing that in our framework the time variable and the path play two distinct roles, as emphasized in the introduction. This, also, allows us to focus only on the definition of horizontal and vertical derivatives. Clearly, everything can be extended in an obvious way to the time-dependent case $\Uc\colon[0,T]\times C([-T,0])\rightarrow\R$, as we shall illustrate later.
\ep
}
\end{Remark}

Consider a map $\Uc\colon C([-T,0])\rightarrow\R$. Our aim is to derive a functional It\^o's formula for $\Uc$. To do this, we are led to define, in the spirit of \cite{dupire} and \cite{contfournie10}, the functional (i.e., horizontal and vertical) derivatives for $\Uc$. Since the definition of functional derivatives necessitates of discontinuous paths, in \cite{contfournie10} the idea is to consider functionals defined on the space of c\`adl\`ag trajectories $\D([-T,0])$. However, we can not, in general, extend in a unique way a functional $\Uc$ defined on $C([-T,0])$ to $\D([-T,0])$. Our idea, instead, is to consider an intermediate space between $C([-T,0])$ and $\D([-T,0])$, denoted $\mathscr C([-T,0])$, which is the space of bounded trajectories on $[-T,0]$, continuous on $[-T,0[$ and with possibly a jump at $0$. We endow $\mathscr C([-T,0])$ with a (inductive) topology such that $C([-T,0])$ is dense in $\mathscr C([-T,0])$ with respect to this topology. Therefore, if $\Uc$ is continuous with respect to the topology of $\mathscr C([-T,0])$, then it admits a unique continuous extension $u\colon \mathscr C([-T,0])\rightarrow\R$.

\begin{Definition}
\label{D:scrC}
We denote by $\mathscr C([-T,0])$ the set of bounded functions $\eta\colon[-T,0]\rightarrow\R$ such that $\eta$ is continuous on $[-T,0[$, equipped with the topology we now describe.\\
\textbf{Convergence.} We endow $\mathscr C([-T,0])$ with a topology inducing the following convergence: $(\eta_n)_n$ converges to $\eta$ in $\mathscr C([-T,0])$ as $n$ tends to infinity if:
\begin{enumerate}
\item[\textup{(i)}] $\|\eta_n\|_\infty \leq C$, for any $n\in\N$, for some positive constant $C$ independent of $n$$;$
\item[\textup{(ii)}] $\sup_{x\in K}|\eta_n(x)-\eta(x)|\rightarrow0$ as $n$ tends to infinity, for any compact set $K\subset[-T,0[$$;$
\item[\textup{(iii)}] $\eta_n(0)\rightarrow\eta(0)$ as $n$ tends to infinity.
\end{enumerate}
\textbf{Topology.} For each compact $K\subset [-T,0[$ define the seminorm $p_K$ on $\mathscr C([-T,0])$ by
\[
p_K(\eta) \ = \ \sup_{x\in K}|\eta(x)| + |\eta(0)|, \qquad \forall\,\eta\in\mathscr C([-T,0]).
\]
Let $M>0$ and $\mathscr C_M([-T,0])$ be the set of functions in $\mathscr C([-T,0])$ which are bounded by $M$. Still denote $p_K$ the restriction of $p_K$ to $\mathscr C_M([-T,0])$ and consider the topology on $\mathscr C_M([-T,0])$ induced by the collection of seminorms $(p_K)_K$. Then, we endow $\mathscr C([-T,0])$ with the smallest topology $($inductive topology$)$ turning all the inclusions $i_M\colon \mathscr C_M([-T,0])\rightarrow\mathscr C([-T,0])$ into continuous maps.
\end{Definition}

\begin{Remark}
\label{R:Density}
{\rm
(i) Notice that $C([-T,0])$ is dense in $\mathscr C([-T,0])$, when endowed with the topology of $\mathscr C([-T,0])$. As a matter of fact, let $\eta\in\mathscr C([-T,0])$ and define, for any $n\in\N\backslash\{0\}$,
\[
\varphi_n(x)=
\begin{cases}
\eta(x), \qquad &-T\leq x\leq-1/n, \\
n(\eta(0)-\eta(-1/n))x + \eta(0), &-1/n<x\leq0.
\end{cases}
\]
Then, we see that $\varphi_n\in C([-T,0])$ and $\varphi_n\rightarrow\eta$ in $\mathscr C([-T,0])$.\\
Now, for any $a\in\R$ define
\begin{align*}
C_a([-T,0]) \ &:= \ \{\eta\in C([-T,0])\colon\eta(0)=a\}, \\\mathscr C_a([-T,0]) \ &:= \ \{\eta\in\mathscr C([-T,0])\colon\eta(0)=a\}.
\end{align*}
Then, $C_a([-T,0])$ is dense in $\mathscr C_a([-T,0])$ with respect to the topology of $\mathscr C([-T,0])$. \\
(ii) We provide two examples of functionals $\Uc\colon C([-T,0])\rightarrow\R$, continuous with respect to the topology of $\mathscr C([-T,0])$, and necessarily with respect to the topology of $C([-T,0])$ (the proof is straightforward and not reported):
\begin{enumerate}
\item[(a)] $\Uc(\eta) = g(\eta(t_1),\ldots,\eta(t_n))$, for all $\eta\in C([-T,0])$, with $-T\leq t_1<\cdots<t_n\leq0$ and $g\colon\R^n\rightarrow\R$ continuous.
\item[(b)] $\Uc(\eta) = \int_{[-T,0]}\varphi(x)d^-\eta(x)$, for all $\eta\in C([-T,0])$, with $\varphi\colon[0,T]\rightarrow\R$ a c\`adl\`ag bounded variation function. Concerning this example, keep in mind that, using the integration by parts formula, $\Uc(\eta)$ admits the representation \eqref{E:IbyP-}.
\end{enumerate}
On the other hand, consider the functional $\Uc(\eta) = 
\sup_{x\in[-T,0]}\eta(x)$, for all $\eta\in C([-T,0])$.
It is obviously continuous, but it
 is not continuous with 
respect to the topology of $\mathscr C([-T,0])$. As a matter of fact, for any $n\in\N$ consider $\eta_n\in C([-T,0])$ given by
\[
\eta_n(x) \ = \
\begin{cases}
0, \qquad\qquad\qquad &-T\leq x\leq-\frac{T}{2^n}, \\
\frac{2^{n+1}}{T}x+2, &-\frac{T}{2^n}<x\leq-\frac{T}{2^{n+1}}, \\
-\frac{2^{n+1}}{T}x, &-\frac{T}{2^{n+1}}<x\leq0.
\end{cases}
\]
Then, $\Uc(\eta_n)=\sup_{x\in[-T,0]}\eta_n(x)=1$, for any $n$. However, $\eta_n$ converges to the zero function in $\mathscr C([-T,0])$, as $n$ tends to infinity. This example  will play an important role in Section \ref{S:StrongViscositySolutions} to justify a weaker notion of solution to the path-dependent nonlinear Kolmogorov equation.
\ep
}
\end{Remark}

To define the functional derivatives, we shall need to separate the
 ``past'' from the ``present'' of $\eta\in\mathscr C([-T,0])$. Indeed, roughly speaking, the horizontal derivative calls in the past values of $\eta$, namely $\{\eta(x)\colon x\in[-T,0[\}$, while the vertical derivative calls in the present value of $\eta$, namely $\eta(0)$. To this end, it is useful to introduce the space $\mathscr C([-T,0[)$.

\begin{Definition}
\label{D:scrC2}
We denote by $\mathscr C([-T,0[)$ the set of bounded continuous functions $\gamma\colon[-T,0[\rightarrow\R$, equipped with the topology we now describe.\\
\textbf{Convergence.} We endow $\mathscr C([-T,0[)$ with a topology inducing the following convergence: $(\gamma_n)_n$ converges to $\gamma$ in $\mathscr C([-T,0[)$ as $n$ tends to infinity if:
\begin{enumerate}
\item[\textup{(i)}] $\sup_{x\in[-T,0[}|\gamma_n(x)| \leq C$, for any $n\in\N$, for some positive constant $C$ independent of $n$$;$
\item[\textup{(ii)}] $\sup_{x\in K}|\gamma_n(x)-\gamma(x)|\rightarrow0$ as $n$ tends to infinity, for any compact set $K\subset[-T,0[$.
\end{enumerate}
\textbf{Topology.} For each compact $K\subset [-T,0[$ define the seminorm $q_K$ on $\mathscr C([-T,0[)$ by
\[
q_K(\gamma) \ = \ \sup_{x\in K}|\gamma(x)|, \qquad \forall\,\gamma\in\mathscr C([-T,0[).
\]
Let $M>0$ and $\mathscr C_M([-T,0[)$ be the set of functions in $\mathscr C([-T,0[)$ which are bounded by $M$. Still denote $q_K$ the restriction of $q_K$ to $\mathscr C_M([-T,0[)$ and consider the topology on $\mathscr C_M([-T,0[)$ induced by the collection of seminorms $(q_K)_K$. Then, we endow $\mathscr C([-T,0[)$ with the smallest topology $($inductive topology$)$ turning all the inclusions $i_M\colon \mathscr C_M([-T,0[)\rightarrow\mathscr C([-T,0[)$ into continuous maps.
\end{Definition}

\begin{Remark}
\label{R:Isomorphism}
{\rm
(i) Notice that $\mathscr C([-T,0])$ is isomorphic to $\mathscr C([-T,0[)\times\R$. As a matter of fact, it is enough to consider the map
\begin{align*}
J \colon \mathscr C([-T,0]) &\rightarrow \mathscr C([-T,0[)\times\R \\
\eta &\mapsto (\eta_{|[-T,0[},\eta(0)).
\end{align*}
Observe that $J^{-1}\colon \mathscr C([-T,0[)\times\R \rightarrow \mathscr C([-T,0])$ is given by $J^{-1}(\gamma,a) = \gamma1_{[-T,0[} + a1_{\{0\}}$.\\
(ii) $\mathscr C([-T,0])$ is a space which contains $C([-T,0])$ as a subset and it has the property of 
separating ``past'' from ``present''. Another space having the same property is $L^2([-T,0]; d \mu)$ where $\mu$ is the  sum of
the Dirac measure at zero and Lebesgue measure. Similarly as for item (i), that space is isomorphic to $L^2([-T,0]) \times \R $, which is a very popular space appearing in the analysis of functional dependent (as delay) equations, starting from \cite{Cho}.
\ep
}
\end{Remark}

For every $u\colon \mathscr C([-T,0])\rightarrow\R$, we can now exploit the space $\mathscr C([-T,0[)$ to define a map $\tilde u\colon \mathscr C([-T,0[)\times\R\rightarrow\R$ where ``past'' and ``present'' are separated.

\begin{Definition}
\label{D:tildeu}
Let $u\colon \mathscr C([-T,0])\rightarrow\R$ and define $\tilde u\colon \mathscr C([-T,0[)\times\R\rightarrow\R$ as
\begin{equation}
\label{E:tildeu}
\tilde u(\gamma,a) \ := \ u(\gamma 1_{[-T,0[} + a1_{\{0\}}), \qquad \forall\,(\gamma,a)\in \mathscr C([-T,0[)\times\R.
\end{equation}
In particular, we have $u(\eta) = \tilde u(\eta_{|[-T,0[},\eta(0))$, for all $\eta\in \mathscr C([-T,0])$.
\end{Definition}

We conclude this subsection with a characterization of the dual spaces of $\mathscr C([-T,0])$ and $\mathscr C([-T,0[)$, which has an independent interest. Firstly, we need to introduce the set $\Mc([-T,0])$ of finite signed Borel measures on $[-T,0]$. We also denote $\Mc_0([-T,0])\subset\Mc([-T,0])$ the set of measures $\mu$ such that $\mu(\{0\})=0$.

\begin{Proposition}
\label{P:DualSpace}
Let $\Lambda\in \mathscr C([-T,0])^*$, the dual space of $\mathscr C([-T,0])$. Then, there exists a unique $\mu\in\Mc([-T,0])$ such that
\[
\Lambda\eta \ = \ \int_{[-T,0]} \eta(x) \mu(dx), \qquad \forall\,\eta\in \mathscr C([-T,0]).
\]
\end{Proposition}
\textbf{Proof.}
Let $\Lambda\in \mathscr C([-T,0])^*$ and define
\[
\tilde\Lambda\varphi \ := \ \Lambda\varphi, \qquad \forall\,\varphi\in C([-T,0]).
\]
Notice that $\tilde\Lambda\colon C([-T,0])\rightarrow\R$ is a continuous functional on the Banach space $C([-T,0])$ endowed with the supremum norm $\|\cdot\|_\infty$. Therefore $\tilde\Lambda\in C([-T,0])^*$ and it follows from Riesz representation theorem (see, e.g., Theorem 6.19 in \cite{rudin}) that there exists a unique $\mu\in\Mc([-T,0])$ such that
\[
\tilde\Lambda\varphi \ = \ \int_{[-T,0]} \varphi(x)\mu(dx), \qquad \forall\,\varphi\in C([-T,0]).
\]
Obviously $\tilde\Lambda$ is also continuous with respect to the topology of $\mathscr C([-T,0])$. Since $C([-T,0])$ is dense in $\mathscr C([-T,0])$ with respect to the topology of $\mathscr C([-T,0])$, we deduce that there exists a unique continuous extension of $\tilde\Lambda$ to $\mathscr C([-T,0])$, which is clearly given by
\[
\Lambda\eta \ = \ \int_{[-T,0]} \eta(x)\mu(dx), \qquad \forall\,\eta\in \mathscr C([-T,0]).
\]
\ep

\begin{Proposition}
\label{P:DualSpace[}
Let $\Lambda\in \mathscr C([-T,0[)^*$, the dual space of $\mathscr C([-T,0[)$. Then, there exists a unique $\mu\in\Mc_0([-T,0])$ such that
\[
\Lambda\gamma \ = \ \int_{[-T,0[} \gamma(x) \mu(dx), \qquad \forall\,\gamma\in \mathscr C([-T,0[).
\]
\end{Proposition}
\textbf{Proof.}
Let $\Lambda\in \mathscr C([-T,0[)^*$ and define
\begin{equation}
\label{E:DualSpace[Proof}
\tilde\Lambda\eta \ := \ \Lambda(\eta_{|[-T,0[}), \qquad \forall\,\eta\in\mathscr C([-T,0]).
\end{equation}
Notice that $\tilde\Lambda\colon\mathscr C([-T,0])\rightarrow\R$ is a continuous functional on $\mathscr C([-T,0])$. It follows from Proposition \ref {P:DualSpace} that there exists a unique $\mu\in\Mc([-T,0])$ such that
\begin{equation}
\label{E:DualSpace[Proof2}
\tilde\Lambda\eta \ = \ \int_{[-T,0]} \eta(x)\mu(dx) \ = \ \int_{[-T,0[} \eta(x)\mu(dx) + \eta(0)\mu(\{0\}), \qquad \forall\,\eta\in\mathscr C([-T,0]).
\end{equation}
Let $\eta_1,\eta_2\in\mathscr C([-T,0])$ be such that $\eta_1 1_{[-T,0[}=\eta_2 1_{[-T,0[}$. Then, we see from \eqref{E:DualSpace[Proof} that $\tilde\Lambda\eta_1=\tilde\Lambda\eta_2$, which in turn implies from \eqref{E:DualSpace[Proof2} that $\mu(\{0\})=0$. In conclusion, $\mu\in\Mc_0([-T,0])$ and $\Lambda$ is given by
\[
\Lambda\gamma \ = \ \int_{[-T,0[} \gamma(x)\mu(dx), \qquad \forall\,\gamma\in \mathscr C([-T,0[).
\]
\ep

\subsection{Functional derivatives and functional It\^o's formula}
\label{SubS:PathwiseDerivatives}

In the present section we shall prove one of the main result of this section, namely the functional It\^o's formula for $\Uc\colon C([-T,0])\rightarrow\R$ and, more generally, for $\Uc\colon[0,T]\times C([-T,0])\rightarrow\R$. We begin introducing the functional derivatives in the spirit of Dupire \cite{dupire}, firstly for a functional $u\colon\mathscr C([-T,0])\rightarrow\R$, and then for $\Uc\colon C([-T,0])\rightarrow\R$.

\begin{Definition}
\label{D:DupireDerivatives}
\quad\\ Consider $u\colon\mathscr C([-T,0])\rightarrow\R$ and $\eta\in\mathscr C([-T,0])$. \\
\textup{(i)} We say that $u$ admits \textbf{horizontal derivative} at $\eta$ if the following limit exists and it is finite:
\begin{equation}
\label{E:DH}
D^H u(\eta) \ := \ \lim_{\eps\rightarrow0^+} \frac{u(\eta(\cdot)1_{[-T,0[}+\eta(0)1_{\{0\}}) - u(\eta(\cdot-\eps)1_{[-T,0[}+\eta(0)1_{\{0\}})}{\eps}.
\end{equation}
\textup{(i)'} Let $\tilde u$ be as in \eqref{E:tildeu}, then we say that $\tilde u$ admits \textbf{horizontal derivative} at $(\gamma,a)\in\mathscr C([-T,0[)\times\R$ if the following limit exists and it is finite:
\begin{equation}
\label{E:DHtildeu}
D^H\tilde u(\gamma,a) \ := \ \lim_{\eps\rightarrow0^+} \frac{\tilde u(\gamma(\cdot),a) - \tilde u(\gamma(\cdot-\eps),a)}{\eps}.
\end{equation}
Notice that if $D^H u(\eta)$ exists then $D^H\tilde u(\eta_{|[-T,0[},\eta(0))$ exists and they are equal; viceversa, if $D^H\tilde u(\gamma,a)$ exists then $D^H u(\gamma 1_{[-T,0[} + a 1_{\{0\}})$ exists and they are equal. \\
\textup{(ii)} We say that $u$ admits \textbf{first-order vertical
 derivative} at $\eta$ if the first-order partial derivative at $(\eta_{|[-T,0[},\eta(0))$ of $\tilde u$ with respect to its second argument, denoted by $\partial_a\tilde u(\eta_{|[-T,0[},\eta(0))$, exists and we set
\[
D^V u(\eta) \ := \ \partial_a \tilde u(\eta_{|[-T,0[},\eta(0)).
\]
\textup{(iii)} We say that $u$ admits \textbf{second-order vertical derivative} at $\eta$ if the second-order partial derivative at $(\eta _{|[-T,0[},\eta(0))$ of $\tilde u$ with respect to its second argument, denoted by $\partial_{aa}^2\tilde u(\eta_{|[-T,0[},\eta(0))$, exists and we set
\[
D^{VV} u(\eta) \ := \ \partial_{aa}^2 \tilde u(\eta_{|[-T,0[},\eta(0)).
\]
\end{Definition}

\begin{Definition}
\label{D:C12}
We say that $u\colon \mathscr C([-T,0])\rightarrow\R$ is of class $\mathscr C^{1,2}(\textup{past}\times\textup{present})$ if:
\begin{enumerate}
\item[\textup{(i)}] $u$ is continuous$;$
\item[\textup{(ii)}] $D^H u$ exists everywhere on $\mathscr C([-T,0])$ and for every $\gamma\in \mathscr C([-T,0[)$ the map
\[
(\eps,a)\longmapsto D^H\tilde u(\gamma(\cdot-\eps),a), \qquad (\eps,a)\in[0,\infty[\times\R
\]
is continuous on $[0,\infty[\times\R$$;$
\item[\textup{(iii)}] $D^V u$ and $D^{VV}u$ exist everywhere on 
$\mathscr C([-T,0])$ and are continuous.
\end{enumerate}
\end{Definition}

\begin{Remark}
{\rm
Notice that in Definition \ref{D:C12} we still obtain the same class of functions $\mathscr C^{1,2}(\textup{past}\times\textup{present})$ if we substitute point (ii) with:
\begin{enumerate}
\item[\textup{(ii')}] $D^H u$ exists everywhere on $\mathscr C([-T,0])$ and for every $\gamma\in \mathscr C([-T,0[)$ there exists $\delta(\gamma)>0$ such that the map
\begin{equation}
\label{E:map}
(\eps,a)\longmapsto D^H\tilde u(\gamma(\cdot-\eps),a), \qquad (\eps,a)\in[0,\infty[\times\R
\end{equation}
is continuous on $[0,\delta(\gamma))\times\R$.
\end{enumerate}
In particular, if (ii') holds then we can always take $\delta(\gamma) = \infty$ for any $\gamma\in \mathscr C([-T,0[)$, which implies (ii). To prove this last statement, let us proceed by contradiction assuming that
\[
\delta^*(\gamma) \ = \ \sup\big\{\delta(\gamma)>0\colon\text{the map \eqref{E:map} is continuous on }[0,\delta(\gamma)[\times\R\big\} \ < \ \infty.
\]
Notice that $\delta^*(\gamma)$ is in fact a \emph{max}, therefore the map \eqref{E:map} is continuous on $[0,\delta^*(\gamma)[\times\R$. Now, define $\bar{\gamma}(\cdot) := \gamma(\cdot-\delta^*(\gamma))$. Then, by condition (ii') there exists $\delta(\bar{\gamma})>0$ such that the map
\[
(\eps,a)\longmapsto D^H\tilde u(\bar{\gamma}(\cdot-\eps),a) = D^H\tilde u(\gamma(\cdot-\eps-\delta^*(\gamma)),a)
\]
is continuous on $[0,\delta(\bar{\gamma})[\times\R$. This shows that the map \eqref{E:map} is continuous on $[0,\delta^*(\gamma)+\delta(\bar{\gamma})[\times\R$, a contradiction with the definition of $\delta^*(\gamma)$.
\ep
}
\end{Remark}

\noindent We can now provide the definition of functional derivatives for a map $\Uc\colon$ $C([-T,0])$ $\rightarrow$ $\R$.

\begin{Definition}
\label{D:DupireDerivativesUc}
Let $\Uc\colon C([-T,0])\rightarrow\R$ and $\eta\in C([-T,0])$. Suppose that there exists a unique extension $u\colon\mathscr C([-T,0])\rightarrow\R$ of $\Uc$ $($e.g., if $\Uc$ is continuous with respect to the topology of $\mathscr C([-T,0])$$)$. Then we define:\\
\textup{(i)} The \textbf{horizontal derivative} of $\Uc$ at $\eta$ as:
\[
D^H \Uc(\eta) \ := \ D^H u(\eta).
\]
\textup{(ii)} The \textbf{first-order vertical derivative} of $\Uc$ at $\eta$ as:
\[
D^V \Uc(\eta) \ := \ D^V u(\eta).
\]
\textup{(iii)} The \textbf{second-order vertical derivative} of $\Uc$ at $\eta$ as:
\[
D^{VV} \Uc(\eta) \ := \ D^{VV} u(\eta).
\]
\end{Definition}

\begin{Definition}
\label{D:C12Uc}
We say that $\Uc\colon C([-T,0])\rightarrow\R$ is $C^{1,2}(\textup{past}\times\textup{present})$ if $\Uc$ admits a $($necessarily unique$)$ extension $u\colon \mathscr C([-T,0])\rightarrow\R$ of class $\mathscr C^{1,2}(\textup{past}\times\textup{present})$.
\end{Definition}

\begin{Theorem}
\label{T:Ito}
Let $\Uc\colon C([-T,0])\rightarrow\R$ be of class $C^{1,2}(\textup{past}\times\textup{present})$ and $X=(X_t)_{t\in[0,T]}$ be a real 
continuous finite quadratic variation process. Then, the following \textbf{functional It\^o's formula} holds, $\P$-a.s.,
\begin{equation}
\label{E:Ito}
\Uc(\X_t) \ = \ \Uc(\X_0) + \int_0^t D^H \Uc(\X_s)ds + \int_0^t D^V \Uc(\X_s) d^- X_s + \frac{1}{2}\int_0^t D^{VV}\Uc(\X_s)d[X]_s,
\end{equation}
for all $0 \leq t \leq T$.
\end{Theorem}
\textbf{Proof.}
Fix $t\in[0,T]$ and consider the quantity
\[
I_0(\eps,t) \ = \ \int_0^t \frac{\Uc(\X_{s+\eps}) - \Uc(\X_s)}{\eps} ds \ = \ \frac{1}{\eps} \int_t^{t+\eps} \Uc(\X_s) ds - \frac{1}{\eps} \int_0^\eps \Uc(\X_s) ds, \qquad \eps>0.
\]
Since $(\Uc(\X_s))_{s\geq0}$ is continuous, $I_0(\eps,t)$ converges ucp to $\Uc(\X_t) - \Uc(\X_0)$, i.e., $\sup_{0 \leq t \leq T}|I_0(\eps,t)-(\Uc(\X_t) - \Uc(\X_0))|$ converges to zero in probability when $\eps\rightarrow0^+$. On the other hand, we can write $I_0(\eps,t)$ in terms of the function $\tilde u$, defined in \eqref{E:tildeu}, as follows
\[
I_0(\eps,t) \ = \ \int_0^t \frac{\tilde u(\X_{s+\eps|[-T,0[},X_{s+\eps}) - \tilde u(\X_{s|[-T,0[},X_s)}{\eps} ds.
\]
Now we split $I_0(\eps,t)$ into two terms:
\begin{align}
I_1(\eps,t) \ &= \ \int_0^t \frac{\tilde u(\X_{s+\eps|[-T,0[},X_{s+\eps}) - \tilde u(\X_{s|[-T,0[},X_{s+\eps})}{\eps} ds, \label{E:I1} \\
I_2(\eps,t) \ &= \ \int_0^t \frac{\tilde u(\X_{s|[-T,0[},X_{s+\eps}) - \tilde u(\X_{s|[-T,0[},X_s)}{\eps} ds. \label{E:I2}
\end{align}
We begin proving that
\begin{equation}
\label{E:I1ucp}
I_1(\eps,t) \underset{\eps\rightarrow0^+}{\overset{\text{ucp}}{\longrightarrow}} \int_0^t D^H \Uc(\X_s) ds.
\end{equation}
Firstly, fix $\gamma\in \mathscr C([-T,0[)$ and define
\[
\phi(\eps,a) \ := \ \tilde u(\gamma(\cdot-\eps),a), \qquad (\eps,a)\in[0,\infty[\times\R.
\]
Then, denoting by $\partial_\eps^+ \phi$ the right partial derivative of $\phi$ with respect to $\eps$ and using formula \eqref{E:DHtildeu}, we find
\begin{align*}
\partial_\eps^+ \phi(\eps,a) \ &= \ \lim_{r\rightarrow0^+} \frac{\phi(\eps+r,a) - \phi(\eps,a)}{r} \\
&= \ -\lim_{r\rightarrow0^+} \frac{\tilde u(\gamma(\cdot-\eps),a) - \tilde u(\gamma(\cdot-\eps-r),a)}{r} \\
&= \ -D^H \tilde u(\gamma(\cdot-\eps),a), \qquad \forall\,(\eps,a)\in[0,\infty[\times\R.
\end{align*}
Since $u\in\mathscr C^{1,2}(\textup{past}\times\textup{present})$, we see from Definition \ref{D:C12}(ii), that $\partial_\eps^+ \phi$ is continuous on $[0,\infty[\times\R$. It follows from a standard differential calculus' result (see for example Corollary 1.2, Chapter 2, in \cite{pazy83}) that $\phi$ is continuously differentiable on $[0,\infty[\times\R$ with respect to its first argument. Then, for every $(\eps,a)\in[0,\infty[\times\R$, from the fundamental theorem of calculus, we have
\[
\phi(\eps,a) - \phi(0,a) \ = \ \int_0^\eps \partial_\eps\phi (r,a) dr,
\]
which in terms of $\tilde u$ reads
\begin{equation}
\label{E:tildeuFundThmCalc}
\tilde u(\gamma(\cdot),a) - \tilde u(\gamma(\cdot-\eps),a) \ = \ \int_0^\eps D^H \tilde u(\gamma(\cdot-r),a) dr.
\end{equation}
Now, we rewrite, by means of a shift in time, the term $I_1(\eps,t)$ in \eqref{E:I1} as follows:
\begin{align}
\label{E:I1bis}
I_1(\eps,t) \ &= \ \int_0^t \frac{\tilde u(\X_{s|[-T,0[},X_s) - \tilde u(\X_{s-\eps|[-T,0[},X_s)}{\eps} ds \notag \\
&\quad + \int_t^{t+\eps} \frac{\tilde u(\X_{s|[-T,0[},X_s) - \tilde u(\X_{s-\eps|[-T,0[},X_s)}{\eps} ds \notag \\
&\quad - \int_0^{\eps} \frac{\tilde u(\X_{s|[-T,0[},X_s) - \tilde u(\X_{s-\eps|[-T,0[},X_s)}{\eps} ds.
\end{align}
Plugging \eqref{E:tildeuFundThmCalc} into \eqref{E:I1bis}, setting $\gamma =  \X_s, a = X_s$, we obtain
\begin{align}
\label{E:I1R1}
I_1(\eps,t) \ &= \ \int_0^t \frac{1}{\eps} \bigg(\int_0^\eps D^H \tilde u(\X_{s-r|[-T,0[},X_s)dr\bigg)ds \notag \\
&\quad + \int_t^{t+\eps} \frac{1}{\eps} \bigg(\int_0^\eps D^H \tilde u(\X_{s-r|[-T,0[},X_s)dr\bigg)ds \notag \\
&\quad - \int_0^\eps \frac{1}{\eps} \bigg(\int_0^\eps D^H \tilde u(\X_{s-r|[-T,0[},X_s)dr\bigg)ds.
\end{align}
Observe that
\[
\int_0^t \frac{1}{\eps} \bigg(\int_0^\eps D^H \tilde u(\X_{s-r|[-T,0[},X_s)dr\bigg)ds \underset{\eps\rightarrow0^+}{\overset{\text{ucp}}{\longrightarrow}} \int_0^t D^H u(\X_s) ds.
\]
Similarly, we see that the other two terms in \eqref{E:I1R1} converge ucp to zero. As a consequence, we get \eqref{E:I1ucp}.

Regarding $I_2(\eps,t)$ in \eqref{E:I2}, it can be written, by means of the following standard Taylor's expansion for a function $f\in C^2(\R)$:
\begin{align*}
f(b) \ &= \ f(a)+f'(a)(b-a)+\frac{1}{2}f''(a)(b-a)^2 \\
&\quad + \int_0^1(1-\alpha)\big(f''(a+\alpha(b-a))-f''(a)\big)(b-a)^2d\alpha,
\end{align*}
as the sum of the following three terms:
\begin{align*}
I_{21}(\eps,t) \ &= \ \int_0^t \partial_a \tilde u(\X_{s|[-T,0[},X_s) \frac{X_{s+\eps}-X_s}{\eps} ds \\
I_{22}(\eps,t) \ &= \ \frac{1}{2}\int_0^t \partial_{aa}^2 \tilde u(\X_{s|[-T,0[},X_s) \frac{(X_{s+\eps}-X_s)^2}{\eps} ds \\
I_{23}(\eps,t) \ &= \ \int_0^t \bigg(\int_0^1(1-\alpha) \big( \partial_{aa}^2 \tilde u(\X_{s|[-T,0[},X_s + \alpha(X_{s+\eps}-X_s)) \\
&\quad\;\, - \partial_{aa}^2 \tilde u(\X_{s|[-T,0[},X_s) \big) \frac{(X_{s+\eps}-X_s)^2}{\eps} d\alpha \bigg) ds.
\end{align*}
By similar arguments as in Proposition 1.2 of \cite{russovallois95}, we have 
\[
I_{22}(\eps,t) \underset{\eps\rightarrow0^+}{\overset{\text{ucp}}{\longrightarrow}} \frac{1}{2}\int_0^t\partial_{aa}^2 \tilde u(\X_{s|[-T,0[},X_s) d[X]_s = \frac{1}{2}\int_0^t D^{VV} u(\X_s) d[X]_s.
\]
Regarding $I_{23}(\eps,t)$, for every $\omega\in\Omega$, define $\psi_\omega\colon[0,T]\times[0,1]\times[0,1]\rightarrow\R$ as
\[
\psi_\omega(s,\alpha,\eps) \ := \ (1-\alpha) \partial_{aa}^2 \tilde u\big(\X_{s|[-T,0[}(\omega),X_s(\omega) + \alpha(X_{s+\eps}(\omega)-X_s(\omega))\big),
\]
for all $(s,\alpha,\eps)\in[0,T]\times[0,1]\times[0,1]$. Notice that $\psi_\omega$ is uniformly continuous. Denote $\rho_{\psi_\omega}$ its continuity modulus, then
\[
\sup_{t\in[0,T]} |I_{23}(\eps,t)| \ \leq \ \int_0^T\rho_{\psi_\omega}(\eps) \frac{(X_{s+\eps}-X_s)^2}{\eps}ds.
\]
Since $X$ has finite quadratic variation, we deduce that $I_{23}(\eps,t)\rightarrow0$ ucp as $\eps\rightarrow0^+$. Finally, because of $I_0(\eps,t)$, $I_1(\eps,t)$, $I_{22}(\eps,t)$, and $I_{23}(\eps,t)$ converge ucp, it follows that the forward integral exists:
\[
I_{21}(\eps,t) \underset{\eps\rightarrow0^+}{\overset{\text{ucp}}{\longrightarrow}} \int_0^t \partial_a \tilde u(\X_{s|[-T,0[},X_s) d^- X_s = \int_0^t D^V u(\X_s) d^- X_s,
\]
from which the thesis follows.
\ep

\begin{Remark}
{\rm
Notice that, under the hypotheses of Theorem \ref{T:Ito}, the forward integral $\int_0^t D^V\Uc(\X_s)d^-X_s$ exists as a ucp limit, which is generally not required.
\ep
}
\end{Remark}

\begin{Remark}
{\rm
\emph{The definition of horizontal derivative}. Notice that our definition of horizontal derivative differs from that introduced in \cite{dupire}, since it is based on a limit on the left, while the definition proposed in \cite{dupire} would conduct to the following formula:
\begin{equation}
\label{E:DHright}
D^{H,+} u(\eta) \ := \ \lim_{\eps\rightarrow0^+} \frac{\tilde u(\eta(\cdot+\eps) 1_{[-T,0[},\eta(0)) - \tilde u(\eta(\cdot) 1_{[-T,0[},\eta(0))}{\eps}.
\end{equation}
To give an insight into the difference between \eqref{E:DH} and \eqref{E:DHright}, let us consider a real continuous finite quadratic variation process $X$ with associated window process $\X$. Then, in the definition \eqref{E:DHright} of $D^{H,+}u(\X_t)$ we consider the increment $\tilde u(\X_{t|[-T,0[}(\cdot+\eps),X_t) - \tilde u(\X_{t|[-T,0[},X_t)$, comparing the present value of $u(\X_t)=\tilde u(\X_{t|[-T,0[},X_t)$ with an hypothetical future value $\tilde u(\X_{t|[-T,0[}(\cdot+\eps),X_t)$, obtained assuming a constant time evolution for $X$. On the other hand, in our definition \eqref{E:DH} we consider the increment $\tilde u(\X_{t|[-T,0[},X_t) - \tilde u(\X_{t-\eps|[-T,0[},X_t)$, where only the present and past values of $X$ are taken into account, and where we also extend in a constant way the trajectory of $X$ before time $0$. In particular, unlike \eqref{E:DHright}, since we do not call in the future in our formula \eqref{E:DH}, we do not have to specify a future time evolution for $X$, but only a past evolution before time $0$. This difference between \eqref{E:DH} and \eqref{E:DHright} is crucial for the proof of the functional It\^o's formula. In particular, the adoption of \eqref{E:DHright} as definition for the horizontal derivative would require an additional regularity condition on $u$ in order to prove an It\^o's formula for the process $t\mapsto u(\X_t)$. Indeed, as it can be seen from the proof of Theorem \ref{T:Ito}, to prove It\^o's formula we are led to consider the following term:
\[
I_1(\eps,t) \ = \ \int_0^t \frac{\tilde u(\X_{s+\eps|[-T,0[},X_{s+\eps}) - \tilde u(\X_{s|[-T,0[},X_{s+\eps})}{\eps} ds.
\]
When adopting definition \eqref{E:DHright} it is convenient to write $I_1(\eps,t)$ as the sum of two integrals:
\begin{align*}
I_{11}(\eps,t) \ &= \ \int_0^t \frac{\tilde u(\X_{s+\eps|[-T,0[},X_{s+\eps}) - \tilde u(\X_{s|[-T,0[}(\cdot+\eps),X_{s+\eps})}{\eps} ds, \\
I_{12}(\eps,t) \ &= \ \int_0^t \frac{\tilde u(\X_{s|[-T,0[}(\cdot+\eps),X_{s+\eps}) - \tilde u(\X_{s|[-T,0[},X_{s+\eps})}{\eps} ds.
\end{align*}
It can be shown quite easily that, under suitable regularity conditions on $u$ (more precisely, if $u$ is continuous, $D^{H,+}u$ exists everywhere on $\mathscr C([-T,0])$, and for every $\gamma\in \mathscr C([-T,0[)$ the map $(\eps,a)\longmapsto D^{H,+}\tilde u(\gamma(\cdot+\eps),a)$ is continuous on $[0,\infty)\times\R$), we have
\[
I_{12}(\eps,t) \underset{\eps\rightarrow0^+}{\overset{\text{ucp}}{\longrightarrow}} \int_0^t D^{H,+} u(\X_s) ds.
\]
To conclude the proof of It\^o's formula along the same lines as in Theorem \ref{T:Ito}, we should prove
\begin{equation}
\label{E:I11->0}
I_{11}(\eps,t) \underset{\eps\rightarrow0^+}{\overset{\text{ucp}}{\longrightarrow}} 0.
\end{equation}
In order to guarantee \eqref{E:I11->0}, we need to impose some additional regularity condition on $\tilde u$, and hence on $u$. As an example, \eqref{E:I11->0} is satisfied if we assume the following condition on $\tilde u$: there exists a constant $C>0$ such that, for every $\eps>0$,
\[
|\tilde u(\gamma_1,a) - \tilde u(\gamma_2,a)| \ \leq \ C \eps \sup_{x\in[-\eps,0[}|\gamma_1(x)-\gamma_2(x)|,
\]
for all $\gamma_1,\gamma_2\in \mathscr C([-T,0[)$ and $a\in\R$, with $\gamma_1(x) = \gamma_2(x)$ for any $x\in[-T,-\eps]$. This last condition is verified if, for example, $\tilde u$ is uniformly Lipschitz continuous with respect to the $L^1([-T,0])$-norm on $\mathscr C([-T,0[)$, namely: there exists a constant $C>0$ such that
\[
|\tilde u(\gamma_1,a) - \tilde u(\gamma_2,a)| \ \leq \ C \int_{[-T,0[} |\gamma_1(x) - \gamma_2(x)| dx,
\]
for all $\gamma_1,\gamma_2\in \mathscr C([-T,0[)$ and $a\in\R$.
\ep
}
\end{Remark}

We conclude this subsection providing the functional It\^o's formula for a map $\Uc\colon[0,T]\times C([-T,0])\rightarrow\R$ depending also on the time variable. Firstly, we notice that for a map $\Uc\colon[0,T]\times C([-T,0])\rightarrow\R$ (resp. $u\colon[0,T]\times\mathscr C([-T,0])\rightarrow\R$) the functional derivatives $D^H\Uc$, $D^V\Uc$, and $D^{VV}\Uc$ (resp. $D^Hu$, $D^Vu$, and $D^{VV}u$) are defined in an obvious way as in Definition \ref{D:DupireDerivativesUc} (resp. Definition \ref{D:DupireDerivatives}). Moreover, given $u\colon[0,T]\times\mathscr C([-T,0])\rightarrow\R$ we can define, as in Definition \ref{D:tildeu}, a map $\tilde u\colon[0,T]\times\mathscr C([-T,0[)\times\R\rightarrow\R$. Then, we can give the following definitions.

\begin{Definition}
\label{D:C12Time}
Let $I$ be $[0,T[$ or $[0,T]$. We say that $u\colon I\times\mathscr C([-T,0])\rightarrow\R$ is of class $\mathscr C^{1,2}((I\times\textup{past})\times\textup{present})$ if the properties below hold.
\begin{enumerate}
\item[\textup{(i)}] $u$ is continuous$;$
\item[\textup{(ii)}] $\partial_tu$ exists everywhere on $I\times\mathscr C([-T,0])$ and is continuous$;$
\item[\textup{(iii)}] $D^H u$ exists everywhere on $I\times\mathscr C([-T,0])$ and for every $\gamma\in \mathscr C([-T,0[)$ the map
\[
(t,\eps,a)\longmapsto D^H\tilde u(t,\gamma(\cdot-\eps),a), \qquad (t,\eps,a)\in I\times[0,\infty[\times\R
\]
is continuous on $I\times[0,\infty[\times\R$$;$
\item[\textup{(iv)}] $D^V u$ and $D^{VV}u$ exist everywhere on $I\times\mathscr C([-T,0])$ and are continuous.
\end{enumerate}
\end{Definition}

\begin{Definition}
\label{D:C12UcTime}
Let $I$ be $[0,T[$ or $[0,T]$. We say that $\Uc\colon I\times C([-T,0])\rightarrow\R$ is $C^{1,2}((I\times\textup{past})\times\textup{present}))$ if $\Uc$ admits a $($necessarily unique$)$ extension $u\colon I\times\mathscr C([-T,0])\rightarrow\R$ of class $\mathscr C^{1,2}((I\times\textup{past})\times\textup{present})$.
\end{Definition}

We can now state the functional It\^o's formula, whose proof is not reported, since it can be done along the same lines as Theorem \ref{T:Ito}.

\begin{Theorem}
\label{T:ItoTime}
Let $\Uc\colon[0,T]\times C([-T,0])\rightarrow\R$ be of class $C^{1,2}(([0,T]\times\textup{past})\times\textup{present})$ and $X=(X_t)_{t\in[0,T]}$ be a real continuous finite quadratic variation process. Then, the following \textbf{functional It\^o's formula} holds, $\P$-a.s.,
\begin{align}
\label{E:ItoTime}
\Uc(t,\X_t) \ &= \ \Uc(0,\X_0) + \int_0^t \big(\partial_t\Uc(s,\X_s) + D^H \Uc(s,\X_s)\big)ds + \int_0^t D^V \Uc(s,\X_s) d^- X_s \notag \\
&\quad \ + \frac{1}{2}\int_0^t D^{VV}\Uc(s,\X_s)d[X]_s,
\end{align}
for all $0 \leq t \leq T$.
\end{Theorem}

\begin{Remark}
{\rm
Notice that, as a particular case, choosing $\Uc(t,\eta)=F(t,\eta(0))$, for any $(t,\eta)\in[0,T]\times C([-T,0])$, with $F\in C^{1,2}([0,T]\times\R)$, we retrieve the classical It\^o's formula 
for finite quadratic variation processes, i.e. \eqref{E:ITOFQV}. More precisely, in this case $\Uc$ admits as unique continuous extension the map $u\colon[0,T]\times\mathscr C([-T,0])\rightarrow\R$ given by $u(t,\eta)=F(t,\eta(0))$, for all $(t,\eta)\in[0,T]\times\mathscr C([-T,0])$. Moreover, we see that $D^H\Uc\equiv0$, while $D^V\Uc=\partial_x F$ and $D^{VV}\Uc=\partial_{xx}^2 F$, where $\partial_x F$ (resp. $\partial_{xx}^2F$) denotes the first-order (resp. second-order) partial derivative of $F$ with respect to its second argument.
\ep
}
\end{Remark}

\subsection{Comparison with Banach space valued calculus via regularization}
\label{SubS:Comparison}

In the present subsection our aim is to make a link between functional It\^o calculus, as derived in this paper, and Banach space valued stochastic calculus via regularization for window processes, which has been conceived in \cite{DGR}, see also \cite{DGR2, digirrusso12, DGRnote},
and  \cite{digirfabbrirusso13} for more recent developments. More precisely, our purpose is to identify the building blocks of our functional It\^o's formula \eqref{E:Ito} with the terms appearing in the It\^o's formula derived in Theorem 6.3 and Section 7.2 in \cite{digirfabbrirusso13}. While it is expected that the vertical derivative $D^V\Uc$ can be identified with the term $D_{dx}^{\delta_0}\Uc$ of the Fr\'echet derivative, it is more difficult to guess to which terms the horizontal derivative $D^H\Uc$ corresponds. To clarify this latter point, in this subsection we derive two formulae which express $D^H\Uc$ in terms of Fr\'echet derivatives of $\Uc$.

\vspace{3mm}

Let us introduce some useful notations. We denote by $BV([-T,0])$ the set of c\`adl\`ag bounded variation functions on $[-T,0]$, which is a Banach space when equipped with the norm
\[
\|\eta\|_{BV([-T,0])} \ := \ |\eta(0)| + \|\eta\|_{\textup{Var}([-T,0])}, \qquad \eta\in BV([-T,0]),
\]
where $\|\eta\|_{\textup{Var}([-T,0])}=|d\eta|([-T,0])$ and $|d\eta|$ is the total variation measure associated to the measure $d\eta\in\Mc([-T,0])$ generated by $\eta$: $d\eta([-T,-t])=\eta(-t)-\eta(-T)$, $t\in[-T,0]$. We recall from subsection \ref{SubS:Background} that we extend $\eta\in BV([-T,0])$ to all $x\in\R$ setting $\eta(x) = 0$, $x<-T$, and $\eta(x)=\eta(0)$, $x\geq0$. Let us now introduce some useful facts about tensor products of Banach spaces.

\begin{Definition}
\label{D:Tensor}
Let $(E,\|\cdot\|_E)$ and $(F,\|\cdot\|_F)$ be two Banach spaces.\\
\textup{(i)} We shall denote by $E\otimes F$ the \textbf{algebraic tensor product} of $E$ and $F$, defined as the set of elements of the form $v = \sum_{i=1}^n e_i\otimes f_i$, for some positive integer $n$, where $e\in E$ and $f\in F$. The map $\otimes\colon E\times F\rightarrow E\otimes F$ is bilinear.\\
\textup{(ii)} We endow $E\otimes F$ with the \textbf{projective norm} $\pi$:
\[
\pi(v) \ := \ \inf\bigg\{\sum_{i=1}^n \|e_i\|_E\|f_i\|_F \ \colon \ v = \sum_{i=1}^n e_i\otimes f_i\bigg\}, \qquad \forall\,v\in E\otimes F.
\]
\textup{(iii)} We denote by $E\hat\otimes_\pi F$ the Banach space obtained as the completion of $E\otimes F$ for the norm $\pi$. We shall refer to $E\hat\otimes_\pi F$ as the \textbf{tensor product of the Banach spaces $E$ and $F$}.\\
\textup{(iv)} If $E$ and $F$ are Hilbert spaces, we denote $E\hat\otimes_h F$ the \textbf{Hilbert tensor product}, which is still a Hilbert space obtained as the completion of $E\otimes F$ for the scalar product $\langle e'\otimes f',e''\otimes f''\rangle := \langle e',e''\rangle_E\langle f',f''\rangle_F$, for any $e',e''\in E$ and $f',f''\in F$.\\
\textup{(v)} The symbols $E\hat\otimes_\pi^2$ and $e\otimes^2$ denote, respectively, the Banach space $E\hat\otimes_\pi E$ and the element $e\otimes e$ of the algebraic tensor product $E\otimes E$.
\end{Definition}

\begin{Remark}
\label{R:Identification}
{\rm
(i) The projective norm $\pi$ belongs to the class of the so-called \emph{reasonable crossnorms} $\alpha$ on $E\otimes F$, verifying $\alpha(e\otimes f)=\|e\|_E\|f\|_F$.\\
(ii) We notice, proceeding for example as in \cite{digirrusso12} (see, in particular, formula (2.1) in \cite{digirrusso12}; for more information on this subject we refer to \cite{ryan02}), that the dual $(E\hat\otimes_\pi F)^*$ of $E\hat\otimes_\pi F$ is isomorphic to the space of continuous bilinear forms $\Bc i(E,F)$, equipped with the norm $\|\cdot\|_{E,F}$ defined as
\[
\|\Phi\|_{E,F} \ := \ \sup_{\substack{e\in E, f\in F \\ \|e\|_E,\|f\|_F \leq 1}} |\Phi(e,f)|, \qquad \forall\,\Phi\in\Bc i(E,F).
\]
Moreover, there exists a canonical isomorphism between $\Bc i(E,F)$ and $L(E,F^*)$, the space of bounded linear operators from $E$ into $F^*$. Hence, we have the following chain of identifications: $(E\hat\otimes_\pi F)^* \cong \Bc i(E,F) \cong L(E;F^*)$.
\ep
}
\end{Remark}

\begin{Definition}
\label{D:C12Frechet}
Let $E$ be a Banach space. We say that $\Uc\colon E\rightarrow\R$ is of class $C^2(E)$ if
\begin{enumerate}
\item[(i)] $D \Uc$, the first Fr\'echet derivative of $\Uc$, belongs to $C(E; E^*)$ and
\item[(ii)] $D^2 \Uc$, the second Fr\'echet derivative of $\Uc$, belongs to $C(E; L(E;E^*))$.
\end{enumerate}
\end{Definition}

\begin{Remark}
{\rm
Take $E = C([-T,0])$ in Definition \ref{D:C12Frechet}.\\
(i) \emph{First Fr\'echet derivative $D\Uc$.} We have
\[
D\Uc \colon C([-T,0]) \ \longrightarrow \ (C([-T,0]))^* \cong \Mc([-T,0]).
\]
For every $\eta\in C([-T,0])$, we shall denote $D_{dx}\Uc(\eta)$ the unique measure in $\Mc([-T,0])$ such that
\[
D\Uc(\eta)\varphi \ = \ \int_{[-T,0]} \varphi(x) D_{dx}\Uc(\eta), \qquad \forall\,\varphi\in C([-T,0]).
\]
Notice that $\Mc([-T,0])$ can be represented as the direct sum: $\Mc([-T,0]) = \Mc_0([-T,0])\oplus\Dc_0$, where we recall that $\Mc_0([-T,0])$ is the subset of $\Mc([-T,0])$ of measures $\mu$ such that $\mu(\{0\})=0$, instead $\Dc_0$ (which is a shorthand for $\Dc_0([-T,0])$) denotes the one-dimensional space of measures which are multiples of the Dirac measure $\delta_0$. For every $\eta\in C([-T,0])$ we denote by $(D_{dx}^\perp \Uc(\eta),D_{dx}^{\delta_0}\Uc(\eta))$ the unique pair in $\Mc_0([-T,0])\oplus\Dc_0$ such that
\[
D_{dx}\Uc(\eta) \ = \ D_{dx}^\perp \Uc(\eta) + D_{dx}^{\delta_0}\Uc(\eta).
\]
(ii) \emph{Second Fr\'echet derivative $D^2\Uc$.} We have
\begin{align*}
D^2\Uc \colon C([-T,0]) \ \longrightarrow \ L(C([-T,0]);(C([-T,0]))^*) &\cong \Bc i(C([-T,0]),C([-T,0])) \\
&\cong (C([-T,0])\hat\otimes_\pi C([-T,0]))^*,
\end{align*}
where we used the identifications of Remark \ref{R:Identification}(iii). Let $\eta\in C([-T,0])$; 
a typical situation arises when there exists $D_{dx\,dy} \Uc(\eta)\in\Mc([-T,0]^2)$ such that $D^2\Uc(\eta)\in L(C([-T,0]);(C([-T,0]))^*)$ admits the representation
\[
D^2\Uc(\eta)(\varphi,\psi) \ = \ \int_{[-T,0]^2} \varphi(x)\psi(y) D_{dx\,dy}\Uc(\eta), \qquad \forall\,\varphi,\psi\in C([-T,0]).
\]
Moreover, $D_{dx\,dy}\Uc(\eta)$ is uniquely determined.
\ep
}
\end{Remark}
The definition below was given in \cite{DGR}.
\begin{Definition}
\label{D:Chi-subspace}
Let $E$ be a Banach space. A Banach subspace $(\chi,\|\cdot\|_\chi)$ continuously injected into $(E\hat\otimes_\pi^2)^*$, i.e., $\|\cdot\|_\chi\geq\|\cdot\|_{(E\hat\otimes_\pi^2)^*}$, will be called a \textbf{Chi-subspace} $($of $(E\hat\otimes_\pi^2)^*$$)$.
\end{Definition}

\begin{Remark}
\label{R:Chi-subspace}
{\rm
Take $E=C([-T,0])$ in Definition \ref{D:Chi-subspace}. 
As indicated in \cite{DGR}, a typical example of Chi-subspace of $C([-T,0])\hat\otimes_\pi^2$ is $\Mc([-T,0]^2)$ equipped with the usual total variation norm, denoted by $\|\cdot\|_{\text{Var}}$. Another important Chi-subspace of $C([-T,0])\hat\otimes_\pi^2$ is the following, which is also a Chi-subspace of $\Mc([-T,0]^2)$:
\begin{align*}
\chi_0 \ &:= \ \big\{\mu\in\Mc([-T,0]^2)\colon\mu(dx,dy) = g_1(x,y)dxdy + \lambda_1\delta_0(dx)\otimes\delta_0(dy) \\
&\quad \ + g_2(x)dx\otimes\lambda_2\delta_0(dy) + \lambda_3\delta_0(dx)\otimes g_3(y)dy + g_4(x)\delta_y(dx)\otimes dy, \\
&\quad \ g_1\in L^2([-T,0]^2),\,g_2,g_3\in L^2([-T,0]),\,g_4\in L^\infty([-T,0]),\,\lambda_1,\lambda_2,\lambda_3\in\R\big\}.
\end{align*}
Using the notations of Example 3.4 and Remark 3.5 in \cite{digirrusso12}, to which we refer for more details on this subject, we notice that $\chi_0$ is indeed given by the direct sum $\chi_0 = L^2([-T,0]^2) \oplus \big(L^2([-T,0])\hat\otimes_h \Dc_0\big) \oplus \big(\Dc_0 \hat\otimes_h L^2([-T,0])\big) \oplus \Dc_{0,0}([-T,0]^2) \oplus Diag([-T,0]^2)$. In the sequel, we shall refer to the term $g_4(x)\delta_y(dx)\otimes dy$ as the {\bf diagonal component} and to $g_4(x)$ as the {\bf diagonal element} of $\mu$. 
\ep
}
\end{Remark}

\noindent We can now state our first representation result for $D^H\Uc$.

\begin{Proposition}
\label{P:DH=Dacdeta}
Let $\Uc\colon C([-T,0])\rightarrow\R$ be continuously Fr\'{e}chet differentiable. Suppose the following.
\begin{enumerate}
\item[\textup{(i)}] For any $\eta\in C([-T,0])$ there exists $D_x^{\textup{ac}}\Uc(\eta)\in BV([-T,0])$ such that
\[
D_{dx}^\perp \Uc(\eta) \ = \ D_x^{\textup{ac}}\Uc(\eta)dx.
\]
\item[\textup{(ii)}] There exist continuous extensions $($necessarily unique$)$$:$
\[
u\colon\mathscr C([-T,0])\rightarrow\R, \qquad\qquad D_x^{\textup{ac}}u\colon\mathscr C([-T,0])\rightarrow BV([-T,0])
\]
of $\Uc$ and $D_x^{\textup{ac}}\Uc$, respectively.
\end{enumerate}
Then, for any $\eta\in C([-T,0])$,
\begin{equation}
\label{E:DH=Dacdeta}
D^H \Uc(\eta) \ = \ \int_{[-T,0]} D_x^{\textup{ac}} \Uc(\eta) d^+ \eta(x),
\end{equation}
where we recall that previous deterministic integral has been defined
in Section \ref{S211}.
In particular, the horizontal derivative $D^H \Uc(\eta)$ and the backward integral in \eqref{E:DH=Dacdeta} exist.
\end{Proposition}

\noindent\textbf{Proof.}
Let $\eta\in C([-T,0])$, then starting from the left-hand side of \eqref{E:DH=Dacdeta}, using the definition of $D^H\Uc(\eta)$, we are led to consider the following increment for the function $u$:
\begin{equation}
\label{E:FirstOrderDH}
\frac{u(\eta) - u(\eta(\cdot-\eps)1_{[-T,0[}+\eta(0)1_{\{0\}})}{\eps}.
\end{equation}
We shall expand \eqref{E:FirstOrderDH} using a Taylor's formula. Firstly, notice that, since $\Uc$ is $C^1$ Fr\'echet on $C([-T,0])$, for every $\eta_1\in C([-T,0])$, with $\eta_1(0)=\eta(0)$,  from the fundamental theorem of calculus we have
\[
\Uc(\eta) - \Uc(\eta_1) \ = \ \int_0^1 \bigg( \int_{-T}^0 D_x^{\text{ac}}\Uc(\eta + \lambda(\eta_1-\eta))(\eta(x)-\eta_1(x)) dx\bigg) d\lambda.
\]
Recalling from Remark \ref{R:Density} the density of $C_{\eta(0)}([-T,0])$ in $\mathscr C_{\eta(0)}([-T,0])$ with respect to the topology of $\mathscr C([-T,0])$, we deduce the following Taylor's formula for $u$:
\begin{equation}
\label{E:u-u1}
u(\eta) - u(\eta_1) \ = \ \int_0^1 \bigg( \int_{-T}^0 D_x^{\text{ac}}u(\eta + \lambda(\eta_1-\eta))(\eta(x)-\eta_1(x)) dx\bigg) d\lambda,
\end{equation}
for all $\eta_1\in\mathscr C_{\eta(0)}([-T,0])$. As a matter of fact, for any $\delta\in]0,T/2]$ let (similarly to Remark \ref{R:Density}(i))
\[
\eta_{1,\delta}(x) \ := \
\begin{cases}
\eta_1(x), \qquad &-T\leq x\leq-\delta, \\
\frac{1}{\delta}(\eta_1(0)-\eta_1(-\delta))x + \eta_1(0), &-\delta<x\leq0
\end{cases}
\]
and $\eta_{1,0}:=\eta_1$. Then $\eta_{1,\delta}\in C([-T,0])$, for any $\delta\in]0,T/2]$, and $\eta_{1,\delta}\rightarrow\eta_1$ in $\mathscr C([-T,0])$, as $\delta\rightarrow0^+$. Now, define $f\colon[-T,0]\times[0,1]\times[0,T/2]\rightarrow\R$ as follows
\[
f(x,\lambda,\delta) \ := \ D_x^{\text{ac}}u(\eta + \lambda(\eta_{1,\delta}-\eta))(\eta(x)-\eta_{1,\delta}(x)),
\]
for all $(x,\lambda,\delta)\in[-T,0]\times[0,1]\times[0,T/2]$. Notice that $f$ is continuous and hence bounded, since its domain is a compact set. Then, it follows from Lebesgue's dominated convergence theorem that
\begin{align*}
&\int_0^1 \bigg( \int_{-T}^0 D_x^{\text{ac}}\Uc(\eta + \lambda(\eta_{1,\delta}-\eta))(\eta(x)-\eta_{1,\delta}(x)) dx\bigg) d\lambda \\
&= \ \int_0^1 \bigg( \int_{-T}^0 f(x,\lambda,\delta) dx\bigg) d\lambda \ \overset{\delta\rightarrow0^+}{\longrightarrow} \ \int_0^1 \bigg( \int_{-T}^0 f(x,\lambda,0) dx\bigg) d\lambda \\
&= \ \int_0^1 \bigg( \int_{-T}^0 D_x^{\text{ac}}u(\eta + \lambda(\eta_1-\eta))(\eta(x)-\eta_1(x)) dx\bigg) d\lambda,
\end{align*}
from which we deduce \eqref{E:u-u1}, since $\Uc(\eta_{1,\delta})\rightarrow u(\eta_1)$ as $\delta\rightarrow0^+$. Taking $\eta_1(\cdot)=\eta(\cdot-\eps)1_{[-T,0[}+\eta(0)1_{\{0\}}$, we obtain
\begin{align*}
&\frac{u(\eta) - u(\eta(\cdot-\eps)1_{[-T,0[}+\eta(0)1_{\{0\}})}{\eps} \\
&= \int_0^1 \bigg( \int_{-T}^0 D_x^{\text{ac}} u\big(\eta + \lambda \big(\eta(\cdot-\eps)-\eta(\cdot)\big) 1_{[-T,0[}\big) \frac{\eta(x)-\eta(x-\eps)}{\eps} dx \bigg) d\lambda \\
&= \ I_1(\eta,\eps) + I_2(\eta,\eps)+ I_3(\eta,\eps),
\end{align*}
where
\begin{align*}
I_1(\eta,\eps) \ &:= \ \int_0^1 \bigg(\int_{-T}^0 \eta(x)\frac{1}{\eps} \Big( D_x^{\text{ac}} u\big(\eta + \lambda \big(\eta(\cdot-\eps)-\eta(\cdot)\big) 1_{[-T,0[}\big) \\
&\quad \ - D_{x+\eps}^{\text{ac}} u\big(\eta + \lambda \big(\eta(\cdot-\eps)-\eta(\cdot)\big) 1_{[-T,0[}\big) \Big) dx \bigg) d\lambda, \\
I_2(\eta,\eps) \ &:= \ \frac{1}{\eps}\int_0^1 \bigg( \int_{-\eps}^0 \eta(x) D_{x+\eps}^{\text{ac}} u\big(\eta + \lambda \big(\eta(\cdot-\eps)-\eta(\cdot)\big) 1_{[-T,0[}\big) dx \bigg) d\lambda, \\
I_3(\eta,\eps) \ &:= \ - \frac{1}{\eps} \int_0^1 \bigg( \int_{-T-\eps}^{-T} \eta(x) D_{x+\eps}^{\text{ac}} u\big(\eta + \lambda \big(\eta(\cdot-\eps)-\eta(\cdot)\big) 1_{[-T,0[}\big) dx \bigg) d\lambda.
\end{align*}
Notice that, since $\eta(x)=0$ for $x<-T$, we see that $I_2(\eta,\eps)=0$. Moreover, since $D_x^{\text{ac}} u(\cdot)=D_0^{\text{ac}} u(\cdot)$, for $x\geq 0$, and $\eta + \lambda (\eta(\cdot-\eps)-\eta(\cdot)) 1_{[-T,0[}\rightarrow\eta$ in $\mathscr C([-T,0])$ as $\eps\rightarrow0^+$, it follows that (using the continuity of $D_x^{\text{ac}}u$ from $\mathscr C([-T,0])$ into $BV([-T,0])$, which implies that $D_0^{\text{ac}}u(\eta + \lambda (\eta(\cdot-\eps)-\eta(\cdot)) 1_{[-T,0[})\rightarrow D_0^{\text{ac}}u(\eta)$ as $\eps\rightarrow0^+$)
\begin{align*}
&\frac{1}{\eps}\int_{-\eps}^0 \eta(x)D_{x+\eps}^{\text{ac}} u\big(\eta + \lambda \big(\eta(\cdot-\eps)-\eta(\cdot)\big) 1_{[-T,0[}\big) dx \notag \\
&= \ \frac{1}{\eps}\int_{-\eps}^0 \eta(x) dx\,D_0^{\text{ac}}u\big(\eta + \lambda \big(\eta(\cdot-\eps)-\eta(\cdot)\big) 1_{[-T,0[}\big) \ \overset{\eps\rightarrow0^+}{\longrightarrow} \ \eta(0)D_0^{\text{ac}}u(\eta).
\end{align*}
Finally, concerning $I_1(\eta,\eps)$, from Fubini's theorem we obtain (denoting $\eta_{\eps,\lambda} := \eta + \lambda (\eta(\cdot-\eps)-\eta(\cdot)) 1_{[-T,0[}$)
\begin{align*}
I_1(\eta,\eps) \ &= \ \int_0^1 \bigg(\int_{-T}^0 \eta(x)\frac{1}{\eps} \Big( D_x^{\text{ac}} u(\eta_{\eps,\lambda}) - D_{x+\eps}^{\text{ac}} u(\eta_{\eps,\lambda}) \Big) dx \bigg) d\lambda \\
&= \ -\int_0^1 \bigg(\int_{-T}^0 \eta(x)\frac{1}{\eps} \bigg( \int_{]x,x+\eps]} D_{dy}^{\text{ac}} u(\eta_{\eps,\lambda})\bigg) dx \bigg) d\lambda \\
&= \ -\int_0^1 \bigg(\int_{]-T,\eps]} \frac{1}{\eps} \bigg( \int_{(-T)\vee(y-\eps)}^{0\wedge y} \eta(x) dx\bigg) D_{dy}^{\text{ac}} u(\eta_{\eps,\lambda}) \bigg) d\lambda \\
&= \ I_{11}(\eta,\eps) + I_{12}(\eta,\eps),
\end{align*}
where
\begin{align*}
I_{11}(\eta,\eps) \ &:= \ -\int_0^1 \bigg(\int_{]-T,\eps]} \frac{1}{\eps} \bigg( \int_{(-T)\vee(y-\eps)}^{0\wedge y} \eta(x) dx\bigg) \Big(D_{dy}^{\text{ac}} u(\eta_{\eps,\lambda}) - D_{dy}^{\text{ac}} u(\eta)\Big) \bigg) d\lambda, \\
I_{12}(\eta,\eps) \ &:= \ -\int_0^1 \bigg(\int_{]-T,\eps]} \frac{1}{\eps} \bigg( \int_{(-T)\vee(y-\eps)}^{0\wedge y} \eta(x) dx\bigg) D_{dy}^{\text{ac}} u(\eta) \bigg) d\lambda \\
&= \ -\bigg(\int_{]-T,\eps]} \frac{1}{\eps} \bigg( \int_{(-T)\vee(y-\eps)}^{0\wedge y} \eta(x) dx\bigg) D_{dy}^{\text{ac}} u(\eta).
\end{align*}
Recalling that $D_x^{\text{ac}} u(\cdot)=D_0^{\text{ac}} u(\cdot)$, for $x\geq 0$, we see that in $I_{11}(\eta,\eps)$ and $I_{12}(\eta,\eps)$ the integrals on $]-T,\eps]$ are equal to the same integrals on $]-T,0]$, i.e.,
\begin{align*}
I_{11}(\eta,\eps) \ &= \ -\int_0^1 \bigg(\int_{]-T,0]} \frac{1}{\eps} \bigg( \int_{(-T)\vee(y-\eps)}^{0\wedge y} \eta(x) dx\bigg) \Big(D_{dy}^{\text{ac}} u(\eta_{\eps,\lambda}) - D_{dy}^{\text{ac}} u(\eta)\Big) \bigg) d\lambda \\
&= \ -\int_0^1 \bigg(\int_{]-T,0]} \frac{1}{\eps} \bigg( \int_{y-\eps}^y \eta(x) dx\bigg) \Big(D_{dy}^{\text{ac}} u(\eta_{\eps,\lambda}) - D_{dy}^{\text{ac}} u(\eta)\Big) \bigg) d\lambda, \\
I_{12}(\eta,\eps) \ &= \ -\int_{]-T,0]} \frac{1}{\eps} \bigg( \int_{(-T)\vee(y-\eps)}^{0\wedge y} \eta(x) dx\bigg) D_{dy}^{\text{ac}} u(\eta) \\
&= \ -\int_{]-T,0]} \frac{1}{\eps} \bigg( \int_{y-\eps}^y \eta(x) dx\bigg) D_{dy}^{\text{ac}} u(\eta).
\end{align*}
Now, observe that
\[
|I_{11}(\eta,\eps)| \ \leq \ \|\eta\|_\infty \|D_\cdot^\text{ac}u(\eta_{\eps,\lambda})-D_\cdot^\text{ac}u(\eta)\|_{\text{Var}([-T,0])} \ \overset{\eps\rightarrow0^+}{\longrightarrow} \ 0.
\]
Moreover, since $\eta$ is continuous at $y\in]-T,0]$, we deduce that $\int_{y-\eps}^y \eta(x) dx/\eps\rightarrow\eta(y)$ as $\eps\rightarrow0^+$. Therefore, by Lebesgue's dominated convergence theorem, we get
\[
I_{12}(\eta,\eps) \ \overset{\eps\rightarrow0^+}{\longrightarrow} \ -\int_{]-T,0]} \eta(y) D_{dy}^{\text{ac}} u(\eta).
\]
In conclusion, we have
\[
D^H\Uc(\eta) \ = \ \eta(0)D_0^{\text{ac}}u(\eta) - \int_{]-T,0]} \eta(y) D_{dy}^{\text{ac}} u(\eta),
\]
which gives \eqref{E:DH=Dacdeta} using the integration by parts formula \eqref{E:IbyP}.
\ep

\vspace{3mm}

For our second representation result of $D^H\Uc$ we need the following generalization of the deterministic backward integral when the integrand is a measure.

%


\begin{Definition}
\label{D:DeterministicIntegralMeasure}
Let $f\colon[-T,0]\rightarrow\R$ be a c\`{a}dl\`{a}g function and $g\in\Mc([-T,0])$. Suppose that the following limit
\[
\int_{[-T,0]}g(ds)d^+f(s) \ := \ \lim_{\eps\rightarrow0^+}\int_{[-T,0]} g(ds)\frac{f(s)-f(s-\eps)}{\eps},
\]
exists and it is finite. Then, the obtained quantity is denoted by $\int_{[-T,0]} gd^+f$ and called \textbf{$($deterministic, definite$)$ backward integral of $g$ with respect to $f$ $($on $[a,b]$$)$}.
\end{Definition}

\begin{Proposition}
\label{P:DH_SecondOrder}
Let $\eta\in C([-T,0])$ be such that the quadratic variation on $[-T,0]$ exists. Let $\Uc\colon C([-T,0])\rightarrow\R$ be twice continuously Fr\'echet differentiable such that
\[
D^2\Uc\colon C([-T,0]) \ \longrightarrow \ \chi_0\subset(C([-T,0])\hat\otimes_\pi C([-T,0]))^*\text{ continuously with respect to $\chi_0$.}
\]
Let us also suppose the following.
\begin{enumerate}
\item[\textup{(i)}] $D_x^{2,Diag}\Uc(\eta)$, the diagonal element of the second-order derivative at $\eta$, has a set of discontinuity which has null measure with respect to $[\eta]$ $($in particular, if it is countable$)$.
\item[\textup{(ii)}] There exist continuous extensions $($necessarily unique$)$$:$
\[
u\colon\mathscr C([-T,0])\rightarrow\R, \qquad\qquad D_{dx\,dy}^2u\colon\mathscr C([-T,0])\rightarrow\chi_0
\]
of $\Uc$ and $D_{dx\,dy}^2\Uc$, respectively.
\item[\textup{(iii)}] The horizontal derivative $D^H\Uc(\eta)$ exists at $\eta\in C([-T,0])$.
\end{enumerate}
Then
\begin{equation}
\label{E:DH=SecondOrder}
D^H \Uc(\eta) \ = \ \int_{[-T,0]} D_{dx}^\perp \Uc(\eta) d^+ \eta(x) - \frac{1}{2}\int_{[-T,0]} D_x^{2,Diag}\Uc(\eta) d[\eta](x).
\end{equation}
In particular, the backward integral in \eqref{E:DH=SecondOrder} exists.
\end{Proposition}
\textbf{Proof.} Let $\eta\in C([-T,0])$, then using the definition of $D^H\Uc(\eta)$ we are led to consider the following increment for the function $u$:
\begin{equation}
\label{E:ProofDH}
\frac{u(\eta) - u(\eta(\cdot-\eps)1_{[-T,0[}+\eta(0)1_{\{0\}})}{\eps},
\end{equation}
with $\eps>0$. Our aim is to expand \eqref{E:ProofDH} using a Taylor's formula. To this end, we begin noting that, since $\Uc$ is $C^2$ Fr\'echet, for every $\eta_1\in C([-T,0])$ the following standard Taylor's expansion holds:
\begin{align*}
\Uc(\eta_1) \ &= \ \Uc(\eta) + \int_{[-T,0]} D_{dx}\Uc(\eta) \big(\eta_1(x) - \eta(x)\big) \\
&\quad \ + \frac{1}{2}\int_{[-T,0]^2} D_{dx\,dy}^2\Uc(\eta) \big(\eta_1(x)-\eta(x)\big) \big(\eta_1(y)-\eta(y)\big) \notag \\
&\quad \ + \int_0^1 (1-\lambda)\bigg( \int_{[-T,0]^2} \Big( D_{dx\,dy}^2\Uc(\eta + \lambda (\eta_1-\eta)) \notag \\
&\quad \ - D_{dx\,dy}^2\Uc(\eta) \Big) \big(\eta_1(x)-\eta(x)\big) \big(\eta_1(y)-\eta(y)\big) \bigg) d\lambda. \notag
\end{align*}
Now, using the density of $C_{\eta(0)}([-T,0])$ into $\mathscr C_{\eta(0)}([-T,0])$ with respect to the topology of $\mathscr C([-T,0])$ and proceeding as in the proof of Proposition \ref{P:DH=Dacdeta}, we deduce the following Taylor's formula for $u$:
\begin{align}
\label{E:SecondOrder1}
& \frac{u(\eta) - u(\eta(\cdot-\eps)1_{[-T,0[}+\eta(0)1_{\{0\}})}{\eps} \\
&= \ \int_{[-T,0]} D_{dx}^\perp \Uc(\eta) \frac{\eta(x) - \eta(x-\eps)}{\eps} \notag \\
&\quad \ - \frac{1}{2}\int_{[-T,0]^2} D_{dx\,dy}^2 \Uc(\eta) \frac{(\eta(x)-\eta(x-\eps)) (\eta(y)-\eta(y-\eps))}{\eps} 1_{[-T,0[\times[-T,0[}(x,y) \notag \\
&\quad \ - \int_0^1 (1-\lambda)\bigg( \int_{[-T,0]^2} \Big( D_{dx\,dy}^2 u(\eta + \lambda (\eta(\cdot-\eps)-\eta(\cdot))1_{[-T,0[}) \notag \\
&\quad \ - D_{dx\,dy}^2 \Uc(\eta) \Big) \frac{(\eta(x)-\eta(x-\eps)) (\eta(y)-\eta(y-\eps))}{\eps} 1_{[-T,0[\times[-T,0[}(x,y) \bigg) d\lambda. \notag
\end{align}
Recalling the definition of $\chi_0$ given in Remark \ref{R:Chi-subspace}, we notice that (due to the presence of the indicator function $1_{[-T,0[\times[-T,0[}$)
\begin{align*}
&\int_{[-T,0]^2} D_{dx\,dy}^2 \Uc(\eta) \frac{(\eta(x)-\eta(x-\eps)) (\eta(y)-\eta(y-\eps))}{\eps} 1_{[-T,0[\times[-T,0[}(x,y) \\
&= \ \int_{[-T,0]^2} D_{x\,y}^{2,L^2} \Uc(\eta) \frac{(\eta(x)-\eta(x-\eps)) (\eta(y)-\eta(y-\eps))}{\eps} dx\,dy \\
&\quad \ + \int_{[-T,0]} D_x^{2,Diag} \Uc(\eta) \frac{(\eta(x)-\eta(x-\eps))^2}{\eps} dx,
\end{align*}
where, by hypothesis, the maps $\eta\in\mathscr C([-T,0])\mapsto D_{x\,y}^{2,L^2} u(\eta)\in L^2([-T,0]^2)$ and $\eta\in\mathscr C([-T,0])\mapsto D_{x}^{2,Diag}u(\eta)\in L^\infty([-T,0])$ are continuous. In particular, \eqref{E:SecondOrder1} becomes
\begin{equation}
\label{E:SecondOrder2}
\frac{u(\eta) - u(\eta(\cdot-\eps)1_{[-T,0[}+\eta(0)1_{\{0\}})}{\eps} \ = \ I_1(\eps) + I_2(\eps) + I_3(\eps) + I_4(\eps) + I_5(\eps),
\end{equation}
where
\begin{align*}
I_1(\eps) \ &:= \ \int_{[-T,0]} D_{dx}^\perp \Uc(\eta) \frac{\eta(x) - \eta(x-\eps)}{\eps}, \\
I_2(\eps) \ &:= \ - \frac{1}{2}\int_{[-T,0]^2} D_{x\,y}^{2,L^2} \Uc(\eta) \frac{(\eta(x)-\eta(x-\eps)) (\eta(y)-\eta(y-\eps))}{\eps} dx\,dy, \\
I_3(\eps) \ &:= \ - \frac{1}{2}\int_{[-T,0]} D_x^{2,Diag} \Uc(\eta) \frac{(\eta(x)-\eta(x-\eps))^2}{\eps} dx, \\
I_4(\eps) \ &:= \ - \int_0^1 (1-\lambda)\bigg( \int_{[-T,0]^2} \Big( D_{x\,y}^{2,L^2} u(\eta + \lambda (\eta(\cdot-\eps)-\eta(\cdot))1_{[-T,0[}) \\
&\quad \ - D_{x\,y}^{2,L^2} \Uc(\eta) \Big) \frac{(\eta(x)-\eta(x-\eps)) (\eta(y)-\eta(y-\eps))}{\eps} dx\,dy \bigg) d\lambda, \\
I_5(\eps) \ &:= \ - \int_0^1 (1-\lambda)\bigg( \int_{[-T,0]} \Big( D_x^{2,Diag} u(\eta + \lambda (\eta(\cdot-\eps)-\eta(\cdot))1_{[-T,0[}) \\
&\quad \ - D_x^{2,Diag} \Uc(\eta) \Big) \frac{(\eta(x)-\eta(x-\eps))^2}{\eps} dx \bigg) d\lambda.
\end{align*}
Firstly, we shall prove that
\begin{equation}
\label{E:I2-->0}
I_2(\eps) \ \overset{\eps\rightarrow0^+}{\longrightarrow} \ 0.
\end{equation}
To this end, for every $\eps>0$, define the operator $T_\eps\colon L^2([-T,0]^2)\rightarrow\R$ as follows:
\[
T_\eps \, g \ = \ \int_{[-T,0]^2} g(x,y) \frac{(\eta(x)-\eta(x-\eps)) (\eta(y)-\eta(y-\eps))}{\eps} dx\,dy, \qquad \forall\, g\in L^2([-T,0]^2).
\]
Then $T_\eps\in L^2([-T,0])^*$. Indeed, from Cauchy-Schwarz inequality,
\begin{align*}
|T_\eps\,g| \ &\leq \ \|g\|_{L^2([-T,0]^2)} \sqrt{\int_{[-T,0]^2}\frac{(\eta(x)-\eta(x-\eps))^2 (\eta(y)-\eta(y-\eps))^2}{\eps^2}dx\,dy} \\
&= \ \|g\|_{L^2([-T,0]^2)} \int_{[-T,0]} \frac{(\eta(x)-\eta(x-\eps))^2}{\eps} dx
\end{align*}
and this last quantity is  bounded with respect to $\eps$ since the quadratic variation of $\eta$ on $[-T,0]$ exists. In particular, we have proved that for every $g\in L^2([-T,0]^2)$ there exists a constant $M_g\geq0$ such that
\[
\sup_{0 < \eps < 1} |T_\eps\,g| \ \leq \ M_g.
\]
It follows from Banach-Steinhaus theorem that there exists a constant $M\geq0$ such that
\begin{equation}
\label{E:BanachSteinhaus}
\sup_{0 < \eps < 1} \|T_\eps\|_{L^2([-T,0])^*} \ \leq \ M.
\end{equation}
Now, let us consider the set $\Sc := \{g\in L^2([-T,0]^2)\colon g(x,y) = e(x)f(y),\text{ with }e,f\in C^1([-T,0])\}$, which is dense in $L^2([-T,0]^2)$. Let us show that
\begin{equation}
\label{E:T_eps-->0}
T_\eps\,g \ \overset{\eps\rightarrow0^+}{\longrightarrow} \ 0, \qquad \forall\,g\in\Sc.
\end{equation}
Fix $g\in\Sc$, with $g(x,y)=e(x)f(y)$ for any $(x,y)\in[-T,0]$, then
\begin{equation}
\label{E:T_eps_Sc}
T_\eps\,g \ = \ \frac{1}{\eps} \int_{[-T,0]} e(x) \big(\eta(x) - \eta(x-\eps)\big) dx \int_{[-T,0]} f(y) \big(\eta(y) - \eta(y-\eps)\big) dy.
\end{equation}
We have
\begin{align*}
&\bigg|\int_{[-T,0]} e(x) \big(\eta(x) - \eta(x-\eps)\big) dx\bigg| \ = \ \bigg|\int_{[-T,0]} \big(e(x) - e(x+\eps)\big) \eta(x) dx \\
&- \int_{[-T-\eps,-T]}e(x+\eps)\eta(x) dx + \int_{[-\eps,0]} e(x+\eps)\eta(x) dx\bigg| \\
&\leq \ \eps \bigg(\int_{[-T,0]}|\dot e(x)|dx + 2\|e\|_\infty\bigg)\|\eta\|_\infty.
\end{align*}
Similarly,
\[
\bigg|\int_{[-T,0]} f(y) \big(\eta(y) - \eta(y-\eps)\big) dy\bigg| \ \leq \ \eps \bigg(\int_{[-T,0]}|\dot f(y)|dy + 2\|f\|_\infty\bigg)\|\eta\|_\infty.
\]
Therefore, from \eqref{E:T_eps_Sc} we find
\[
|T_\eps\,g| \ \leq \ \eps \bigg(\int_{[-T,0]}|\dot e(x)|dx + 2\|e\|_\infty\bigg) \bigg(\int_{[-T,0]}|\dot f(y)|dy + 2\|f\|_\infty\bigg) \|\eta\|_\infty^2,
\]
which converges to zero as $\eps$ goes to zero and therefore \eqref{E:T_eps-->0} is established. This in turn implies that
\begin{equation}
\label{E:T_eps-->0bis}
T_\eps\,g \ \overset{\eps\rightarrow0^+}{\longrightarrow} \ 0, \qquad \forall\,g\in L^2([-T,0]^2).
\end{equation}
Indeed, fix $g\in L^2([-T,0]^2)$ and let $(g_n)_n\subset\Sc$ be such that $g_n\rightarrow g$ in $L^2([-T,0]^2)$. Then
\[
|T_\eps\,g| \ \leq \ |T_\eps(g-g_n)| + |T_\eps\,g_n| \ \leq \ \|T_\eps\|_{L^2([-T,0]^2)^*} \|g-g_n\|_{L^2([-T,0]^2)} + |T_\eps\,g_n|.
\]
From \eqref{E:BanachSteinhaus} it follows that
\[
|T_\eps\,g| \ \leq \ M \|g-g_n\|_{L^2([-T,0]^2)} + |T_\eps\,g_n|,
\]
which implies $\limsup_{\eps\rightarrow0^+}|T_\eps\,g| \leq M \|g-g_n\|_{L^2([-T,0]^2)}$. Sending $n$ to infinity, we deduce \eqref{E:T_eps-->0bis} and finally \eqref{E:I2-->0}.

Let us now consider the term $I_3(\eps)$ in \eqref{E:SecondOrder2}. Since the quadratic variation $[\eta]$ exists, it follows from Portmanteau's theorem and hypothesis (i) that
\[
I_3(\eps) \ = \ \int_{[-T,0]} D_x^{2,Diag} \Uc(\eta) \frac{(\eta(x)-\eta(x-\eps))^2}{\eps} dx \ \underset{\eps\rightarrow0^+}{\longrightarrow} \ \int_{[-T,0]} D_x^{2,Diag}\Uc(\eta) d[\eta](x).
\]
Regarding the term $I_4(\eps)$ in \eqref{E:SecondOrder2}, let $\phi_\eta\colon[0,1]^2\rightarrow L^2([-T,0]^2)$ be given by
\[
\phi_\eta(\eps,\lambda)(\cdot,\cdot) \ = \ D_{\cdot\,\cdot}^{2,L^2} u\big(\eta + \lambda (\eta(\cdot-\eps)-\eta(\cdot))1_{[-T,0[}\big).
\]
By hypothesis,  $\phi_\eta$ is a continuous map, and hence it is uniformly continuous, since $[0,1]^2$ is a compact set. Let $\rho_{\phi_\eta}$ denote the continuity modulus of $\phi_\eta$, then
\begin{align*}
&\big\|D_{\cdot\,\cdot}^{2,L^2} u\big(\eta + \lambda (\eta(\cdot-\eps)-\eta(\cdot))1_{[-T,0[}\big) - D_{\cdot\,\cdot}^{2,L^2} \Uc(\eta)\big\|_{L^2([-T,0]^2)} \\
&= \ \|\phi_\eta(\eps,\lambda) - \phi_\eta(0,\lambda)\|_{L^2([-T,0]^2)} \ \leq \ \rho_{\phi_\eta}(\eps).
\end{align*}
This implies, by Cauchy-Schwarz inequality,
\begin{align*}
&\bigg|\int_0^1 (1-\lambda)\bigg( \int_{[-T,0]^2} \Big( D_{x\,y}^{2,L^2} u(\eta + \lambda (\eta(\cdot-\eps)-\eta(\cdot))1_{[-T,0[}) \\
&- D_{x\,y}^{2,L^2} \Uc(\eta) \Big) \frac{(\eta(x)-\eta(x-\eps)) (\eta(y)-\eta(y-\eps))}{\eps} dx\,dy \bigg) d\lambda\bigg| \\
&\leq \ \int_0^1(1-\lambda)\big\|D_{\cdot\,\cdot}^{2,L^2}u(\eta+\lambda(\eta(\cdot-\eps)-\eta(\cdot))1_{[-T,0]}) \\
&- D_{\cdot\,\cdot}^{2,L^2}\Uc(\eta)\big\|_{L^2([-T,0]^2)}\sqrt{\int_{[-T,0]^2}\frac{(\eta(x)-\eta(x-\eps))^2 (\eta(y)-\eta(y-\eps))^2}{\eps^2}dx\,dy}\,d\lambda \\
&\leq \ \int_0^1(1-\lambda)\rho_{\phi_\eta}(\eps) \bigg(\int_{[-T,0]} \frac{(\eta(x)-\eta(x-\eps))^2}{\eps} dx \bigg)d\lambda \\
&= \ \frac{1}{2}\rho_{\phi_\eta}(\eps) \int_{[-T,0]} \frac{(\eta(x)-\eta(x-\eps))^2}{\eps} dx \ \overset{\eps\rightarrow0^+}{\longrightarrow} \ 0.
\end{align*}
Finally, we consider the term $I_5(\eps)$ in \eqref{E:SecondOrder2}. Define $\psi_\eta\colon[0,1]^2\rightarrow L^\infty([-T,0])$ as follows:
\[
\psi_\eta(\eps,\lambda)(\cdot) \ = \ D_\cdot^{2,Diag} u\big(\eta + \lambda (\eta(\cdot-\eps)-\eta(\cdot))1_{[-T,0[}\big).
\]
We see that $\psi_\eta$ is uniformly continuous. Let $\rho_{\psi_\eta}$ denote the continuity modulus of $\psi_\eta$, then
\begin{align*}
&\big\|D_\cdot^{2,Diag} u\big(\eta + \lambda (\eta(\cdot-\eps)-\eta(\cdot))1_{[-T,0[}\big) - D_\cdot^{2,Diag} \Uc(\eta)\big\|_{L^\infty([-T,0])} \\
&= \ \|\psi_\eta(\eps,\lambda) - \psi_\eta(0,\lambda)\|_{L^\infty([-T,0])} \ \leq \ \rho_{\psi_\eta}(\eps).
\end{align*}
Therefore, we have
\begin{align*}
&\bigg|\int_0^1 (1-\lambda)\bigg( \int_{[-T,0]} \Big( D_x^{2,Diag} u(\eta + \lambda (\eta(\cdot-\eps)-\eta(\cdot))1_{[-T,0[}) \\
&- D_x^{2,Diag} \Uc(\eta) \Big) \frac{(\eta(x)-\eta(x-\eps))^2}{\eps} dx \bigg) d\lambda\bigg| \\
&\leq \ \int_0^1(1-\lambda)\bigg(\int_{[-T,0]}\rho_{\psi_\eta}(\eps) \frac{(\eta(x)-\eta(x-\eps))^2}{\eps} dx \bigg) d\lambda \\
&= \ \frac{1}{2}\rho_{\psi_\eta}(\eps)\int_{[-T,0]} \frac{(\eta(x)-\eta(x-\eps))^2}{\eps} dx \ \overset{\eps\rightarrow0^+}{\longrightarrow} \ 0.
\end{align*}
In conclusion, we have proved that all the integral terms in the right-hand side of \eqref{E:SecondOrder2}, unless $I_1(\eps)$, admit a limit when $\eps$ goes to zero. Since the left-hand side admits a limit, namely $D^H\Uc(\eta)$, we deduce that the backward integral
\[
I_1(\eps) \ = \ \int_{[-T,0]} D_{dx}^\perp \Uc(\eta) \frac{\eta(x) - \eta(x-\eps)}{\eps} \ \overset{\eps\rightarrow0^+}{\longrightarrow} \ \int_{[-T,0]} D_{dx}^\perp \Uc(\eta) d^+\eta(x)
\]
exists and it is finite, which concludes the proof.
\ep

\section{Strong-viscosity solutions to path-dependent PDEs}
\label{S:StrongViscositySolutions}

\setcounter{equation}{0} \setcounter{Assumption}{0}
\setcounter{Theorem}{0} \setcounter{Proposition}{0}
\setcounter{Corollary}{0} \setcounter{Lemma}{0}
\setcounter{Definition}{0} \setcounter{Remark}{0}

In the present section, we study the path-dependent nonlinear Kolmogorov equation:
\begin{equation}
\label{E:KolmEq}
\begin{cases}
\partial_t\Uc + D^H\Uc + \frac{1}{2}D^{VV}\Uc \ = \ F(t,\eta,\Uc,D^V\Uc), \;\;\; &\forall\,(t,\eta)\in[0,T[\times C([-T,0]), \\
\Uc(T,\eta) \ = \ G(\eta), &\forall\,\eta\in C([-T,0]),
\end{cases}
\end{equation}
where $G\colon C([-T,0])\rightarrow\R$ and $F\colon[0,T]\times C([-T,0])\times\R\times\R\rightarrow\R$ are Borel measurable functions. Firstly, we provide a definition of classical solution to equation \eqref{E:KolmEq}. Then, motivated by a significant hedging example, we introduce a concept of weak (strong-viscosity) solution.

\subsection{Path-dependent Kolmogorov equation: classical solutions}
\label{SubS:ClassicalSolutions}

In this subsection, we give the definition of classical solution to the path-dependent nonlinear Kolmogorov equation \eqref{E:KolmEq} and provide a uniqueness result using BSDE methods. We conclude proving existence in the case $F\equiv0$, arising for example in hedging problems of path-dependent contingent claims.

\begin{Definition}
A function $\Uc\colon[0,T]\times C([-T,0])\rightarrow\R$ in $C^{1,2}(([0,T[\times\text{past})\times\text{present})\cap C([0,T]\times C([-T,0]))$, which solves equation \eqref{E:KolmEq}, is called a \textbf{classical solution} to the path-dependent nonlinear Kolmogorov equation \eqref{E:KolmEq}.
\end{Definition}

\noindent To prove uniqueness we need to introduce some additional notations. Let $(\Omega,\Fc,\P)$ be a complete probability space on which a real Brownian motion $W=(W_t)_{t\geq0}$ is defined. Let $\F=(\Fc_t)_{t\geq0}$ denote the completion of the natural filtration generated by $W$.

\begin{itemize}
\item $\S^p(t,T)$, $p\geq1$, $0 \leq t \leq T$, the set  of real c\`adl\`ag $\F$-predictable processes $Y=(Y_s)_{t\leq s\leq T}$ such that
\[
\|Y\|_{_{\S^p(t,T)}}^p := \ \E\Big[ \sup_{t\leq s\leq T} |Y_s|^p \Big] \ < \ \infty.
\]
\item $\H^p(t,T)^d$, $p$ $\geq$ $1$, $0 \leq t \leq T$, the set of  $\R^d$-valued predictable processes $Z=(Z_s)_{t\leq s\leq T}$ such that
\[
\|Z\|_{_{\H^p(t,T)^d}}^p := \ \E\bigg[\bigg(\int_t^T |Z_s|^2 ds\bigg)^{\frac{p}{2}}\bigg] \ < \ \infty.
\]
We simply write $\H^p(t,T)$ when $d=1$.
\item $\A^{+,2}(t,T)$, $0 \leq t \leq T$, the  set of real nondecreasing predictable processes $K$ $=$ $(K_s)_{t\leq s\leq T}\in\S^2(t,T)$ with $K_t$ $=$ $0$, so that
\[
\|K\|_{_{\S^2(t,T)}}^2 := \ \E\big[|K_T|^2\big].
\]
\item $\L^p(t,T;\R^m)$, $p\geq1$, $0 \leq t \leq T$, the set of $\R^m$-valued $\F$-predictable processes $\phi = (\phi_s)_{t \leq s \leq T}$ such that
\[
\|\phi\|_{_{\L^p(t,T;\R^m)}}^p := \ \E\bigg[\int_t^T |\phi_s|^p ds\bigg] \ < \ \infty.
\]
\end{itemize}

\begin{Definition}
\label{E:StochasticFlow}
Let $t\in[0,T]$ and $\eta\in C([-T,0])$. Then, we define the \textbf{stochastic flow}:
\[
\mathbb W_s^{t,\eta}(x) \ = \
\begin{cases}
\eta(x+s-t), &-T \leq x \leq t-s, \\
\eta(0) + W_{x+s} - W_t, \qquad &t-s < x \leq 0,
\end{cases}
\]
for any $t \leq s \leq T$.
\end{Definition}

\begin{Theorem}
\label{T:UniqClassical}
Let $G\colon C([-T,0])\rightarrow\R$ and $F\colon[0,T]\times C([-T,0])\times\R\times\R\rightarrow\R$ be Borel measurable functions satisfying, for some positive constants $C$ and $m$,
\begin{align*}
|F(t,\eta,y,z) - F(t,\eta,y',z')| \ &\leq \ C\big(|y-y'| + |z-z'|\big), \\
|G(\eta)| + |F(t,\eta,0,0)| \ &\leq \ C\big(1 + \|\eta\|_\infty^m\big),\end{align*}
for all $(t,\eta)\in[0,T]\times C([-T,0])$, $y,y'\in\R$, and $z,z'\in\R$. Let $\Uc\colon[0,T]\times C([-T,0])\rightarrow\R$ be a classical solution to equation \eqref{E:KolmEq}, satisfying the polynomial growth condition:
\begin{equation}
\label{E:PolGrowth_U}
|\Uc(t,\eta)| \ \leq \ C\big(1 + \|\eta\|_\infty^m\big), \qquad \forall\,(t,\eta)\in[0,T]\times C([-T,0]).
\end{equation}
Then, we have
\[
\Uc(t,\eta) \ = \ Y_t^{t,\eta}, \qquad \forall\,(t,\eta)\in[0,T]\times C([-T,0]),
\]
where $(Y_s^{t,\eta},Z_s^{t,\eta})_{s\in[t,T]} = (\Uc(s,\mathbb W_s^{t,\eta}),D^V\Uc(s,\mathbb W_s^{t,\eta})1_{[t,T[}(s))_{s\in[t,T]}\in\S^2(t,T)\times\H^2(t,T)$ is the solution to the backward stochastic differential equation, $\P$-a.s.,
\[
Y_s^{t,\eta} \ = \ G(\mathbb W_T^{t,\eta}) + \int_s^T F(r,\mathbb W_r^{t,\eta},Y_r^{t,\eta},Z_r^{t,\eta}) dr - \int_s^T Z_r^{t,\eta} dW_r, \qquad t \leq s \leq T.
\]
In particular, there exists at most one classical solution to the path-dependent nonlinear Kolmogorov equation \eqref{E:KolmEq}.
\end{Theorem}
\textbf{Proof.}
Fix $(t,\eta)\in[0,T[\times C([-T,0])$ and set, for all $t \leq s \leq T$,
\[
Y_s^{t,\eta} \ = \ \Uc(s,\mathbb W_s^{t,\eta}), \qquad Z_s^{t,\eta} \ = \ D^V\Uc(s,\mathbb W_s^{t,\eta})1_{[t,T[}(s).
\]
Then, for any $T_0\in[t,T[$, applying It\^o's formula \eqref{E:ItoTime} to $\Uc(s,\mathbb W_s^{t,\eta})$ and using the fact that $\Uc$ solves equation \eqref{E:KolmEq}, we find, $\P$-a.s.,
\begin{equation}
\label{E:BSDEIto}
Y_s^{t,\eta} \ = \ Y_{T_0}^{t,\eta} + \int_s^{T_0} F(r,\mathbb W_r^{t,\eta},Y_r^{t,\eta},Z_r^{t,\eta}) dr - \int_s^{T_0} Z_r^{t,\eta} dW_r, \qquad t \leq s \leq T_0.
\end{equation}
The thesis would follow if we could pass to the limit in \eqref{E:BSDEIto} as $T_0\rightarrow T$. To do this, we notice that it follows from Proposition \ref{P:EstimateBSDEAppendix} that there exists a positive constant $c$, depending only on $T$ and the constants $C$ and $m$ appearing in the statement of the present Theorem \ref{T:UniqClassical}, such that
\[
\E\int_t^{T_0} |Z_s^{t,\eta}|^2 ds \ \leq \ c\|Y^{t,\eta}\|_{\S^2(t,T)}^2 + c\E\int_t^T |F(r,\mathbb W_r^{t,\eta},0,0)|^2 dr, \qquad \forall\,T_0\in[t,T[.
\]
We recall that, for any $q\geq1$,
\begin{equation}
\label{E:EstimateSupW}
\E\Big[\sup_{t\leq s\leq T}\|\mathbb W_s^{t,\eta}\|_\infty^q\Big] \ < \ \infty.
\end{equation}
Notice that from \eqref{E:PolGrowth_U} and \eqref{E:EstimateSupW} we have $\|Y^{t,\eta}\|_{\S^2(t,T)}<\infty$, so that $Y\in\S^2(t,T)$. Then, from monotone convergence theorem we find
\[
\E\int_t^T |Z_s^{t,\eta}|^2 ds \ \leq \ c\|Y^{t,\eta}\|_{\S^2(t,T)}^2 + c\E\int_t^T |F(r,\mathbb W_r^{t,\eta},0,0)|^2 dr.
\]
Therefore, it follows from the polynomial growth condition of $F$ and \eqref{E:EstimateSupW} that $Z\in\H^2(t,T)$. This implies, using the Lipschitz
character of $F$ in $(y,z)$, that $\E\int_t^T |F(r,\mathbb W_r^{t,\eta},Y_r^{t,\eta},Z_r^{t,\eta})|^2 dr<\infty$, so that we can pass to the limit in \eqref{E:BSDEIto} and we get the thesis.
\ep

\vspace{3mm}

We conclude this subsection with an existence result for the \emph{path-dependent heat equation}, namely for the path-dependent nonlinear Kolmogorov equation \eqref{E:KolmEq} with $F\equiv0$.

\begin{Theorem}
\label{T:ExistenceClassical}
Let $F\equiv0$ and $G\colon C([-T,0])\rightarrow\R$ be given by,  for all $\eta\in C([-T,0])$
\begin{equation}
\label{E:G_Cylindrical}
G(\eta) \ = \ g\bigg(\int_{[-T,0]}\varphi_1(x+T)d^-\eta(x),\ldots,\int_{[-T,0]}\varphi_N(x+T)d^-\eta(x)\bigg),
\end{equation}
for some functions $g\in C_p^2(\R^N)$ $($$g$ and its first and second derivatives are continuous and have polynomial growth$)$ and $\varphi_1,\ldots,\varphi_N\in C^2([0,T])$, with $N\in\N\backslash\{0\}$ and $\varphi_i(x)=0$ for any $x\in\R\backslash[0,T]$. Then, there exists a unique classical solution $\Uc$ to the path-dependent heat equation \eqref{E:KolmEq}, which is given by
\[
\Uc(t,\eta) \ = \ \E\big[G(\mathbb W_T^{t,\eta})\big], \qquad \forall\,(t,\eta)\in[0,T]\times C([-T,0]).
\]
\end{Theorem}
\textbf{Proof.}
Let us consider the function $\Uc\colon[0,T]\times C([-T,0])\rightarrow\R$ given by, for all $(t,\eta)\in[0,T]\times C([-T,0])$,
\begin{align*}
\Uc(t,\eta) \ &= \ \E\big[G(\mathbb W_T^{t,\eta})\big] \\
&= \ \E\bigg[g\bigg(\int_{[-t,0]}\varphi_1(x+t)d^-\eta(x) + \int_t^T\varphi_1(s)dW_s,\ldots\bigg)\bigg] \\
&= \ \Psi\bigg(t,\int_{[-t,0]}\varphi_1(x+t)d^-\eta(x),\ldots,\int_{[-t,0]}\varphi_N(x+t)d^-\eta(x)\bigg),
\end{align*}
where
\[
\Psi(t,x_1,\ldots,x_N) \ = \ \E\bigg[g\bigg(x_1+\int_t^T\varphi_1(s)dW_s,\ldots,x_N+\int_t^T\varphi_N(s)dW_s\bigg)\bigg],
\]
for any $(t,x_1,\ldots,x_N)\in[0,T]\times\R^N$. Notice that, for any $i,j=1,\ldots,N$,
\begin{align*}
D_{x_i}\Psi(t,x_1,\ldots,x_N) \ &= \ \E\bigg[D_{x_i} g\bigg(x_1+\int_t^T\varphi_1(s)dW_s,\ldots,x_N+\int_t^T\varphi_N(s)dW_s\bigg)\bigg], \\
D_{x_ix_j}^2\Psi(t,x_1,\ldots,x_N) \ &= \ \E\bigg[D_{x_ix_j}^2 g\bigg(x_1+\int_t^T\varphi_1(s)dW_s,\ldots,x_N+\int_t^T\varphi_N(s)dW_s\bigg)\bigg],
\end{align*}
so that $\Psi$ and its first and second spatial derivatives are continuous on $[0,T]\times\R^N$. Let us focus on the time derivative $\partial_t\Psi$ of $\Psi$. We have, for any $h>0$ such that $t+h\in[0,T]$,
\begin{align*}
&\frac{\Psi(t+h,x_1,\ldots,x_N)-\Psi(t,x_1,\ldots,x_N)}{h} \\
&= \ \frac{1}{h}\E\bigg[g\bigg(x_1+\int_{t+h}^T\varphi_1(s)dW_s,\ldots\bigg)-g\bigg(x_1+\int_t^T\varphi_1(s)dW_s,\ldots\bigg)\bigg].
\end{align*}
Then, using a standard Taylor's formula, we find
\begin{align}
\label{E:Malliavin1}
&\frac{\Psi(t+h,x_1,\ldots,x_N)-\Psi(t,x_1,\ldots,x_N)}{h} \\
&= \ -\frac{1}{h}\E\bigg[\int_0^1 \sum_{i=1}^N D_{x_i}g\bigg(x_1+\int_t^T\varphi_1(s)dW_s - \alpha\int_t^{t+h}\varphi_1(s)dW_s,\ldots\bigg)\int_t^{t+h}\varphi_i(s)dW_sd\alpha\bigg]. \notag
\end{align}
Now, it follows from the integration by parts formula of Malliavin calculus, see, e.g., formula (1.42) in \cite{nualart} (taking into account that It\^o integrals are Skorohod integrals), that, for any $i=1,\ldots,N$,
\begin{align}
\label{E:Malliavin1>0}
&\E\bigg[D_{x_i}g\bigg(x_1+\int_t^T\varphi_1(s)\big(1-\alpha1_{[t,t+h]}(s)\big)dW_s,\ldots\bigg)\int_t^{t+h}\varphi_i(s)dW_s\bigg] \\
&= \ (1-\alpha)\E\bigg[\sum_{j=1}^N D_{x_ix_j}^2 g\bigg(x_1+\int_t^T\varphi_1(s)\big(1-\alpha1_{[t,t+h]}(s)\big)dW_s,\ldots\bigg)\int_t^{t+h}\varphi_i(s)\varphi_j(s)ds\bigg]. \notag
\end{align}
Then, plugging \eqref{E:Malliavin1>0} into \eqref{E:Malliavin1} and letting $h\rightarrow0^+$, we get (recalling that $D_{x_ix_j}^2 g$ has polynomial growth, for any $i,j$)
\begin{equation}
\label{E:Repr_partial_t_Psi}
\partial_t^+\Psi(t,x_1,\ldots,x_N) \ = \ -\frac{1}{2}\E\bigg[\sum_{i,j=1}^N D_{x_ix_j}^2 g\bigg(x_1+\int_t^T\varphi_1(s)dW_s,\ldots\bigg)\varphi_i(t)\varphi_j(t)\bigg],
\end{equation}
for any $(t,x_1,\ldots,x_N)\in[0,T[\times\R^N$, where $\partial_t^+\Psi$ denotes the right-time derivative of $\Psi$. Since $\Psi$ and $\partial_t^+\Psi$ are continuous, we deduce that $\partial_t\Psi$ exists and is continuous on $[0,T[$ (see for example Corollary 1.2, Chapter 2, in \cite{pazy83}). Moreover, from the representation formula \eqref{E:Repr_partial_t_Psi} we see that $\partial_t\Psi$ exists and is continuous up to time $T$. Furthermore, from the expression of $D_{x_ix_j}^2 \Psi$, we see that
\[
\partial_t\Psi(t,x_1,\ldots,x_N) \ = \ -\frac{1}{2}\sum_{i,j=1}^N \varphi_i(t)\varphi_j(t) D_{x_ix_j}^2 \Psi(t,x_1,\ldots,x_N).
\]
Therefore, $\Psi\in C^{1,2}([0,T]\times\R^N)$ and is a classical solution to the Cauchy problem:
\begin{equation}
\label{E:HeatEquationPsi}
\begin{cases}
\partial_t \Psi(t,\mathbf x) + \frac{1}{2}\sum_{i,j=1}^N \varphi_i(t)\varphi_j(t) D_{x_ix_j}^2 \Psi(t,\mathbf x) = 0, \qquad\qquad &\forall\,(t,\mathbf x)\in[0,T[\times\R^N, \\
\Psi(t,\mathbf x) = g(\mathbf x), &\forall\,\mathbf x\in\R^N.
\end{cases}
\end{equation}
Now we express the derivatives of $\Uc$ in terms of $\Psi$. We begin noting that, taking into account Proposition \ref{P:BVI}, for each $i$ and $t\in[0,T]$, the linear functional $T_{i,t}\colon C([-T,0])\rightarrow\R$ is given by
\[
T_{i,t}\eta \ = \ \int_{[-t,0]}\varphi_i(x+t)d^-\eta(x) \ = \ \eta(0)\varphi_i(t) - \int_{-t}^0 \eta(x)\dot{\varphi}_i(x+t)dx, \qquad \forall\,\eta\in C([-T,0]).
\]
This shows easily that $T_{i,t}$ is continuous with respect to the topology of $\mathscr C([-T,0])$. This in turn implies that $\Uc$ is continuous with respect to the topology of $\mathscr C([-T,0])$. Therefore, $\Uc$ admits a unique extension $u\colon\mathscr C([-T,0])\rightarrow\R$, which is given by
\[
u(t,\eta) \ = \ \Psi\bigg(t,\int_{[-t,0]}\varphi_1(x+t)d^-\eta(x),\ldots,\int_{[-t,0]}\varphi_N(x+t)d^-\eta(x)\bigg),
\]
for all $(t,\eta)\in[0,T]\times\mathscr C([-T,0])$. We also define the map $\tilde u\colon[0,T]\times\mathscr C([-T,0[)\times\R\rightarrow\R$ as in \eqref{E:tildeu}:
\[
\tilde u(t,\gamma,a) \ = \ u(t,\gamma1_{[-T,0[}+a1_{\{0\}}) \ = \ \Psi\bigg(t,\ldots,a\varphi_i(t) - \int_{-t}^0 \gamma(x)\dot{\varphi}_i(x+t)dx,\ldots\bigg),
\]
for all $(t,\gamma,a)\in[0,T]\times\mathscr C([-T,0[)\times\R$. Let us evaluate the time derivative $\partial_t\Uc(t,\eta)$, for a given $(t,\eta)\in[0,T[\times C([-T,0])$:
\begin{align*}
\partial_t\Uc(t,\eta) \ &= \ \partial_t\Psi\bigg(t,\int_{[-t,0]}\varphi_1(x+t)d^-\eta(x),\ldots,\int_{[-t,0]}\varphi_N(x+t)d^-\eta(x)\bigg) \\
&\quad \ + \sum_{i=1}^N D_{x_i}\Psi\bigg(t,\ldots,\int_{[-t,0]}\varphi_i(x+t)d^-\eta(x),\ldots\bigg)\partial_t\bigg(\int_{[-t,0]}\varphi_i(x+t)d^-\eta(x)\bigg).
\end{align*}
Notice that
\begin{align*}
\partial_t\bigg(\int_{[-t,0]}\varphi_i(x+t)d^-\eta(x)\bigg) \ &= \ \partial_t\bigg(\eta(0)\varphi(t) - \int_{-t}^0\eta(x)\dot\varphi_i(x+t)dx\bigg) \\
&= \ \eta(0)\dot\varphi(t) - \eta(-t)\dot\varphi_i(0^+) - \int_{-t}^0 \eta(x)\ddot\varphi_i(x+t) dx.
\end{align*}
Let us proceed with the horizontal derivative.  We have
\begin{align*}
&D^H\Uc(t,\eta) \ = \ D^H u(t,\eta) \ = \ D^H\tilde u(t,\eta_{|[-T,0[},\eta(0)) \\
&= \ \lim_{\eps\rightarrow0^+} \frac{\tilde u(t,\eta_{|[-T,0[}(\cdot),\eta(0)) - \tilde u(t,\eta_{|[-T,0[}(\cdot-\eps),\eta(0))}{\eps} \\
&= \ \lim_{\eps\rightarrow0^+} \bigg(\frac{1}{\eps}\Psi\left(t,\ldots,\eta(0)\varphi_i(t) - \int_{-t}^0 \eta(x)\dot{\varphi}_i(x+t)dx,\ldots\right) \\
&\quad \ - \frac{1}{\eps}\Psi\left(t,\ldots,\eta(0)\varphi_i(t) - \int_{-t}^0 \eta(x-\eps)\dot{\varphi}_i(x+t)dx,\ldots\right)\bigg).
\end{align*}
From the fundamental theorem of calculus, we obtain
\begin{align*}
&\frac{1}{\eps}\Psi\left(t,\ldots,\eta(0)\varphi_i(t) - \int_{-t}^0 \eta(x)\dot{\varphi}_i(x+t)dx,\ldots\right) \\
&- \frac{1}{\eps}\Psi\left(t,\ldots,\eta(0)\varphi_i(t) - \int_{-t}^0 \eta(x-\eps)\dot{\varphi}_i(x+t)dx,\ldots\right) \\
&= \ \frac{1}{\eps}\int_0^\eps \sum_{i=1}^N D_{x_i}\Psi\left(t,\ldots,\eta(0)\varphi_i(t) - \int_{-t}^0 \eta(x-y)\dot{\varphi}_i(x+t)dx,\ldots\right)\partial_y\bigg(\eta(0)\varphi_i(t) \\
&\quad \ - \int_{-t}^0 \eta(x-y)\dot{\varphi}_i(x+t)dx\bigg) dy.
\end{align*}
Notice that
\begin{align*}
&\partial_y\bigg(\eta(0)\varphi_i(t) - \int_{-t}^0 \eta(x-y)\dot{\varphi}_i(x+t)dx\bigg) \ = \ -\partial_y\bigg(\int_{-t-y}^{-y}\eta(x)\dot{\varphi}_i(x+y+t)dx\bigg) \\
&= \ - \bigg(\eta(-y)\dot{\varphi}_i(t) - \eta(-t-y)\dot{\varphi}_i(0^+) + \int_{-t-y}^{-y}\eta(x)\ddot{\varphi}_i(x+y+t)dx\bigg).
\end{align*}
Therefore
\begin{align*}
&D^H\Uc(t,\eta) \\
&= \ -\lim_{\eps\rightarrow0^+} \frac{1}{\eps}\int_0^\eps \sum_{i=1}^N D_{x_i}\Psi\bigg(t,\ldots,\eta(0)\varphi_i(t) - \int_{-t}^0 \eta(x-y)\dot{\varphi}_i(x+t)dx,\ldots\bigg)\bigg(\eta(-y)\dot{\varphi}_i(t) \\
&\quad \ - \eta(-t-y)\dot{\varphi}_i(0^+) + \int_{-t-y}^{-y}\eta(x)\ddot{\varphi}_i(x+y+t)dx\bigg) dy \\
&=\ - \sum_{i=1}^N D_{x_i}\Psi\bigg(t,\ldots,\eta(0)\varphi_i(t) - \int_{-t}^0 \eta(x)\dot{\varphi}_i(x+t)dx,\ldots\bigg)\bigg(\eta(0)\dot\varphi(t) - \eta(-t)\dot\varphi_i(0^+) \\
&\quad \ - \int_{-t}^0 \eta(x)\ddot\varphi_i(x+t) dx\bigg).
\end{align*}
Finally, concerning the vertical derivative we have
\begin{align*}
D^V\Uc(t,\eta) \ = \ D^Vu(t,\eta) \ &= \ \partial_a\tilde u(t,\eta1_{[-T,0[}+\eta(0)1_{\{0\}}) \\
&= \ \sum_{i=1}^N D_{x_i}\Psi\bigg(t,\int_{[-t,0]}\varphi_1(x+t)d^-\eta(x),\ldots\bigg)\varphi_i(t)
\end{align*}
and
\begin{align*}
D^{VV}\Uc(t,\eta) \ = \ D^{VV}u(t,\eta) \ &= \ \partial_{aa}^2\tilde u(t,\eta1_{[-T,0[}+\eta(0)1_{\{0\}}) \\
&= \ \sum_{i,j=1}^N D_{x_ix_j}^2\Psi\bigg(t,\int_{[-t,0]}\varphi_1(x+t)d^-\eta(x),\ldots\bigg)\varphi_i(t)\varphi_j(t).
\end{align*}
From the regularity of $\Psi$ it follows that $\Uc\in C^{1,2}(([0,T]\times\textup{past})\times\textup{present}))$. Moreover, since $\Psi$ satisfies the Cauchy problem \eqref{E:HeatEquationPsi}, we conclude that $\partial_t\Uc(t,\eta) + D^H\Uc(t,\eta) + \frac{1}{2}D^{VV}\Uc(t,\eta)=0$, for all $(t,\eta)\in[0,T[\times C([-T,0])$, therefore $\Uc$ is a classical solution to the path-dependent heat equation \eqref{E:KolmEq}.
\ep

\subsection{Towards a weaker notion of solution: a significant hedging example}
\label{SubS:HedgingExample}

In the present subsection, we consider the path-dependent nonlinear Kolmogorov equation \eqref{E:KolmEq} in the case $F\equiv0$. This situation is particularly interesting, since it arises, for example, in hedging problems of path-dependent contingent claims. More precisely, consider a real continuous finite quadratic variation process $X$ on $(\Omega,\Fc,\P)$ and  denote $\X$ the window process associated to $X$. Let us assume that $[X]_t = t$, for any $t\in[0,T]$. the hedging problem that we have in mind is the following: given a contingent claim's payoff $G(\mathbb X_T)$, is it possible to have
\begin{equation}
\label{E:ReplForm}
G(\mathbb X_T) \ = \ G_0 + \int_0^T Z_t \, d^- X_t,
\end{equation}
for some $G_0\in\R$ and some $\F$-adapted process $Z = (Z_t)_{t\in[0,T]}$ such that $Z_t=v(t,\mathbb X_t)$, with $v\colon[0,T]\times C([-T,0])\rightarrow\R$? When $X$ is a Brownian motion $W$ and $\int_0^T|Z_t|^2dt<\infty$, $\P$-a.s., the previous forward integral is an It\^o integral. If $G$ is regular enough and it is cylindrical in the sense of \eqref{E:G_Cylindrical}, we know from Theorem \ref{T:ExistenceClassical} that there exists a unique classical solution $\Uc\colon[0,T]\times C([-T,0])\rightarrow\R$ to equation \eqref{E:KolmEq}.\\
Then, we see from It\^o's formula \eqref{E:ItoTime} that $\Uc$ satisfies, $\P$-a.s.,
\begin{equation}
\label{E:Hedging}
\Uc(t,\mathbb X_t) \ = \ \Uc(0,\mathbb X_0) + \int_0^t D^V\Uc(s,\mathbb X_s)\,d^- X_s, \qquad 0 \leq t \leq T.
\end{equation}
In particular, \eqref{E:ReplForm} holds with $Z_t=D^V\Uc(t,\mathbb X_t)$,
 for any $t\in[0,T]$, $G_0 = \Uc(0,\mathbb X_t)$.\\
 However, a significant hedging example is
the \emph{lookback-type payoff}
\[
G(\eta) \ = \ \sup_{x\in[-T,0]}\eta(x), \qquad \forall\,\eta\in C([-T,0]).
\]
We look again for $\Uc\colon[0,T]\times C([-T,0])\rightarrow\R$ which verifies \eqref{E:Hedging}, at least for $X$ being a Brownian motion $W$. Since $\Uc(t,\mathbb W_t)$ has to be a martingale, a candidate for $\Uc$ it is $\Uc(t,\eta)=\E[G(\mathbb W_T^{t,\eta})]$, for all $(t,\eta)\in[0,T]\times C([-T,0])$. However, this latter $\Uc$ can  be shown not to be regular to be a classical solution to equation \eqref{E:KolmEq}, even if it is ``virtually'' a solution to the path-dependent nonlinear Kolmogorov equation \eqref{E:KolmEq}. This will lead us to introduce a weaker notion of solution to equation \eqref{E:KolmEq}. To characterize the map $\Uc$, we notice that it admits the probabilistic representation formula, for all $(t,\eta)\in[0,T]\times C([-T,0])$,
\begin{align*}
\Uc(t,\eta) \ &= \ \E\big[G(\mathbb W_T^{t,\eta})\big] \ = \E\Big[\sup_{-T\leq x\leq0}\mathbb W_T^{t,\eta}(x)\Big] \\
&= \ \E\Big[\Big(\sup_{-t\leq x\leq0}\eta(x)\Big)\vee\Big(\sup_{t\leq x\leq T}\big(W_x-W_t+\eta(0)\big)\Big)\Big] \ = \ f\Big(t,\sup_{-t\leq x\leq0}\eta(x),\eta(0)\Big),
\end{align*}
where the function $f\colon[0,T]\times\R\times\R\rightarrow\R$ is given by
\begin{equation}
\label{E:f}
f(t,m,x) \ = \ \E\big[m \vee (S_{T-t} + x)\big], \qquad \forall\,(t,m,x)\in[0,T]\times\R\times\R,
\end{equation}
with $S_t = \sup_{0 \leq s \leq t} W_s$, for all $t\in[0,T]$. Recalling Remark \ref{R:Density}, it follows from the presence of $\sup_{-t\leq x\leq0}\eta(x)$ among the arguments of $f$, that $\Uc$ is not continuous with respect to the topology of $\mathscr C([-T,0])$, therefore it can not be a classical solution to equation \eqref{E:KolmEq}.  However, we notice that $\sup_{-t\leq x\leq0}\eta(x)$ is Lipschitz on $(C([-T,0]),\|\cdot\|_\infty)$, therefore it will follow from Theorem \ref{T:Exist} that $\Uc$ is a strong-viscosity solution to equation \eqref{E:KolmEq} in the sense of Definition \ref{D:Strong}. Nevertheless, in this particular case, even if $\Uc$ is not a classical solution, we shall prove that it is associated to the classical solution of a certain finite dimensional PDE. To this end, we begin computing an explicit form for $f$, for which it is useful to recall the following standard result.

\begin{Lemma}[Reflection principle]
\label{L:LawS}
For every $a>0$ and $t>0$,
\[
\P(S_t \geq a) \ = \ \P(|B_t| \geq a).
\]
In particular, for each $t$, the random variables $S_t$ and $|B_t|$ have the same law, whose density is given by:
\[
\varphi_t(z) \ = \ \sqrt{\frac{2}{\pi t}} e^{-\frac{z^2}{2t}}1_{[0,\infty)}(z), \qquad \forall\,z\in\R.
\]
\end{Lemma}
\textbf{Proof.}
See Proposition 3.7, Chapter III, in \cite{revuzyor99}.
\ep

\vspace{3mm}

\noindent From Lemma \ref{L:LawS} it follows that, for all $(t,m,x)\in[0,T[\times\R\times\R$,
\[
f(t,m,x) \ = \ \int_0^\infty m \vee (z + x) \, \varphi_{T-t}(z) dz \ = \ \int_0^\infty m \vee (z + x) \frac{2}{\sqrt{T-t}} \varphi\Big(\frac{z}{\sqrt{T-t}}\Big) dz,
\]
where $\varphi(z) = \exp(z^2/2)/\sqrt{2\pi}$, $z\in\R$, is the standard Gaussian density.

\begin{Lemma}
\label{L:f}
The function $f$ defined in \eqref{E:f} is given by, for all $(t,m,x)\in[0,T[\times\R\times\R$,
\[
f(t,m,x) \ = \ 2m \Big(\Phi\Big(\frac{m-x}{\sqrt{T-t}}\Big) - \frac{1}{2}\Big) + 2x \Big(1 - \Phi\Big(\frac{m-x}{\sqrt{T-t}}\Big)\Big) + \sqrt{\frac{2(T-t)}{\pi}} e^{-\frac{(m-x)^2}{2(T-t)}},
\]
for $x \leq m$, and
\[
f(t,x,m) \ = \ x + \sqrt{\frac{2(T-t)}{\pi}},
\]
for $x > m$, where $\Phi(y) = \int_{-\infty}^y \varphi(z)dz$, $y\in\R$, is the standard Gaussian cumulative distribution function.
\end{Lemma}
\textbf{Proof.}
\emph{First case: $x \leq m$.} We have
\begin{equation}
\label{E:f1}
f(t,m,x) \ = \ \int_0^{m-x} m \frac{2}{\sqrt{T-t}} \varphi\Big(\frac{z}{\sqrt{T-t}}\Big) dz + \int_{m-x}^\infty (z+x) \frac{2}{\sqrt{T-t}} \varphi\Big(\frac{z}{\sqrt{T-t}}\Big) dz.
\end{equation}
The first integral on the right-hand side of \eqref{E:f1} becomes
\[
\int_0^{m-x} m \frac{2}{\sqrt{T-t}} \varphi\Big(\frac{z}{\sqrt{T-t}}\Big) dz \ = \ 2m \int_0^{\frac{m-x}{\sqrt{T-t}}} \varphi(z) dz \ = \ 2m \Big(\Phi\Big(\frac{m-x}{\sqrt{T-t}}\Big) - \frac{1}{2}\Big),
\]
where $\Phi(y) = \int_{-\infty}^y \varphi(z)dz$, $y\in\R$, is the standard Gaussian cumulative distribution function. Concerning the second integral in \eqref{E:f1}, we have
\begin{align*}
\int_{m-x}^\infty (z+x) \frac{2}{\sqrt{T-t}} \varphi\Big(\frac{z}{\sqrt{T-t}}\Big) dz \ &= \ 2\sqrt{T-t}\int_{\frac{m-x}{\sqrt{T-t}}}^\infty z\varphi(z) dz + 2x \int_{\frac{m-x}{\sqrt{T-t}}}^\infty \varphi(z) dz \\
&= \ \sqrt{\frac{2(T-t)}{\pi}} e^{-\frac{(m-x)^2}{2(T-t)}} + 2x \Big(1 - \Phi\Big(\frac{m-x}{\sqrt{T-t}}\Big)\Big).
\end{align*}
\emph{Second case: $x>m$.} We have
\begin{align*}
f(t,m,x) \ &= \ \int_0^\infty (z + x) \frac{2}{\sqrt{T-t}} \varphi\Big(\frac{z}{\sqrt{T-t}}\Big) dz \\
&= \ 2\sqrt{T-t}\int_0^\infty z\varphi(z)dz + 2x\int_0^\infty \varphi(z)dz \ = \ \sqrt{\frac{2(T-t)}{\pi}} + x.
\end{align*}
\ep

\vspace{3mm}

\noindent We also have the following regularity result regarding the function $f$.

\begin{Lemma}
The function $f$ defined in \eqref{E:f} is continuous on $[0,T]\times\R\times\R$, moreover it is once $($resp. twice$)$ continuously differentiable in $(t,m)$ $($resp. in $x$$)$ on $[0,T[\times\overline Q$, where $\overline Q$ is the closure of the set $Q := \{(m,x)\in\R\times\R\colon m > x \}$. In addition, the following It\^o's formula holds:
\begin{align}
\label{E:Ito_f}
f(t,S_t,B_t) \ &= \ f(0,0,0) + \int_0^t \Big(\partial_t f(s,S_s,B_s) + \frac{1}{2}\partial_{xx}^2 f(s,S_s,B_s)\Big) ds \\
&\quad \ + \int_0^t \partial_m f(s,S_s,B_s) dS_s + \int_0^t \partial_x f(s,S_s,B_s) dB_s, \qquad 0 \leq t \leq T,\,\P\text{-a.s.} \notag
\end{align}
\end{Lemma}
\textbf{Proof.}
The regularity properties of $f$ are deduced from its explicit form derived in Lemma~\ref{L:f}, after straightforward calculations. Concerning It\^o's formula \eqref{E:Ito_f}, the proof can be done along the same lines as the standard It\^o's formula. We simply notice that, in the present case, only the restriction of $f$ to $\overline Q$ is smooth. However, the process $((S_t,B_t))_t$ is $\overline Q$-valued. It is well-known that if $\overline Q$ would be an open set, then It\^o's formula would hold. In our case, $\overline Q$ is the closure of its interior $Q$. This latter property is enough for the validity of It\^o's formula. In particular, the basic tools for the proof of It\^o's formula are the following Taylor's expansions for the function $f$:
\begin{align*}
f(t',m,x) \ &= \ f(t,m,x) + \partial_t f(t,m,x) (t'-t) \\
&\quad \ + \int_0^1 \partial_t f(t + \lambda(t'-t),m,x)(t' - t) d\lambda, \\
f(t,m',x) \ &= \ f(t,m,x) + \partial_m f(t,m,x) (m'-m) \\
&\quad \ + \int_0^1 \partial_m f(t,m + \lambda(m'-m),x)(m' - m) d\lambda, \\
f(t,m,x') \ &= \ f(t,m,x) + \partial_x f(t,m,x)(x'-x) + \frac{1}{2}\partial_{xx}^2 f(t,m,x)(x'-x)^2 \\
&\quad \ + \int_0^1(1-\lambda)\big(\partial_{xx}^2 f(t,m,x+\lambda(x'-x)) - \partial_{xx}^2 f(t,m,x)\big)(x'-x)^2d\lambda,
\end{align*}
for all $(t,m,x)\in[0,T]\times\overline Q$. To prove the above Taylor's formulae, note that they hold on the open set $Q$, using the regularity of $f$. Then, we can extend them to the closure of $Q$, since $f$ and its derivatives are continuous on $\overline Q$. Consequently, It\^o's formula can be proved in the usual way.
\ep

\vspace{3mm}

Even though, as already observed, $\Uc$ does not belong to $C^{1,2}(([0,T[\times\text{past})\times\text{present})\cap C([0,T]\times C([-T,0]))$, so that it can not be a classical solution to equation \eqref{E:KolmEq}, the function $f$ is a solution to a certain Cauchy problem, as stated in the following proposition.

\begin{Proposition}
The function $f$ defined in \eqref{E:f} solves the backward heat equation:
\[
\begin{cases}
\partial_t f(t,m,x) + \frac{1}{2}\partial_{xx}^2 f(t,m,x) \ = \ 0, \qquad\;\, &\forall\,(t,m,x)\in[0,T[\times\overline Q, \\
f(T,m,x) \ = \ m, &\forall\,(m,x)\in\overline Q.
\end{cases}
\]
\end{Proposition}
\textbf{Proof.}
We provide two distinct proofs.\\
\emph{Direct proof.} Since we know the explicit expression of $f$, we can derive the form of $\partial_t f$ and $\partial_{xx}^2 f$ by direct calculations:
\[
\partial_t f(t,m,x) \ = \ - \frac{1}{\sqrt{T-t}} \varphi\Big(\frac{m-x}{\sqrt{T-t}}\Big), \qquad \partial_{xx}^2 f(t,m,x) \ = \ \frac{2}{\sqrt{T-t}} \varphi\Big(\frac{m-x}{\sqrt{T-t}}\Big),
\]
for all $(t,m,x)\in[0,T[\times\overline Q$, from which the thesis follows.\\
\emph{Probabilistic proof.} By definition, the process $(f(t,S_t,B_t))_{t\in[0,T]}$ is given by:
\begin{align*}
f(t,S_t,B_t) \ = \ \E\big[S_T\big|\Fc_t\big],
\end{align*}
so that it is a uniformly integrable $\F$-martingale. Then, it follows from It\^o's formula \eqref{E:Ito_f} that
\[
\int_0^t \Big(\partial_t f(s,S_s,B_s) + \frac{1}{2}\partial_{xx}^2 f(s,S_s,B_s)\Big) ds + \int_0^t \partial_m f(s,S_s,B_s) dS_s \ = \ 0,
\]
for all $0 \leq t \leq T$, $\P$-almost surely. As a consequence, the thesis follows if we prove that
\begin{equation}
\label{E:LocalTime}
\int_0^t \partial_m f(s,S_s,B_s) dS_s \ = \ 0.
\end{equation}
By direct calculation, we have
\[
\partial_m f(t,m,x) \ = \ 2\Phi\Big(\frac{m-x}{\sqrt{T-t}}\Big) - 1, \qquad \forall(t,m,x)\in[0,T[\times\overline Q.
\]
Therefore, \eqref{E:LocalTime} becomes
\begin{equation}
\label{E:LocalTime2}
\int_0^t \bigg(2\Phi\bigg(\frac{S_s-B_s}{\sqrt{T-s}}\bigg) - 1\bigg) dS_s \ = \ 0.
\end{equation}
Now we observe that the local time of $S_s-B_s$ is equal to $2S_s$, see Exercise 2.14 in \cite{revuzyor99}. It follows that the measure $dS_s$ is carried by $\{s\colon S_s-B_s=0\}$. This in turn implies the validity of \eqref{E:LocalTime2}, since the integrand in \eqref{E:LocalTime2} is zero on the set $\{s\colon S_s-B_s=0\}$.
\ep

\subsection{Path-dependent Kolmogorov equation: strong-viscosity solutions}
\label{SubS:StringViscositySolutions}

Motivated by previous subsection, we now introduce a notion of weak solution for the path-dependent nonlinear Kolmogorov equation \eqref{E:KolmEq}, which we refer to as \emph{strong-viscosity solution}. Firstly, we need the following definition.

\begin{Definition}
Let $\mathscr F$ be a collection of $\R^d$-valued functions on $[0,T]\times X$, where $(X,\|\cdot\|)$ is a normed space. We say that $\mathscr F$ is \textbf{locally equicontinuous} if to any $R,\eps>0$ corresponds a $\delta$ such that $|f(t,x)-f(s,y)|<\eps$ for every $f\in\mathscr F$ and for all pair of points $(t,x),(s,y)$ with $|t-s|,\|x-y\|<\delta$ and $\|x\|,\|y\|<R$.
\end{Definition}

\begin{Definition}
\label{D:Strong}
A function $\Uc\colon[0,T]\times C([-T,0])\rightarrow\R$ is called  \textbf{strong-viscosity solution} to the path-dependent nonlinear Kolmogorov equation \eqref{E:KolmEq} if there exists a sequence $(\Uc_n,G_n,F_n)_n$ satisfying:
\begin{enumerate}
\item[\textup{(i)}] $\Uc_n\colon[0,T]\times C([-T,0])\rightarrow\R$, $G_n\colon C([-T,0])\rightarrow\R$, and $F_n\colon[0,T]\times C([-T,0])\times\R\times\R\rightarrow\R$ are locally equicontinuous functions such that, for some positive constants $C$ and $m$, independent of $n$,
\begin{align*}
|F_n(t,\eta,y,z) - F_n(t,\eta,y',z')| \ &\leq \ C(|y-y'| + |z-z'|), \\
|\Uc_n(t,\eta)| + |G_n(\eta)| + |F_n(t,\eta,0,0)| \ &\leq \ C\big(1 + \|\eta\|_\infty^m\big),
\end{align*}
for all $(t,\eta)\in[0,T]\times C([-T,0])$, $y,y'\in\R$, and $z,z'\in\R$.
\item[\textup{(ii)}] $\Uc_n$ is a classical solution to
\[
\begin{cases}
\partial_t\Uc_n + D^H\Uc_n + \frac{1}{2}D^{VV}\Uc_n \ = \ F_n(t,\eta,\Uc_n,D^V\Uc_n), \;\;\; &\forall\,(t,\eta)\in[0,T)\times C([-T,0]), \\
\Uc_n(T,\eta) \ = \ G_n(\eta), &\forall\,\eta\in C([-T,0]).
\end{cases}
\]
\item[\textup{(iii)}] $(\Uc_n(t,\eta),G_n(\eta),F_n(t,\eta,y,z))\rightarrow(\Uc(t,\eta),G(\eta),F(t,\eta,y,z))$, as $n$ tends to infinity, for any $(t,\eta,y,z)\in[0,T]\times C([-T,0])\times\R\times\R$.
\end{enumerate}
\end{Definition}

\noindent The following uniqueness result for strong-viscosity solution holds.

\begin{Theorem}
Let $\Uc\colon[0,T]\times C([-T,0])\rightarrow\R$ be a strong-viscosity solution to the path-dependent nonlinear Kolmogorov equation \eqref{E:KolmEq}. Then, we have
\[
\Uc(t,\eta) \ = \ Y_t^{t,\eta}, \qquad \forall\,(t,\eta)\in[0,T]\times C([-T,0]),
\]
where $(Y_s^{t,\eta},Z_s^{t,\eta})_{s\in[t,T]}\in\S^2(t,T)\times\H^2(t,T)$, with $Y_s^{t,\eta}=\Uc(s,\mathbb W_s^{t,\eta})$, solves the backward stochastic differential equation, $\P$-a.s.,
\[
Y_s^{t,\eta} \ = \ G(\mathbb W_T^{t,\eta}) + \int_s^T F(r,\mathbb W_r^{t,\eta},Y_r^{t,\eta},Z_r^{t,\eta}) dr - \int_s^T Z_r^{t,\eta} dW_r, \qquad t \leq s \leq T.
\]
In particular, there exists at most one strong-viscosity solution to the path-dependent nonlinear Kolmogorov equation \eqref{E:KolmEq}.
\end{Theorem}
\textbf{Proof.}
Consider a sequence $(\Uc_n,G_n,F_n)_n$ satisfying conditions (i)-(iii) of Definition \ref{D:Strong}. For every $n\in\N$ and any $(t,\eta)\in[0,T]\times C([-T,0])$, we know from Theorem \ref{T:UniqClassical} that $(Y_s^{n,t,\eta},Z_s^{n,t,\eta})_{s\in[t,T]} = (\Uc_n(s,\mathbb W_s^{t,\eta}),D^V\Uc_n(s,\mathbb W_s^{t,\eta}))_{s\in[t,T]}\in\S^2(t,T)\times\H^2(t,T)$ is the solution to the backward stochastic differential equation, $\P$-a.s.,
\[
Y_s^{n,t,\eta} \ = \ G_n(\mathbb W_T^{t,\eta}) + \int_s^T F_n(r,\mathbb W_r^{t,\eta},Y_r^{n,t,\eta},Z_r^{n,t,\eta}) dr - \int_s^T Z_r^{n,t,\eta} dW_r, \qquad t \leq s \leq T.
\]
From the polynomial growth condition of $(\Uc_n)_n$ and estimate \eqref{E:EstimateSupW}, we see that
\[
\sup_n\|Y^{n,t,\eta}\|_{\S^p(t,T)} \ < \ \infty, \qquad \text{for any }p\geq1.
\]
This implies, using Proposition \ref{P:EstimateBSDEAppendix} and the polynomial growth condition of $(F_n)_n$, that
\[
\sup_n\|Z^{n,t,\eta}\|_{\H^2(t,T)} \ < \ \infty.
\]
Let $Y_s^{t,\eta}=\Uc(s,\mathbb W_s^{t,\eta})$, for any $s\in[t,T]$. Then, we see that all the hypotheses of Proposition \ref{P:LimitThmBSDE} follow by assumptions and estimate \eqref{E:EstimateSupW} (notice that, in this case, $K^n\equiv0$ for any $n$, therefore the proof of Proposition \ref{P:LimitThmBSDE} simplifies drastically), so the thesis follows.
\ep

\vspace{3mm}

We now prove an existence result for strong-viscosity solutions to the path-dependent heat equation, namely to equation \eqref{E:KolmEq} in the case $F\equiv0$. To this end, we need the following stability result for strong-viscosity solutions.

\begin{Lemma}
\label{L:Stability}
Consider $(\Uc_n,G_n,F_n)_n$ and $(\Uc_{n,k},G_{n,k},F_{n,k})_{n,k}$ satisfying:
\begin{enumerate}
\item[\textup{(i)}] $\Uc_{n,k}\colon[0,T]\times C([-T,0])\rightarrow\R$, $G_{n,k}\colon C([-T,0])\rightarrow\R$, and $F_{n,k}\colon[0,T]\times C([-T,0])\times\R\times\R\rightarrow\R$ are locally equicontinuous functions such that, for some positive constants $C$ and $m$, independent of $n$ and $k$,
\begin{align*}
|F_{n,k}(t,\eta,y,z) - F_{n,k}(t,\eta,y',z')| \ &\leq \ C(|y-y'| + |z-z'|), \\
|\Uc_{n,k}(t,\eta)| + |G_{n,k}(\eta)| + |F_{n,k}(t,\eta,0,0)| \ &\leq \ C\big(1 + \|\eta\|_\infty^m\big),
\end{align*}
for all $(t,\eta)\in[0,T]\times C([-T,0])$, $y,y'\in\R$, and $z,z'\in\R$.
\item[\textup{(ii)}] $\Uc_{n,k}$ is a classical solution to
\[
\begin{cases}
\partial_t\Uc_{n,k} + D^H\Uc_{n,k} + \frac{1}{2}D^{VV}\Uc_{n,k} = F_{n,k}(t,\eta,\Uc_{n,k},D^V\Uc_{n,k}), &\forall\,(t,\eta)\in[0,T)\times C([-T,0]), \\
\Uc_{n,k}(T,\eta) \ = \ G_{n,k}(\eta), &\forall\,\eta\in C([-T,0]).
\end{cases}
\]
\item[\textup{(iii)}] $(\Uc_{n,k}(t,\eta),G_{n,k}(\eta),F_{n,k}(t,\eta,y,z))\rightarrow(\Uc_n(t,\eta),G_n(\eta),F_n(t,\eta,y,z))$, as $k$ tends to infinity, for any $n\in\N$ and $(t,\eta,y,z)\in[0,T]\times C([-T,0])\times\R\times\R$.
\end{enumerate}
If for every $(t,\eta,y,z)\in[0,T]\times C([-T,0])\times\R\times\R$, $(\Uc_n(t,\eta),G_n(\eta),F_n(t,\eta,y,z))_n$ converges, we define
\[
(\Uc(t,\eta),G(\eta),F(t,\eta,y,z)) \ := \ \lim_{n\rightarrow\infty}(\Uc_n(t,\eta),G_n(\eta),F_n(t,\eta,y,z)).
\]
Then, there exists a subsequence $(\Uc_{n,k(n)},G_{n,k(n)},F_{n,k(n)})_n$ which converges pointwisely to $(\Uc,G,F)$, as $n$ tends to infinity, so that $\Uc$ is a strong solution to
\[
\begin{cases}
\partial_t\Uc + D^H\Uc + \frac{1}{2}D^{VV}\Uc \ = \ F(t,\eta,\Uc,D^V\Uc), \;\;\; &\forall\,(t,\eta)\in[0,T)\times C([-T,0]), \\
\Uc(T,\eta) \ = \ G(\eta), &\forall\,\eta\in C([-T,0]).
\end{cases}
\]
\end{Lemma}
\textbf{Proof.}
The thesis follows from Lemma \ref{L:StabilityApp}.
\ep

\begin{Theorem}
\label{T:Exist}
Let $F\equiv0$ and $G\colon C([-T,0])\rightarrow\R$ be locally uniformly continuous and satisfying the polynomial growth condition
\[
|G(\eta)| \ \leq \ C(1 + \|\eta\|_\infty^m), \qquad \forall\,\eta\in C([-T,0]),
\]
for some positive constants $C$ and $m$. Then, there exists a unique strong solution $\Uc$ to the path-dependent heat equation \eqref{E:KolmEq}, which is given by
\[
\Uc(t,\eta) \ = \ \E\big[G(\mathbb W_T^{t,\eta})\big], \qquad \forall\,(t,\eta)\in[0,T]\times C([-T,0]).
\]
\end{Theorem}
\textbf{Proof.}
Let $(e_i)_{i\geq0}$ be the orthonormal basis of $L^2([-T,0])$ composed by $C^\infty([-T,0])$, periodic, and uniformly bounded functions:
\[
e_0=\frac{1}{\sqrt{T}}, \quad e_{2i-1}=\sqrt{\frac{2}{T}}\sin\bigg(\frac{2\pi}{T}(x+T)i\bigg), \quad e_{2i}=\sqrt{\frac{2}{T}}\cos\bigg(\frac{2\pi}{T}(x+T)i\bigg), \quad i\in\N\backslash\{0\}.
\]
Let us define the linear operator $\Lambda\colon C([-T,0])\rightarrow C([-T,0])$ by
\[
(\Lambda\eta)(x) \ = \ \frac{\eta(0)-\eta(-T)}{T}x, \qquad x\in[-T,0],\,\eta\in C([-T,0]).
\]
Notice that $(\eta-\Lambda\eta)(-T) = (\eta-\Lambda\eta)(0)$, therefore $\eta-\Lambda\eta$ can be extended to the entire real line in a periodic way with period $T$, so that we can expand it in Fourier series. In particular, for each $n\in\N$ and $\eta\in C([-T,0])$, consider the Fourier partial sum
\begin{equation}
\label{E:s_n}
s_n(\eta-\Lambda\eta) \ = \ \sum_{i=0}^n (\eta_i-(\Lambda\eta)_i) e_i, \qquad \forall\,\eta\in C([-T,0]),
\end{equation}
where (denoting $\tilde e_i(x) = \int_{-T}^x e_i(y) dy$, for any $x\in[-T,0]$), by Proposition \ref{P:BVI},
\begin{align}
\label{E:eta_i}
\eta_i \ = \ \int_{-T}^0 \eta(x)e_i(x) dx \ &= \ \eta(0)\tilde e _i(0) - \int_{[-T,0]} \tilde e_i(x)d^-\eta(x) \notag \\
&= \ \int_{[-T,0]} (\tilde e_i(0) - \tilde e_i(x)) d^- \eta(x),
\end{align}
since $\eta(0) = \int_{[-T,0]}d^-\eta(x)$. Moreover we have
\begin{align}
\label{E:Lambda_eta_i}
(\Lambda\eta)_i \ &= \ \int_{-T}^0 (\Lambda\eta)(x)e_i(x) dx \ = \ \frac{1}{T} \int_{-T}^0 xe_i(x) dx \bigg(\int_{[-T,0]} d^-\eta(x) - \eta(-T)\bigg).
\end{align}
Define
\[
\sigma_n \ = \ \frac{s_0 + s_1 + \cdots + s_n}{n+1}.
\]
Then, by \eqref{E:s_n},
\[
\sigma_n(\eta-\Lambda\eta) \ = \ \sum_{i=0}^n \frac{n+1-i}{n+1} (\eta_i-(\Lambda\eta)_i) e_i, \qquad \forall\,\eta\in C([-T,0]).
\]
We know from Fej\'er's theorem on Fourier series (see, e.g., Theorem 3.4, Chapter III, in \cite{zygmund02}) that, for any $\eta\in C([-T,0])$, $\sigma_n(\eta-\Lambda\eta)\rightarrow\eta-\Lambda\eta$ uniformly on $[-T,0]$, as $n$ tends to infinity, and $\|\sigma_n(\eta-\Lambda\eta)\|_\infty \leq \|\eta-\Lambda\eta\|_\infty$. Let us define the linear operator $T_n\colon C([-T,0])\rightarrow C([-T,0])$ by (denoting $e_{-1}(x) = x$, for any $x\in[-T,0]$)
\begin{align}
\label{E:T_n_eta}
T_n\eta \ = \ \sigma_n(\eta-\Lambda\eta) + \Lambda\eta \ &= \ \sum_{i=0}^n \frac{n+1-i}{n+1} (\eta_i-(\Lambda\eta)_i) e_i + \frac{\eta(0) - \eta(-T)}{T}e_{-1} \notag \\
&= \ \sum_{i=0}^n \frac{n+1-i}{n+1} x_ie_i + x_{-1}e_{-1},
\end{align}
where, using \eqref{E:eta_i} and \eqref{E:Lambda_eta_i},
\begin{align*}
x_{-1} \ &= \int_{[-T,0]} \frac{1}{T} d^-\eta(x) - \frac{1}{T}\eta(-T), \\
x_i \ &= \ \int_{[-T,0]}\bigg (\tilde e_i(0)-\tilde e_i(x)- \frac{1}{T}\int_{-T}^0 xe_i(x)dx\bigg )d^-\eta(x) + \frac{1}{T}\int_{-T}^0 xe_i(x)dx\,\eta(-T),
\end{align*}
for $i=0,\ldots,n$. Then, for any $\eta\in C([-T,0])$, $T_n\eta\rightarrow\eta$ uniformly on $[-T,0]$, as $n$ tends to infinity. Furthermore, there exists a positive constant $M$ such that
\begin{equation}
\label{E:UniformBoundT_n}
\|T_n\eta\|_\infty \ \leq \ M\|\eta\|_\infty, \qquad \forall\,n\in\N,\,\forall\,\eta\in C([-T,0]).
\end{equation}
In particular, the family of linear operators $(T_n)_n$ is equicontinuous. Now, let us define $G_n\colon C([-T,0])\rightarrow\R$ as follows
\[
G_n(\eta) \ = \ G(T_n\eta), \qquad \forall\,\eta\in C([-T,0]).
\]
We see from \eqref{E:UniformBoundT_n} that the family $(G_n)_n$ is locally equicontinuous. Moreover, from the polynomial growth condition of $G$ and \eqref{E:UniformBoundT_n} we have
\[
|G_n(\eta)| \ \leq \ C(1 + \|T_n\eta\|_\infty^m) \ \leq \ C(1 + M^m\|\eta\|_\infty^m), \qquad \forall\,n\in\N,\,\forall\,\eta\in C([-T,0]).
\]
Now, we observe that since $\{e_{-1},e_0,e_1,\ldots,e_n\}$ are linearly independent, then we see from \eqref{E:T_n_eta} that $T_n\eta$ is completely characterized by the coefficients of $e_{-1},e_0,e_1,\ldots,e_n$. Therefore, the function $g_n\colon\R^{n+2}\rightarrow\R$ given by
\[
g_n(x_{-1},\ldots,x_n) \ = \ G_n(\eta) \ = \ G\bigg(\sum_{i=0}^n \frac{n+1-i}{n+1} x_ie_i + x_{-1}e_{-1}\bigg), \qquad \forall\,(x_{-1},\ldots,x_n)\in\R^{n+2},
\]
completely characterizes $G_n$. Since the family $(G_n)_n$ is locally equicontinuous, it follows that the family $(g_n)_n$ is locally equicontinuous, as well. Moreover, fix $\eta\in C([-T,0])$ and consider the corresponding coefficients $x_{-1},\ldots,x_n$ with respect to $\{e_{-1},\ldots,e_n\}$ in the expression \eqref{E:T_n_eta} of $T_n\eta$. Set
\begin{align*}
\varphi_{-1}(x) \ &= \ \frac{1}{T}, \qquad \varphi_i(x) \ = \ \tilde e_i(0)-\tilde e_i(x-T)- \frac{1}{T}\int_{-T}^0 xe_i(x)dx, \qquad x\in[0,T], \\
a_{-1} \ &= \ -\frac{1}{T}, \qquad\;\; a_i \ = \ \frac{1}{T}\int_{-T}^0 xe_i(x)dx.
\end{align*}
Notice that $\varphi_{-1},\ldots,\varphi_n\in C^\infty([0,T])$. Then, we have
\[
G_n(\eta) = g_n\bigg(\int_{[-T,0]}\varphi_{-1}(x+T)d^-\eta(x) + a_{-1}\eta(-T),\ldots,\int_{[-T,0]}\varphi_n(x+T)d^-\eta(x) + a_n\eta(-T)\bigg).
\]
Let $\phi(x) = c\exp(1/((x+T)^2-1))1_{[-T,-T+1[}(x)$, $x\in[-T,0]$, with $c>0$ such that $\int_{-T}^0\phi(x)dx=1$. Define, for any $\eps>0$, $\phi_\eps(x) = \phi(x/\eps)/\eps$, $x\in[-T,0]$. Notice that $\phi_\eps\in C^\infty([-T,0])$ and (denoting $\tilde\phi_\eps(x) = \int_{-T}^x\phi_\eps(-T-y)dy$, for any $x\in[-T,0]$)
\[
\int_{-T}^0\eta(x)\phi_\eps(-T-x)dx \ = \ \eta(0)\tilde\phi_\eps(0) - \int_{[-T,0]}\tilde\phi_\eps(x)d^-\eta(x) \ = \ \int_{[-T,0]}\big(\tilde\phi_\eps(0)-\tilde\phi_\eps(x)\big)d^-\eta(x).
\]
Therefore
\[
\lim_{\eps\rightarrow0^+}\int_{[-T,0]}\big(\tilde\phi_\eps(0)-\tilde\phi_\eps(x)\big)d^-\eta(x) \ = \ \lim_{\eps\rightarrow0^+} \int_{-T}^0\eta(x)\phi_\eps(-T-x)dx \ = \ \eta(-T).
\]
For this reason, we introduce the function $G_{n,\eps}\colon C([-T,0])\rightarrow\R$ given by
\[
G_{n,\eps}(\eta) \ = \ g_n\bigg(\ldots,\int_{[-T,0]}\varphi_i(x+T)d^-\eta(x) + a_i\int_{[-T,0]}\big(\tilde\phi_\eps(0)-\tilde\phi_\eps(x)\big)d^-\eta(x),\ldots\bigg).
\]
Now, for any $n\in\N$, let $(g_{n,k})_{k\in\N}$ be a locally equicontinuous sequence of $C^2(\R^{n+2};\R)$ functions, uniformly polynomially bounded, such that $g_{n,k}$ converges pointwise to $g_n$, as $k$ tends to infinity. Define $G_{n,\eps,k}\colon C([-T,0])\rightarrow\R$ as follows:
\[
G_{n,\eps,k}(\eta) \ = \ g_{n,k}\bigg(\ldots,\int_{[-T,0]}\varphi_i(x+T)d^-\eta(x) + a_i\int_{[-T,0]}\big(\tilde\phi_\eps(0)-\tilde\phi_\eps(x)\big)d^-\eta(x),\ldots\bigg).
\]
Then, we know from Theorem \ref{T:ExistenceClassical} that the function $\Uc_{n,\eps,k}\colon[0,T]\times C([-T,0])\rightarrow\R$ given by
\[
\Uc_{n,\eps,k}(t,\eta) \ = \ \E\big[G_{n,\eps,k}(\mathbb W_T^{t,\eta})\big], \qquad \forall\,(t,\eta)\in[0,T]\times C([-T,0])
\]
is a classical solution to the path-dependent heat equation \eqref{E:KolmEq}. Moreover, the family $(\Uc_{n,\eps,k})_{n,\eps,k}$ is locally equicontinuous and uniformly polynomially bounded. Then, using the stability result Lemma \ref{L:Stability}, it follows that $\Uc$ is a strong-viscosity solution to the path-dependent heat equation \eqref{E:KolmEq}.
\ep

\subsection{Strong-viscosity solutions: finite dimensional case}
\label{SubS:FiniteDimensionalCase}

In the previous subsection we provided the definition of strong-viscosity solution for a particular path-dependent nonlinear Kolmogorov equation, driven by the path-dependent heat operator. The definition of strong-viscosity solution can be adapted to more general semilinear path-dependent PDEs and,  in the present subsection, we want to give an idea on how this definition can be generalized, focusing on the more understandable finite dimensional case. In particular, our aim is to emphasize how uniqueness for strong-viscosity solutions can be proved using probabilistic methods, in contrast with real analysis' tools which characterize comparison theorem for viscosity solutions. We provide two definitions of strong-viscosity solution, presented in order of increasing generality, and we prove the comparison theorem. We conclude discussing the relation with the standard notion of viscosity solution.

\vspace{3mm}

Let $b\colon[0,T]\times\R^d\rightarrow\R$, $\sigma\colon[0,T]\times\R^d\rightarrow\R^{d\times d}$, $f\colon[0,T]\times\R^d\times\R\times\R^d\rightarrow\R$, and $g\colon\R^d\rightarrow\R$ be Borel measurable functions. Consider the nonlinear Kolmogorov equation (we denote $A\trans$ the transpose of a matrix $A\in\R^{d\times d}$)
\begin{equation}
\label{E:KolmEqFinite}
\begin{cases}
- \partial_t u(t,x) - \langle b(t,x),D_x u(t,x)\rangle - \frac{1}{2}\textup{tr}(\sigma\sigma\trans(t,x)D_x^2 u(t,x)) & \\
\hspace{2.8cm}-\, f(t,x,u(t,x),\sigma\trans(t,x)D_x u(t,x)) \ = \ 0, &\forall\,(t,x)\in[0,T)\times\R^d, \\
u(T,x) \ = \ g(x), &\forall\,x\in\R^d.
\end{cases}
\end{equation}

\subsubsection{First definition of strong-viscosity solution}

We begin providing the standard definition of classical solution.

\begin{Definition}
A function $u\colon[0,T]\times\R^d\rightarrow\R$, with $u\in C^{1,2}([0,T[\times\R^d)\cap C([0,T]\times\R^d)$, is called a \textbf{classical solution} to the nonlinear Kolmogorov equation \eqref{E:KolmEqFinite} if $u$ solves \eqref{E:KolmEqFinite}.
\end{Definition}
We have the following uniqueness result for classical solutions.

\begin{Proposition}
\label{P:UniqClassicalFiniteStandard}
Suppose that the functions $b$, $\sigma$, $f$, and $g$, appearing in the nonlinear Kolmogorov equation \eqref{E:KolmEqFinite}, satisfy, for some positive constants $C$ and $m$,
\begin{align*}
|b(t,x)-b(t,x')| + |\sigma(t,x)-\sigma(t,x')| \ &\leq \ C|x-x'|, \\
|f(t,x,y,z)-f(t,x,y',z')| \ &\leq \ C\big(|y-y'| + |z-z'|\big), \\
|b(t,0)| + |\sigma(t,0)| \ &\leq \ C, \\
|f(t,x,0,0)| + |g(x)| \ &\leq \ C\big(1 + |x|^m\big),
\end{align*}
for all $t\in[0,T]$, $x,x'\in\R^d$, $y,y'\in\R$, and $z,z'\in\R^d$. Let $u\colon[0,T]\times\R^d\rightarrow\R$ be a classical solution to the nonlinear Kolmogorov equation \eqref{E:KolmEqFinite}, satisfying the polynomial growth condition
\[
|u(t,x)| \ \leq \ C\big(1 + |x|^m\big), \qquad \forall\,(t,x)\in[0,T]\times\R^d.
\]
Then, we have
\[
u(t,x) \ = \ Y_t^{t,x}, \qquad \forall\,(t,x)\in[0,T]\times\R^d,
\]
where $(Y_s^{t,x},Z_s^{t,x})_{s\in[t,T]}=(u(s,X_s^{t,x}),\sigma\trans(s,X_s^{t,x})D_x u(s,X_s^{t,x}))_{s\in[t,T]}\in\S^2(t,T)\times\H^2(t,T)^d$ is the solution to the backward stochastic differential equation, $\P$-a.s.,
\[
Y_s^{t,x} \ = \ g(X_T^{t,x}) + \int_s^T f(r,X_r^{t,x},Y_r^{t,x},Z_r^{t,x}) dr - \int_s^T Z_r^{t,x} dW_r, \qquad t \leq s \leq T.
\]
In particular, there exists at most one classical solution to the nonlinear Kolmogorov equation \eqref{E:KolmEqFinite}.
\end{Proposition}
\textbf{Proof.}
The proof can be done along the lines of Theorem \ref{T:UniqClassical}.
\ep

\vspace{3mm}

\noindent We can now present our first definition of strong-viscosity solution to equation \eqref{E:KolmEqFinite}.

\begin{Definition}
\label{D:ViscosityFinite}
A function $u\colon[0,T]\times\R^d\rightarrow\R$ is called a \textbf{strong-viscosity solution} to the nonlinear Kolmogorov equation \eqref{E:KolmEqFinite} if there exists a sequence $(u_n,g_n,f_n,b_n,\sigma_n)_n$ satisfying:
\begin{enumerate}
\item[\textup{(i)}] $u_n\colon[0,T]\times\R^d\rightarrow\R$, $g_n\colon\R^d\rightarrow\R$, $f_n\colon[0,T]\times\R^d\times\R\times\R^d\rightarrow\R$, $b_n\colon[0,T]\times\R^d\rightarrow\R$, and $\sigma_n\colon[0,T]\times\R^d\rightarrow\R^{d\times d}$ are locally equicontinuous functions such that, for some positive constants $C$ and $m$, independent of $n$,
\begin{align*}
|b_n(t,x)-b_n(t,x')| + |\sigma_n(t,x)-\sigma_n(t,x')| \ &\leq \ C|x-x'|, \\
|f_n(t,x,y,z)-f_n(t,x,y',z')| \ &\leq \ C\big(|y-y'| + |z-z'|\big), \\
|b_n(t,0)| + |\sigma_n(t,0)| \ &\leq \ C, \\
|u_n(t,x)| + |g_n(x)| + |f_n(t,x,0,0)| \ &\leq \ C\big(1 + |x|^m\big),
\end{align*}
for all $t\in[0,T]$, $x,x'\in\R^d$, $y,y'\in\R$, and $z,z'\in\R^d$.
\item[\textup{(ii)}] $u_n$ is a classical solution to
\[
\begin{cases}
- \partial_t u_n(t,x) - \langle b_n(t,x),D_x u_n(t,x)\rangle - \frac{1}{2}\textup{tr}(\sigma_n\sigma_n\trans(t,x)D_x^2 u_n(t,x)) & \\
\hspace{2.8cm}-\, f_n(t,x,u_n(t,x),\sigma_n\trans(t,x)D_x u_n(t,x)) \ = \ 0, &\!\!\!\!\forall\,(t,x)\in[0,T)\times\R^d, \\
u_n(T,x) \ = \ g_n(x), &\!\!\!\!\forall\,x\in\R^d.
\end{cases}
\]
\item[\textup{(iii)}] $(u_n(t,x),g_n(x),f_n(t,x,y,z),b_n(t,x),\sigma_n(t,x))\rightarrow(u(t,x),g(x),f(t,x,y,z),b(t,x),\sigma(t,x))$, as $n$ tends to infinity, for any $(t,x,y,z)\in[0,T]\times\R^d\times\R\times\R^d$.
\end{enumerate}
\end{Definition}

\begin{Theorem}
Let $u\colon[0,T]\times\R^d\rightarrow\R$ be a strong-viscosity solution to the nonlinear Kolmogorov equation \eqref{E:KolmEqFinite}. Then, we have
\[
u(t,x) \ = \ Y_t^{t,x}, \qquad \forall\,(t,x)\in[0,T]\times\R^d,
\]
where $(Y_s^{t,x},Z_s^{t,x})_{s\in[t,T]}\in\S^2(t,T)\times\H^2(t,T)^d$, with $Y_s^{t,x}=u(s,X_s^{t,x})$, is the solution to the backward stochastic differential equation
\[
Y_s^{t,x} \ = \ g(X_T^{t,x}) + \int_s^T f(r,X_r^{t,x},Y_r^{t,x},Z_r^{t,x}) dr - \int_s^T Z_r^{t,x} dW_r,
\]
for all $t \leq s \leq T$, $\P$-almost surely. In particular, there exists at most one strong solution to the nonlinear Kolmogorov equation \eqref{E:KolmEqFinite}.
\end{Theorem}
\textbf{Proof.}
Since a strong-viscosity solution is in particular a generalized strong-viscosity solution (see Definition \ref{D:StrongSuperSub} below), the thesis follows from Corollary \ref{C:CompThm}.
\ep

\subsubsection{Second definition of strong-viscosity solution}

Our second definition of strong-viscosity solution to equation \eqref{E:KolmEqFinite} is more in the spirit of the standard definition of viscosity solution, which is usually required to be both a viscosity subsolution and a viscosity supersolution. Indeed, we introduce the concept of \emph{generalized strong-viscosity solution}, which has to  be both a strong-viscosity subsolution and a strong-viscosity supersolution. As it will be clear from the definition, this new notion of solution is more general (in other words, weaker), than the concept of strong-viscosity solution given earlier in Definition \ref{D:ViscosityFinite}. For this reason, we added the adjective \emph{generalized} to its name.

\vspace{3mm}

\noindent Firstly, we need to introduce the standard concepts of classical sub and supersolution.
\begin{Definition}
A function $u\colon[0,T]\times\R^d\rightarrow\R$, with $u\in C^{1,2}([0,T[\times\R^d)\cap C([0,T]\times\R^d)$, is called a \textbf{classical supersolution} $($resp. \textbf{classical subsolution}$)$ to the nonlinear Kolmogorov equation \eqref{E:KolmEqFinite} if $u$ solves
\[
\begin{cases}
- \partial_t u(t,x) - \langle b(t,x),D_x u(t,x)\rangle - \frac{1}{2}\textup{tr}(\sigma\sigma\trans(t,x)D_x^2 u(t,x)) & \\
\hspace{1.1cm}-\, f(t,x,u(t,x),\sigma\trans(t,x)D_x u(t,x)) \ \geq \ (\text{resp. $\leq$}) \ 0, &\forall\,(t,x)\in[0,T)\times\R^d, \\
u(T,x) \ \geq \ (\text{resp. $\leq$}) \ g(x), &\forall\,x\in\R^d.
\end{cases}
\]
\end{Definition}

\noindent The following probabilistic representation result for classical sub and supersolutions holds.

\begin{Proposition}
\label{P:UniqClassicalFinite}
Suppose that the functions $b$, $\sigma$, $f$, and $g$, appearing in the nonlinear Kolmogorov equation \eqref{E:KolmEqFinite}, satisfy, for some positive constants $C$ and $m$,
\begin{align*}
|b(t,x)-b(t,x')| + |\sigma(t,x)-\sigma(t,x')| \ &\leq \ C|x-x'|, \\
|f(t,x,y,z)-f(t,x,y',z')| \ &\leq \ C\big(|y-y'| + |z-z'|\big), \\
|b(t,0)| + |\sigma(t,0)| \ &\leq \ C, \\
|f(t,x,0,0)| + |g(x)| \ &\leq \ C\big(1 + |x|^m\big),
\end{align*}
for all $t\in[0,T]$, $x,x'\in\R^d$, $y,y'\in\R$, and $z,z'\in\R^d$.\\
\textup{(i)} Let $u\colon[0,T]\times\R^d\rightarrow\R$ be a classical supersolution to the nonlinear Kolmogorov equation \eqref{E:KolmEqFinite}, satisfying the polynomial growth condition:
\[
|u(t,x)| \ \leq \ C\big(1 + |x|^m\big), \qquad \forall\,(t,x)\in[0,T]\times\R^d.
\]
Then, we have
\[
u(t,x) \ = \ Y_t^{t,x}, \qquad \forall\,(t,x)\in[0,T]\times\R^d,
\]
for some uniquely determined $(Y_s^{t,x},Z_s^{t,x},K_s^{t,x})_{s\in[t,T]}\in\S^2(t,T)\times\H^2(t,T)^d\times\A^{+,2}(t,T)$, with $(Y_s^{t,x},Z_s^{t,x})=(u(s,X_s^{t,x}),\sigma\trans(s,X_s^{t,x})D_x u(s,X_s^{t,x})1_{[t,T[}(s))$, solving the backward stochastic differential equation, $\P$-a.s.,
\[
Y_s^{t,x} \ = \ Y_T^{t,x} + \int_s^T f(r,X_r^{t,x},Y_r^{t,x},Z_r^{t,x}) dr + K_T^{t,x} - K_s^{t,x} - \int_s^T \langle Z_r^{t,x},dW_r\rangle, \quad t \leq s \leq T.
\]
\textup{(ii)} Let $u\colon[0,T]\times\R^d\rightarrow\R$ be a classical subsolution to the nonlinear Kolmogorov equation \eqref{E:KolmEqFinite}, satisfying the polynomial growth condition:
\[
|u(t,x)| \ \leq \ C\big(1 + |x|^m\big), \qquad \forall\,(t,x)\in[0,T]\times\R^d.
\]
Then, we have
\[
u(t,x) \ = \ Y_t^{t,x}, \qquad \forall\,(t,x)\in[0,T]\times\R^d,
\]
for some uniquely determined $(Y_s^{t,x},Z_s^{t,x},K_s^{t,x})_{s\in[t,T]}\in\S^2(t,T)\times\H^2(t,T)^d\times\A^{+,2}(t,T)$, with $(Y_s^{t,x},Z_s^{t,x})=(u(s,X_s^{t,x}),\sigma\trans(s,X_s^{t,x})D_x u(s,X_s^{t,x})1_{[t,T[}(s))$, solving the backward stochastic differential equation, $\P$-a.s.,
\[
Y_s^{t,x} \ = \ Y_T^{t,x} + \int_s^T f(r,X_r^{t,x},Y_r^{t,x},Z_r^{t,x}) dr - (K_T^{t,x} - K_s^{t,x}) - \int_s^T \langle Z_r^{t,x},dW_r\rangle, \quad t \leq s \leq T.
\]
\end{Proposition}
\textbf{Proof.}
The proof can be done along the lines of Theorem \ref{T:UniqClassical}, using Proposition \ref{P:EstimateBSDEAppendix} in the full general case with the presence of the process $K$.
\ep

\vspace{3mm}

\noindent We can now provide the definition of generalized strong-viscosity solution.

\begin{Definition}
\label{D:StrongSuperSub}
A function $u\colon[0,T]\times\R^d\rightarrow\R$ is called a \textbf{strong-viscosity supersolution} $($resp. \textbf{strong-viscosity subsolution}$)$ to the nonlinear Kolmogorov equation \eqref{E:KolmEqFinite} if there exists a sequence $(u_n,g_n,f_n,b_n,\sigma_n)_n$ satisfying:
\begin{enumerate}
\item[\textup{(i)}] $u_n\colon[0,T]\times\R^d\rightarrow\R$, $g_n\colon\R^d\rightarrow\R$, $f_n\colon[0,T]\times\R^d\times\R\times\R^d\rightarrow\R$, $b_n\colon[0,T]\times\R^d\rightarrow\R$, and $\sigma_n\colon[0,T]\times\R^d\rightarrow\R^{d\times d}$ are locally equicontinuous functions such that, for some positive constants $C$ and $m$, independent of $n$,
\begin{align*}
|b_n(t,x)-b_n(t,x')| + |\sigma_n(t,x)-\sigma_n(t,x')| \ &\leq \ C|x-x'|, \\
|f_n(t,x,y,z)-f_n(t,x,y',z')| \ &\leq \ C\big(|y-y'| + |z-z'|\big), \\
|b_n(t,0)| + |\sigma_n(t,0)| \ &\leq \ C, \\
|u_n(t,x)| + |g_n(x)| + |f_n(t,x,0,0)| \ &\leq \ C\big(1 + |x|^m\big),
\end{align*}
for all $t\in[0,T]$, $x,x'\in\R^d$, $y,y'\in\R$, and $z,z'\in\R^d$.
\item[\textup{(ii)}] $u_n$ is a classical supersolution $($resp. classical subsolution$)$ to
\[
\begin{cases}
- \partial_t u_n(t,x) - \langle b_n(t,x),D_x u_n(t,x)\rangle - \frac{1}{2}\textup{tr}(\sigma_n\sigma_n\trans(t,x)D_x^2 u_n(t,x)) & \\
\hspace{2.8cm}-\, f_n(t,x,u_n(t,x),\sigma_n\trans(t,x)D_x u_n(t,x)) \ = \ 0, &\!\!\!\!\!\!\!\forall\,(t,x)\in[0,T)\times\R^d, \\
u_n(T,x) \ = \ g_n(x), &\!\!\!\!\!\!\!\forall\,x\in\R^d.
\end{cases}
\]
\item[\textup{(iii)}] $(u_n(t,x),g_n(x),f_n(t,x,y,z),b_n(t,x),\sigma_n(t,x))\rightarrow(u(t,x),g(x),f(t,x,y,z),b(t,x),\sigma(t,x))$, as $n$ tends to infinity, for any $(t,x,y,z)\in[0,T]\times\R^d\times\R\times\R^d$.
\end{enumerate}
A function $u\colon[0,T]\times\R^d\rightarrow\R$ is called a \textbf{generalized strong-viscosity solution} to the nonlinear Kolmogorov equation \eqref{E:KolmEqFinite} if it is both a strong-viscosity supersolution and a strong-viscosity subsolution to \eqref{E:KolmEqFinite}.
\end{Definition}

We can now state the following probabilistic representation result for strong-viscosity sub and supersolutions, that is one of the main results of this paper, from which the comparison theorem will follow as corollary.

\begin{Theorem}
\label{T:RepresentationSuperSub}
\textup{(1)} Let $u\colon[0,T]\times\R^d\rightarrow\R$ be a strong-viscosity supersolution to the nonlinear Kolmogorov equation \eqref{E:KolmEqFinite}. Then, we have
\[
u(t,x) \ = \ Y_t^{t,x}, \qquad \forall\,(t,x)\in[0,T]\times\R^d,
\]
for some uniquely determined $(Y_s^{t,x},Z_s^{t,x},K_s^{t,x})_{s\in[t,T]}\in\S^2(t,T)\times\H^2(t,T)^d\times\A^{+,2}(t,T)$, with $Y_s^{t,x}=u(s,X_s^{t,x})$, solving the backward stochastic differential equation, $\P$-a.s.,
\begin{align}
\label{E:BSDE_StrongNonlinear_Super}
Y_s^{t,x} \ &= \ Y_T^{t,x} + \int_s^T f(r,X_r^{t,x},Y_r^{t,x},Z_r^{t,x}) dr \\
&\quad \ + K_T^{t,x} - K_s^{t,x} - \int_s^T \langle Z_r^{t,x},dW_r\rangle, \qquad\qquad\qquad t \leq s \leq T. \notag
\end{align}
\textup{(2)} Let $u\colon[0,T]\times\R^d\rightarrow\R$ be a strong-viscosity subsolution to the nonlinear Kolmogorov equation \eqref{E:KolmEqFinite}. Then, we have
\[
u(t,x) \ = \ Y_t^{t,x}, \qquad \forall\,(t,x)\in[0,T]\times\R^d,
\]
for some uniquely determined $(Y_s^{t,x},Z_s^{t,x},K_s^{t,x})_{s\in[t,T]}\in\S^2(t,T)\times\H^2(t,T)^d\times\A^{+,2}(t,T)$, with $Y_s^{t,x}=u(s,X_s^{t,x})$, solving the backward stochastic differential equation, $\P$-a.s.,
\begin{align}
\label{E:BSDE_StrongNonlinear_Sub}
Y_s^{t,x} \ &= \ Y_T^{t,x} + \int_s^T f(r,X_r^{t,x},Y_r^{t,x},Z_r^{t,x}) dr \\
&\quad \ - \big(K_T^{t,x} - K_s^{t,x}\big) - \int_s^T \langle Z_r^{t,x},dW_r\rangle, \qquad\qquad\qquad t \leq s \leq T. \notag
\end{align}
\end{Theorem}
\textbf{Proof.}
We shall prove statement (1), since (2) can be proved similarly. To prove (1), consider a sequence $(u_n,g_n,f_n,b_n,\sigma_n)_n$ satisfying conditions (i)-(iii) of Definition \ref{D:StrongSuperSub}. For every $n\in\N$ and any $(t,x)\in[0,T]\times\R^d$, consider the stochastic equation, $\P$-a.s.,
\[
X_s \ = \ x + \int_t^s b_n(r,X_r) dr + \int_t^s \sigma_n(r,X_r) dW_r, \qquad t \leq s \leq T.
\]
It is well known that there exists a unique solution $(X_s^{n,t,x})_{s\in[t,T]}$ to the above equation. Moreover, from Proposition \ref{P:UniqClassicalFinite} we know that $u_n(t,x) = Y_t^{n,t,x}$, $(t,x)\in[0,T]\times\R^d$, for some $(Y_s^{n,t,x},Z_s^{n,t,x},K_s^{n,t,x})_{s\in[t,T]}\in\S^2(t,T)\times\H^2(t,T)^d\times\A^{+,2}(t,T)$ solving the backward stochastic differential equation, $\P$-a.s.,
\begin{align*}
Y_s^{n,t,x} \ &= \ Y_T^{n,t,x} + \int_s^T f_n(r,X_r^{n,t,x},Y_r^{n,t,x},Z_r^{n,t,x}) dr \\
&\quad \ + K_T^{n,t,x} - K_s^{n,t,x} - \int_s^T \langle Z_r^{n,t,x}, dW_r\rangle, \qquad\qquad\qquad t \leq s \leq T.
\end{align*}
Notice that, from the uniform polynomial growth condition of $(u_n)_n$ and estimate \eqref{E:EstimateX^n} we have, for any $p\geq1$,
\[
\sup_{n\in\N} \|Y^{n,t,x}\|_{\S^p(t,T)} \ < \ \infty.
\]
Then, it follows from Proposition \ref{P:EstimateBSDEAppendix}, the polynomial growth condition of $(f_n)_n$ in $x$, and the linear growth condition of $(f_n)_n$ in $(y,z)$, that
\[
\sup_n\big(\|Z^{n,t,x}\|_{\H^2(t,T)^d} + \|K^{n,t,x}\|_{\S^2(t,T)}\big) \ < \ \infty.
\]
Set $Y_s^{t,x}=u(s,X_s^{t,x})$, for any $s\in[t,T]$. Then, from the polynomial growth condition that $u$ inherits from the sequence $(u_n)_n$, and using estimate \eqref{E:EstimateX^n}, we deduce that $\|Y^{t,x}\|_{\S^p(t,T)}<\infty$, for any $p\geq1$. In particular, $Y\in\S^2(t,T)$ and it is continuous process. We also have, using
the convergence result \eqref{E:X^n-->X}, that there exists a subsequence of $(X^{n,t,x})_n$, which we still denote $(X^{n,t,x})_n$, such that
\begin{equation}
\label{E:supX^n-X-->0}
\sup_{t\leq s\leq T}|X_s^{n,t,x}(\omega)-X_s^{t,x}(\omega)| \ \overset{n\rightarrow\infty}{\longrightarrow} \ 0, \qquad \forall\,\omega\in\Omega\backslash N,
\end{equation}
for some null measurable set $N\subset\Omega$. Moreover, from  estimate \eqref{E:EstimateX^n} it follows that, possibly enlarging $N$, $\sup_{t\leq s\leq T}(|X_s^{n,t,x}(\omega)| + |X_s^{t,x}(\omega)|)<\infty$, for any $n\in\N$ and any $\omega\in\Omega\backslash N$. Now, fix $\omega\in\Omega\backslash N$, then
\begin{align*}
&|Y_s^{n,t,x}(\omega)-Y_s^{t,x}(\omega)| \ = \ |u_n(s,X_s^{n,t,x}(\omega))-u(s,X_s^{t,x}(\omega))| \\
&= \ |u_n(s,X_s^{n,t,x}(\omega))-u_n(s,X_s^{t,x}(\omega))| + |u_n(s,X_s^{t,x}(\omega))-u(s,X_s^{t,x}(\omega))|.
\end{align*}
For any $\eps>0$, from point (iii) of Definition \ref{D:StrongSuperSub} it follows that there exists $n'\in\N$ such that
\[
|u_n(s,X_s^{t,x}(\omega))-u(s,X_s^{t,x}(\omega))| \ < \ \frac{\eps}{2}, \qquad \forall\,n\geq n'.
\]
On the other hand, from the local equicontinuity of $(u_n)_n$, we see that there exists $\delta>0$, independent of $n$, such that
\[
|u_n(s,X_s^{n,t,x}(\omega))-u_n(s,X_s^{t,x}(\omega))| \ < \ \frac{\eps}{2}, \qquad \text{if }|X_s^{n,t,x}(\omega)-X_s^{t,x}(\omega)|<\delta.
\]
Using \eqref{E:supX^n-X-->0}, we can find $n''\in\N$, $n''\geq n'$, such that
\[
\sup_{t\leq s\leq T}|X_s^{n,t,x}(\omega)-X_s^{t,x}(\omega)| \ < \ \delta, \qquad \forall\,n\geq n''.
\]
In conclusion, for any $\omega\in\Omega\backslash N$ and any $\eps>0$ there exists $n''\in\N$ such that
\[
|Y_s^{n,t,x}(\omega)-Y_s^{t,x}(\omega)| \ < \ \eps, \qquad \forall\,n\geq n''.
\]
Therefore, $Y_s^{n,t,x}(\omega)$ converges to $Y_s^{t,x}(\omega)$, as $n$ tends to infinity, for any $(s,\omega)\in[t,T]\times(\Omega\backslash N)$. In a similar way, we can prove that there exists a null measurable set $N'\subset\Omega$ such that $f_n(s,X_s^{n,t,x}(\omega),y,z)\rightarrow f(s,X_s^{t,x}(\omega),y,z)$, for any $(s,\omega,y,z)\in[t,T]\times(\Omega\backslash N')\times\R\times\R^d$. As a consequence, the thesis follows from Proposition \ref{P:LimitThmBSDE}.
\ep

\vspace{3mm}

We can finally state a comparison theorem for strong-viscosity sub and supersolutions, which follows directly from the comparison theorem for BSDEs.

\begin{Corollary}[Comparison Theorem]
\label{C:CompThm}
Let $\check u\colon[0,T]\times\R^d\rightarrow\R$ $($resp. $\hat u\colon[0,T]\times\R^d\rightarrow\R$$)$ be a strong-viscosity subsolution $($resp. strong-viscosity supersolution$)$ to the nonlinear Kolmogorov equation \eqref{E:KolmEq}. Then $\check u \leq \hat u$ on $[0,T]\times\R^d$. In particular, there exists at most one generalized strong-viscosity solution to the nonlinear Kolmogorov equation \eqref{E:KolmEqFinite}.
\end{Corollary}
\textbf{Proof.}
We know that $\check u(T,x) \leq g(x) \leq \hat u(T,x)$, for all $x\in\R^d$. Moreover, from Theorem \ref{T:RepresentationSuperSub} we have
\[
\check u(t,x) \ = \ \check Y_t^{t,x}, \qquad \hat u(t,x) \ = \ \hat Y_t^{t,x}, \qquad \text{for all }(t,x)\in[0,T]\times\R^d,
\]
for some $(\check Y_s^{t,x},\check Z_s^{t,x},\check K_s^{t,x})_{s\in[t,T]},(\hat Y_s^{t,x},\hat Z_s^{t,x},\hat K_s^{t,x})_{s\in[t,T]}\in\S^2(t,T)\times\H^2(t,T)^d\times\A^{+,2}(t,T)$ satisfying \eqref{E:BSDE_StrongNonlinear_Sub} and \eqref{E:BSDE_StrongNonlinear_Super}, respectively. Let us denote $\overline{Y}:=\hat Y^{t,x}-\check Y^{t,x}$, $\overline{Z}:=\hat Z^{t,x}-\check Z^{t,x}$, $\overline K:=\hat K^{t,x}+\check K^{t,x}$, and $\bar{f}_{s}:=f(s,X_s^{t,x},\hat Y_s^{t,x},\hat Z_s^{t,x})-f(s,X_s^{t,x},\check Y_s^{t,x},\check Z_s^{t,x})$. Then, $\P$-a.s.,
\begin{equation}
\label{BSDEComp}
\overline{Y}_s \ = \ \overline{Y}_T + \int_s^T \bar{f}_r dr + \overline K_T - \overline K_s - \int_s^T  \overline{Z}_r dW_r, \qquad t \leq s \leq T.
\end{equation}
Now we introduce the real predictable process $a$ given by, $\P$-a.s.,
\[
a_s \ = \ \frac{f(s,X_s^{t,x},\hat Y_s^{t,x},\hat Z_s^{t,x})-f(s,X_s^{t,x},\check Y_s^{t,x},\hat Z_s^{t,x})}{\hat Y_s^{t,x}-\check Y_s^{t,x}}
1_{\{\hat Y_s^{t,x}-\check Y_s^{t,x}\neq0\}}, \qquad   t \leq s \leq T
\]
and the $\mathbb{R}^{d}$-valued predictable process $b$ defined componentwise by, $\P$-a.s.,
\[
b_s^k \ = \ \frac{f(s,X_s^{t,x},\check Y_s^{t,x},Z_s^{(k-1)})-f(s,X_s^{t,x},\check Y_s^{t,x},Z_s^{(k)})}{V_s^k}
1_{\{V_s^k\neq0\}}, \qquad   t \leq s \leq T,
\]
for all $k=1,\ldots,d$, where $Z^{(k)}$ is the $\mathbb{R}^{d}$-valued process whose $k$ first components are those of $\check Z^{t,x}$ and whose $(d-k)$ lasts are those of $\hat Z^{t,x}$, and $V^k$ is the $k$-th component of $Z^{(k-1)}-Z^{(k)}$. Notice that, since $f$ is uniformly Lipschitz in $(y,z)$, the processes $a$ and $b$ are bounded.

Equation \eqref{BSDEComp} can be rewritten as
\[
\overline{Y}_s \ = \ \overline{Y}_T + \int_s^T \big(a_r\overline{Y}_r + \langle b_r,\overline{Z}_r\rangle\big)dr + \overline{K}_T - \overline{K}_s - \int_s^T \overline{Z}_r dW_r.
\]
Consider now the process $\Gamma=(\Gamma_s)_{t \leq s \leq T}$ satisfying, $\P$-a.s.,
\[
\Gamma_s \ = \ 1 + \int_t^s \Gamma_r a_r dr + \int_t^s \Gamma_r\langle b_r,dW_r\rangle, \qquad t \leq s \leq T.
\]
Notice that $\Gamma\in\S^2(t,T)$, since $a$ and $b$ are bounded. Moreover, $\Gamma$ is strictly positive. An application of It\^o's formula yields, $\P$-a.s.,
\begin{equation}
\label{E:GammaY}
\Gamma_s\overline{Y}_s \ = \ \overline Y_t - \int_t^s \Gamma_r d\overline{K}_r + \int_t^s \Gamma_r\langle\overline{Y}_r b_r + \overline{Z}_r,dW_r\rangle, \qquad t \leq s \leq T.
\end{equation}
From Burkholder-Davis-Gundy inequality it follows that there exists a positive constant $C$ such that
\begin{align*}
\E\bigg[\sup_{t \leq s \leq T}\bigg|\int_t^s \Gamma_r\langle\overline{Y}_r b_r + \overline{Z}_r,dW_r\rangle\bigg|\bigg] \ &\leq \ C\E\bigg[\bigg(\int_t^T \Gamma_s^2|\overline{Y}_s b_s + \overline{Z}_s|^2 ds\bigg)^{\frac{1}{2}}\bigg] \\
&\leq \ \frac{C}{2} \E\bigg[\sup_{t \leq s \leq T}\Gamma_s^2 +  2b_\infty^2\int_t^T |\overline{Y}_s|^2 ds + 2\int_t^T |\overline{Z}_s|^2 ds\bigg] \\
&< \ \infty,
\end{align*}
where $b_\infty$ is an upper-bound for $b$. Therefore, the local martingale $(\int_t^s \Gamma_r\langle\overline{Y}_r b_r + \overline{Z}_r,dW_r\rangle)_{t \leq s \leq T}$ is indeed a martingale. Hence, we see from \eqref{E:GammaY} that the process $(\Gamma_s\overline{Y}_s)_{t \leq s \leq T}$ is a supermartingale, which implies (recalling that $\hat u(T,x) - \check u(T,x) \geq 0$, for all $x\in\R^d$)
\[
\hat u(t,x) - \check u(t,x) \ = \ \overline{Y}_t \ \geq \ \mathbb{E}\big[\Gamma_T\overline{Y}_T\big] \ \geq \ 0.
\]
\ep

\subsubsection{Relation with the standard definition of viscosity solution}

We now conclude this brief digression concerning strong-viscosity solutions, investigating the equivalence between the concept of strong-viscosity solution and the standard notion of viscosity solution, for which we refer, e.g., to \cite{crandishiilions92}. Let us begin recalling the definition of viscosity solution for equation \eqref{E:KolmEqFinite}.

\begin{Definition}
\textup{(i)} A lower $($resp. upper$)$ semicontinuous function $u\colon[0,T]\times\R^d\rightarrow\R$ is called a \textbf{viscosity supersolution} $($resp. \textbf{viscosity subsolution}$)$ to the nonlinear Kolmogorov equation \eqref{E:KolmEq} if
\[
u(T,x) \ \geq \ (\text{resp. $\leq$}) \ g(x), \qquad \forall\,x\in\R^d
\]
and
\begin{align*}
- \partial_t \varphi(t,x) - \langle b(t,x),D_x \varphi(t,x)\rangle - \frac{1}{2}\textup{tr}\big(\sigma\sigma\trans(t,x)D_x^2 \varphi(t,x)\big) & \\
- f\big(t,x,u(t,x),\sigma\trans(t,x)D_x \varphi(t,x)\big)& \ \geq \ (\text{resp. $\leq$}) \ 0,
\end{align*}
for any $(t,x)\in[0,T)\times\R^d$ and any $\varphi\in C^{1,2}([0,T]\times\R^d)$ such that $u-\varphi$ has a local minimum $($resp. maximum$)$ at $(t,x)$.\\
\textup{(ii)} A continuous function $u\colon[0,T]\times\R^d\rightarrow\R$ is called a \textbf{viscosity solution} to the nonlinear Kolmogorov equation \eqref{E:KolmEq} if it is both a viscosity supersolution and a viscosity subsolution to \eqref{E:KolmEq}.
\end{Definition}

\begin{Theorem}
Suppose that the functions $b$, $\sigma$, $f$, and $g$, appearing in the nonlinear Kolmogorov equation \eqref{E:KolmEqFinite}, are bounded and satisfy, for some positive constant $C$ and continuity modulus $\rho$,
\begin{align*}
|b(t,x) - b(t',x')| + |\sigma(t,x) - \sigma(t',x')| \ &\leq \ \rho(|t-t'|) + C|x-x'|, \\
|f(t,x,y,z) - f(t',x',y',z')| \ &\leq \ \rho(|t-t'| + |x-x'|) + C(|y-y'| + |z-z'|), \\
|g(x) - g(x')| \ &\leq \ \rho(|x-x'|),
\end{align*}
for all $(t,x,y,z),(t',x',y',z')\in[0,T]\times\R^d\times\R\times\R^d$. Suppose, moreover, that $\sigma(t,x)$ is a positive semidefinite matrix, for all $(t,x)\in[0,T]\times\R^d$. Let $u\colon[0,T]\times\R^d\rightarrow\R$ be bounded and uniformly continuous with a continuity modulus $\rho$. Then $u$ is a viscosity solution to the nonlinear Kolmogorov equation \eqref{E:KolmEqFinite} if and only if $u$ is a strong-viscosity solution to \eqref{E:KolmEqFinite}.
\end{Theorem}
\textbf{Proof.}
The \emph{if} part follows from the method of half-relaxed limits of Barles and Perthame (see, e.g., Lemma 6.1 and Remark 6.3 in \cite{crandishiilions92}). Let us focus on the \emph{only if} part.

It is well-known that, under the present assumptions, a uniqueness result for viscosity solutions to equation \eqref{E:KolmEqFinite} holds, see, e.g., Theorem 5.1 in \cite{pardoux_pradeilles_rao97} for the case $b$ and $\sigma$ independent of $t$. Let $\phi(x) = c\exp(1/(|x|^2-1))1_{\{|x|<1\}}(x)$, $x\in\R^d$, with $c>0$ such that $\int_{\R^d}\phi(x)dx=1$. Define, for any $n\in\N\backslash\{0\}$, $\phi_n(x) = n\phi(nx)$, $x\in\R^d$. Now, define (denoting by $I$ the $d\times d$ identity matrix)
\begin{align*}
b_n(t,x) \ := \ \int_{\R^d} \phi_n(x-x')b(t,x') dx', \qquad \sigma_n(t,x) \ &:= \ \int_{\R^d} \phi_n(x-x')\sigma(t,x') dx' + \frac{1}{n}I, \\
g_n(x) \ &:= \ \int_{\R^d} \phi_n(x-x')g(x') dx'.
\end{align*}
Similarly, let $\psi(x,y,z) = \bar c\exp(1/(|(x,y,z)|^2-1))1_{\{|(x,y,z)|<1\}}(x,y,z)$, $(x,y,z)\in\R^d\times\R\times\R^d$, with $\bar c>0$ such that $\int_{\R^d\times\R\times\R^d}\psi(x,y,z)dxdydz=1$. Define, for any $n\in\N\backslash\{0\}$, $\psi_n(x,y,z) = n\psi(nx,ny,nz)$, $(x,y,z)\in\R^d\times\R\times\R^d$. Now, define
\[
f_n(t,x,y,z) \ := \ \int_{\R^d\times\R\times\R^d} \psi_n(x-x',y-y',z-z')f(t,x',y',z') dx'dy'dz'.
\]
Notice that $b_n,\sigma_n\in C^{0,\infty}([0,T]\times\R^d)$, $g_n\in C^\infty(\R^d)$, and $f_n\in C^{0,\infty}([0,T]\times\R^{2d+1})$, with
\[
(b_n(t,x),\sigma_n(t,x),g_n(x),f_n(t,x,y,z)) \ \overset{n\rightarrow\infty}{\longrightarrow} \ (b(t,x),\sigma(t,x),g(x),f(t,x,y,z)),
\]
for all $(t,x,y,z)\in[0,T]\times\R^d\times\R\times\R^d$. Moreover, $b_n$, $\sigma_n$, $g_n$, and $f_n$ are bounded and satisfy
\begin{align*}
|b_n(t,x) - b_n(t',x')| + |\sigma_n(t,x) - \sigma_n(t',x')| \ &\leq \ \rho(|t-t'|) + C|x-x'|, \\
|f_n(t,x,y,z) - f_n(t',x',y',z')| \ &\leq \ \rho(|t-t'| + |x-x'|) + C(|y-y'| + |z-z'|), \\
|g_n(x) - g_n(x')| \ &\leq \ \rho(|x-x'|),
\end{align*}
for all $(t,x,y,z),(t',x',y',z')\in[0,T]\times\R^d\times\R\times\R^d$, with the same constant $C$ and continuity modulus $\rho$ as in the statement of the theorem. Let us now consider, for each $n\in\N\backslash\{0\}$, the nonlinear Kolmogorov equation:
\begin{equation}
\label{E:KolmEqFinite_n}
\begin{cases}
- \partial_t u_n(t,x) - \langle b_n(t,x),D_x u_n(t,x)\rangle - \frac{1}{2}\textup{tr}(\sigma_n\sigma_n\trans(t,x)D_x^2 u_n(t,x)) & \\
\hspace{2cm}-\, f_n(t,x,u_n(t,x),\sigma_n\trans(t,x)D_x u_n(t,x)) \ = \ 0, &\hspace{-1cm}\forall\,(t,x)\in[0,T)\times\R^d, \\
u_n(T,x) \ = \ g_n(x), &\hspace{-1cm}\forall\,x\in\R^d.
\end{cases}
\end{equation}
Since $\sigma(t,x)$ is a positive semidefinite matrix, we see that $\sigma_n(t,x)$ is a positive definite matrix and equation \eqref{E:KolmEqFinite_n} is uniformly elliptic. Then, it follows from classical results on regularity theory for parabolic equations (see, e.g., Theorem 8.1 in \cite{ladyzenskaja}) that there exists a unique classical solution $u_n\in C^{1,2}([0,T]\times\R^d)$ to equation \eqref{E:KolmEqFinite_n}, with $\sup_{[0,T]\times\R^d}|u_n|\leq M$ for some positive constant $M$, independent of $n$, as it can be seen using the uniform boundedness of $b_n$, $\sigma_n$, $g_n$, and $f_n$. Clearly, $u_n$ is also a viscosity solution to equation \eqref{E:KolmEqFinite_n}. Then, using the notations of Section 6 in \cite{crandishiilions92}, set
\[
\overline u(t,x) \ := \ \limsup_{n\rightarrow\infty}\!{}^*\,u_n(t,x), \qquad \underline u(t,x) \ := \ \liminf_{n\rightarrow\infty}\!{}_*\,u_n(t,x).
\]
Notice that $|\overline u|,|\underline u|\leq M$. From Remark 6.3 in \cite{crandishiilions92}, we know that $\overline u$ (resp. $\underline u$) is a viscosity subsolution (resp. supersolution) to equation \eqref{E:KolmEqFinite}. From the comparison theorem for viscosity solutions, this implies that $\overline u \leq \underline u$. Since $\underline u \leq \overline u$ by definition, we get $\overline u$ and $\underline u$ are equal and are both viscosity solutions to equation \eqref{E:KolmEqFinite}. From uniqueness, we must have $u = \overline u = \underline u$. Moreover, it follows from Remark 6.4 in \cite{crandishiilions92} that $u_n$ converges to $u$ uniformly on compact sets. In conclusion, $u$ is a strong-viscosity solution to equation \eqref{E:KolmEqFinite}.
\ep

\appendix

\setcounter{equation}{0} \setcounter{Assumption}{0}
\setcounter{Theorem}{0} \setcounter{Proposition}{0}
\setcounter{Corollary}{0} \setcounter{Lemma}{0}
\setcounter{Definition}{0} \setcounter{Remark}{0}

\renewcommand\thesection{Appendix}

\section{}

In the present appendix we fix a complete probability space $(\Omega,\Fc,\P)$ on which a $d$-dimensional Brownian motion $W=(W_t)_{t\geq0}$ is defined. We denote $\F=(\Fc_t)_{t\geq0}$ the completion of the natural filtration generated by $W$.

\renewcommand\thesection{\Alph{subsection}}

\renewcommand\thesubsection{\Alph{subsection}.}

\subsection{Estimates for supersolutions to BSDEs}

We derive estimates for the norm of the $Z$ and $K$ components for supersolutions to backward stochastic differential equations, in terms of the norm of the $Y$ component. These results are standard, but seemingly not at disposal in the following form in the literature. Firstly, let us introduce a generator function $F\colon[0,T]\times\Omega\times\R\times\R^d\rightarrow\R$ satisfying the usual assumptions:
\begin{enumerate}
\item[(A.a)] $F(\cdot,y,z)$ is $\F$-predictable for every $(y,z)\in\R\times\R^d$.
\item[(A.b)] There exists a positive constant $C_F$ such that
\[
|F(s,y,z) - F(s,y',z')| \ \leq \ C_F\big(|y-y'| + |z-z'|\big),
\]
for all $y,y'\in\R$, $z,z'\in\R^d$, $ds\otimes d\P$-a.e.
\item[(A.c)] Integrability condition:
\[
\E\bigg[\int_t^T |F(s,0,0)|^2 ds\bigg] \ \leq \ M_F,
\]
for some positive constant $M_F$.
\end{enumerate}

\begin{Proposition}
\label{P:EstimateBSDEAppendix}
For any $t,T\in\R_+$, $t<T$, consider $(Y_s,Z_s,K_s)_{s\in[t,T]}$ satisfying:
\begin{enumerate}
\item[\textup{(i)}] $Y\in\S^2(t,T)$ and it is continuous.
\item[\textup{(ii)}] $Z$ is an $\R^d$-valued $\F$-predictable process such that $\P(\int_t^T |Z_s|^2 ds<\infty)=1$.
\item[\textup{(iii)}] $K$ is a real nondecreasing $($resp. nonincreasing$)$ continuous $\F$-predictable process such that $K_t = 0$.
\end{enumerate}
Suppose that $(Y_s,Z_s,K_s)_{s\in[t,T]}$ solves the BSDE, $\P$-a.s.,
\begin{equation}
\label{E:BSDE_Appendix}
Y_s \ = \ Y_T + \int_s^T F(r,Y_r,Z_r) dr + K_T - K_s - \int_s^T \langle Z_r, dW_r\rangle, \qquad t \leq s \leq T,
\end{equation}
for some generator function $F$ satisfying conditions \textup{(A.b)-(A.c)}. Then $(Z,K)\in\H^2(t,T)^d\times\A^{+,2}(t,T)$ and
\[
\|Z\|_{\H^2(t,T)^d}^2 + \|K\|_{\S^2(t,T)}^2 \ \leq \ C(1+T^3)\bigg(\|Y\|_{\S^2(t,T)}^2 + \int_t^T |F(s,0,0)|^2 ds\bigg),
\]
for some positive constant $C$ depending only on $C_F$, the Lipschitz constant of $F$.
\end{Proposition}
\textbf{Proof.}
Let us consider the case where $K$ is nondecreasing. For every $k\in\N$, define the stopping time
\[
\tau_k \ = \ \inf\bigg\{s\geq t\colon \int_t^s |Z_r|^2 dr \geq k \bigg\} \wedge T.
\]
Then, the local martingale $(\int_t^s Y_r\langle 1_{[t,\tau_k]}(r)Z_r,dW_r\rangle)_{s\in[t,T]}$ satisfies, using Burkholder-Davis-Gundy inequality,
\[
\E\bigg[\sup_{t \leq s \leq T}\bigg|\int_t^s Y_r\langle 1_{[t,\tau_k]}(r)Z_r,dW_r\rangle\bigg|\bigg] \ < \ \infty,
\]
therefore it is a martingale. As a consequence, an application of It\^o's formula to $|Y_s|^2$ between $t$ and $\tau_k$ yields
\begin{align}
\label{E:Proof_ItoAppendix}
\E\big[|Y_t|^2\big] + \E\int_t^{\tau_k} |Z_r|^2 dr \ &= \ \E\big[|Y_{\tau_k}|^2\big] + 2\E\int_t^{\tau_k} Y_r F(r,Y_r,Z_r) dr + 2 \E\int_t^{\tau_k} Y_r dK_r.
\end{align}
In the sequel $c$ and $c'$ will be two strictly positive constants depending only on $C_F$, the Lipschitz constant of $F$. Using (A.b) and recalling the standard inequality $ab \leq a^2 + b^2/4$, for any $a,b\in\R$, we see that
\begin{align}
\label{E:Proof_LipschitzAppendix}
&2\E\int_t^{\tau_k} Y_r F(r,Y_r,Z_r) dr \notag \\
&\leq \ cT\|Y\|_{\S^2(t,T)}^2 + \frac{1}{4}\E\int_t^{\tau_k} |Z_r|^2 dr + \E\int_t^T |F(r,0,0)|^2 dr.
\end{align}
Regarding the last term on the right-hand side in \eqref{E:Proof_ItoAppendix}, for every $\eps>0$ we have (recalling the standard inequality $2ab \leq \eps a^2 + b^2/\eps$, for any $a,b\in\R$)
\begin{equation}
\label{E:YdK_Appendix}
2\E\int_t^{\tau_k} Y_r dK_r \ \leq \ \frac{1}{\eps}\|Y\|_{\S^2(t,T)}^2 + \eps\E\big[|K_{\tau_k}|^2\big].
\end{equation}
Now, from \eqref{E:BSDE_Appendix} we get
\[
K_{\tau_k} \ = \ Y_t - Y_{\tau_k} - \int_t^{\tau_k} F(r,Y_r,Z_r) dr + \int_t^{\tau_k} \langle Z_r,dW_r\rangle.
\]
Therefore (recalling that $(x_1+\cdots+x_4)\leq4(x_1^2+\cdots+x_4^2)$, for any $x_1,\ldots,x_4\in\R$)
\[
\E\big[|K_{\tau_k}|^2\big] \ \leq \ 8\|Y\|_{\S^2(t,T)}^2 + 4T\E\int_t^{\tau_k} |F(r,Y_r,Z_r)|^2 dr + 4\E\bigg|\int_t^{\tau_k} \langle Z_r,dW_r\rangle\bigg|^2.
\]
From It\^o's isometry and (A.b), we obtain
\begin{align}
\label{E:K<C_Appendix}
\E\big[|K_{\tau_k}|^2\big] \ &\leq \ c'(1+T^2)\|Y\|_{\S^2(t,T)}^2 + c'(1+T)\E\int_t^{\tau_k} |Z_r|^2 dr \notag \\
&\quad \ + c'T\E\int_t^T |F(r,0,0)|^2 dr.
\end{align}
Then, taking $\eps=1/(4c'(1+T))$ in \eqref{E:YdK_Appendix} we get
\begin{align}
\label{E:YdK2_Appendix}
&2\E\int_t^{\tau_k} Y_r dK_r \notag \\
&\leq \ \frac{16c'(1+T)^2+1+T^2}{4(1+T)}\|Y\|_{\S^2(t,T)}^2 + \frac{1}{4}\E\int_t^{\tau_k} |Z_r|^2 dr + \frac{T}{4(1+T)}\E\int_t^T |F(r,0,0)|^2 dr \notag \\
&\leq \ c(1+T^2)\|Y\|_{\S^2(t,T)}^2 + \frac{1}{4}\E\int_t^{\tau_k} |Z_r|^2 dr + cT\E\int_t^T |F(r,0,0)|^2 dr.
\end{align}
Plugging \eqref{E:Proof_LipschitzAppendix} and \eqref{E:YdK2_Appendix} into \eqref{E:Proof_ItoAppendix}, we end up with
\[
\E\big[|Y_{\tau_k}|^2\big] + \frac{1}{2}\E\int_t^{\tau_k} |Z_r|^2 dr \ \leq \ c(1+T^2)\|Y\|_{\S^2(t,T)}^2 + c(1+T)\E\int_t^T |F(r,0,0)|^2 dr.
\]
Then, from monotone convergence theorem,
\begin{equation}
\label{E:Z<C_Appendix}
\E\int_t^T |Z_r|^2 dr \ \leq \ c(1+T^2)\|Y\|_{\S^2(t,T)}^2 + c(1+T)\E\int_t^T |F(r,0,0)|^2 dr.
\end{equation}
Plugging \eqref{E:Z<C_Appendix} into \eqref{E:K<C_Appendix}, and using again monotone convergence theorem, we finally obtain
\[
\|K\|_{\S^2(t,T)}^2 \ = \ \E\big[|K_T|^2\big] \ \leq \ c(1+T^3)\|Y\|_{\S^2(t,T)}^2 + c(1+T^2)\E\int_t^T |F(r,0,0)|^2 dr.
\]
When $K$ is nonincreasing, the proof can be done along the same lines.
\ep

\setcounter{equation}{0} \setcounter{Assumption}{0}
\setcounter{Theorem}{0} \setcounter{Proposition}{0}
\setcounter{Corollary}{0} \setcounter{Lemma}{0}
\setcounter{Definition}{0} \setcounter{Remark}{0}

\subsection{Estimates for stochastic differential equations}

We shall report here a result about stochastic differential equations, whose proof is standard.

\begin{Proposition}
\label{P:Estimate}
For any $n\in\N$, let $b_n\colon[0,T]\times\R^d\rightarrow\R$ and $\sigma_n\colon[0,T]\times\R^d\rightarrow\R^{d\times d}$ be Borel measurable functions, satisfying, for some positive constant $C$, independent of $n$,
\begin{align*}
|b_n(t,x)-b_n(t,x')| + |\sigma_n(t,x)-\sigma_n(t,x')| \ &\leq \ C|x-x'|, \\
|b_n(t,0)| + |\sigma_n(t,0)| \ &\leq \ C,
\end{align*}
for all $t\in[0,T]$, $x,x'\in\R^d$. Then, for any $n\in\N$ and $(t,x)\in[0,T]\times\R^d$ there exists a unique solution $(X_s^{n,t,x})_{s\in[t,T]}$ to the equation
\begin{equation}
\label{E:EquationX^n}
X_s \ = \ x + \int_t^s b_n(r,X_r) dr + \int_t^s \sigma_n(r,X_r) dW_r, \qquad t \leq s \leq T,\,\P\text{-a.s.}
\end{equation}
Moreover, suppose that, for every $(t,x)\in[0,T]\times\R^d$, the sequence $\{(b_n(t,x),\sigma_n(t,x))\}_n$ converges as $n$ goes to infinity and define
\[
(b(t,x),\sigma(t,x)) \ := \ \lim_{n\rightarrow\infty}(b_n(t,x),\sigma_n(t,x)).
\]
Then, for any $(t,x)\in[0,T]\times\R^d$ there exists a unique solution $(X_s^{t,x})_{s\in[t,T]}$ to the equation
\begin{equation}
\label{E:EquationX}
X_s \ = \ x + \int_t^s b(r,X_r) dr + \int_t^s \sigma(r,X_r) dW_r, \qquad t \leq s \leq T,\,\P\text{-a.s.}
\end{equation}
Furthermore, for any $p\geq1$ there exists a positive constant $C_p$, independent of $n$ and $(t,x)$, such that
\begin{equation}
\label{E:EstimateX^n}
\E\Big[\sup_{t \leq s \leq T}\big(|X_s^{n,t,x}|^p + |X_s^{t,x}|^p\big)\Big] \ \leq \ C_p\big(1 + |x|^p\big).
\end{equation}
Finally, for any $p\geq1$ we have
\begin{equation}
\label{E:X^n-->X}
\lim_{n\rightarrow\infty}\E\Big[\sup_{t \leq s \leq T}|X_s^{n,t,x} - X_s^{t,x}|^p\Big] \ = \ 0.
\end{equation}
\end{Proposition}
\textbf{Proof.}
It is well known that, under the assumptions of Proposition \ref{P:Estimate}, for any $n\in\N$ and $(t,x)\in[0,T]\times\R^d$, there exists a unique solution $(X_s^{n,t,x})_{s\in[t,T]}$ to equation \eqref{E:EquationX^n}, satisfying estimate \eqref{E:EstimateX^n}. Notice that the constant $C_p$ in \eqref{E:EstimateX^n} depends only on $T$, $p$, and the Lipschitz constants of $b_n$ and $\sigma_n$, which are uniformly bounded in $n$, so that $C_p$ does not depend on $n$. Now, we see that $b$ and $\sigma$ are Borel measurable and satisfy:
\begin{align*}
|b(t,x)-b(t,x')| + |\sigma(t,x)-\sigma(t,x')| \ &\leq \ C|x-x'|, \\
|b(t,0)| + |\sigma(t,0)| \ &\leq \ C,
\end{align*}
for all $t\in[0,T]$, $x,x'\in\R^d$. As a consequence, for any $(t,x)\in[0,T]\times\R^d$, there exists a unique solution $(X_s^{t,x})_{s\in[t,T]}$ to equation \eqref{E:EquationX}, satisfying estimate \eqref{E:EstimateX^n}. It remains to prove \eqref{E:X^n-->X}. Observe that
\[
X_s^{n,t,x} - X_s^{t,x} \ = \ \int_t^s \big(b_n(r,X_r^{n,t,x}) - b(r,X_r^{t,x})\big) dr + \int_t^s \big(\sigma_n(r,X_r^{n,t,x}) - \sigma(r,X_r^{t,x})\big) dW_r.
\]
Then, taking the $p$-th power, we get (recalling the standard inequality $(a+b)^p\leq2^{p-1}(a^p+b^p)$, for any $a,b\in\R$)
\begin{align*}
|X_s^{n,t,x} - X_s^{t,x}|^p \ &\leq \ 2^{p-1}\bigg|\int_t^s \big(b_n(r,X_r^{n,t,x}) - b(r,X_r^{t,x})\big) dr\bigg|^p \\
&\quad \ + 2^{p-1}\bigg|\int_t^s \big(\sigma_n(r,X_r^{n,t,x}) - \sigma(r,X_r^{t,x})\big) dW_r\bigg|^p.
\end{align*}
Taking the supremum over the time variable $s$, and applying H\"older's inequality to the drift term, we get (in the sequel we shall denote $c_p$ a generic positive constant, independent of $n$, depending only on $T$, $p$, and on the constant $C$ appearing in the statement of Proposition~\ref{P:Estimate})
\begin{align}
\label{E:SupProof2}
\sup_{t \leq s \leq T}|X_s^{n,t,x} - X_s^{t,x}|^p \ &\leq \ c_p\int_t^T \big|b_n(r,X_r^{n,t,x}) - b(r,X_r^{t,x})\big|^p dr \notag \\
&\quad \ + 2^{p-1}\sup_{t \leq s \leq T}\bigg|\int_t^s \big(\sigma_n(r,X_r^{n,t,x}) - \sigma(r,X_r^{t,x})\big) dW_r\bigg|^p.
\end{align}
Notice that
\begin{align}
\label{E:bn-b}
&\int_t^T \big|b_n(r,X_r^{n,t,x}) - b(r,X_r^{t,x})\big|^p dr \notag \\
&\leq \ 2^{p-1}\int_t^T \big|b_n(r,X_r^{n,t,x}) - b_n(r,X_r^{t,x})\big|^p dr + 2^{p-1}\int_t^T \big|b_n(r,X_r^{t,x}) - b(r,X_r^{t,x})\big|^p dr \notag \\
&\leq \ c_p\int_t^T \sup_{t \leq s \leq r}|X_s^{n,t,x} - X_s^{t,x}|^p dr + 2^{p-1}\int_t^T \big|b_n(r,X_r^{t,x}) - b(r,X_r^{t,x})\big|^p dr.
\end{align}
In addition, from Burkholder-Davis-Gundy inequality we have
\begin{align}
\label{E:sigman-sigma}
&\E\bigg[\sup_{t \leq s \leq T}\bigg|\int_t^s \big(\sigma_n(r,X_r^{n,t,x}) - \sigma(r,X_r^{t,x})\big) dW_r\bigg|^p\bigg] \notag \\
&\leq \ c_p \E\bigg[\int_t^T \big|\text{tr}\big((\sigma_n(r,X_r^{n,t,x})-\sigma(r,X_r^{t,x}))(\sigma_n\trans(r,X_r^{n,t,x})-\sigma\trans(r,X_r^{t,x}))\big)\big|^{\frac{p}{2}} dr\bigg] \notag \\
&\leq \ c_p \E\bigg[\int_t^T \big|\text{tr}\big((\sigma_n(r,X_r^{n,t,x})-\sigma_n(r,X_r^{t,x}))(\sigma_n\trans(r,X_r^{n,t,x})-\sigma_n\trans(r,X_r^{t,x}))\big)\big|^{\frac{p}{2}} dr\bigg] \notag \\
&\quad \ + c_p \E\bigg[\int_t^T \big|\text{tr}\big((\sigma_n(r,X_r^{t,x})-\sigma(r,X_r^{t,x}))(\sigma_n\trans(r,X_r^{t,x})-\sigma\trans(r,X_r^{t,x}))\big)\big|^{\frac{p}{2}} dr\bigg] \notag \\
&\leq \ c_p\int_t^T \sup_{t \leq s \leq r}|X_s^{n,t,x} - X_s^{t,x}|^p dr \notag \\
&\quad \ + c_p \E\bigg[\int_t^T \big|\text{tr}\big((\sigma_n(r,X_r^{t,x})-\sigma(r,X_r^{t,x}))(\sigma_n\trans(r,X_r^{t,x})-\sigma\trans(r,X_r^{t,x}))\big)\big|^{\frac{p}{2}} dr\bigg].
\end{align}
We remind that $\|A\|_{\text{tr}}:=\sqrt{\text{tr}(AA\trans)}$, with $A\in\R^{d\times d}$, defines a norm on the space of $d\times d$ matrices, which is indeed the Frobenius norm. Taking the expectation in \eqref{E:SupProof2}, and using \eqref{E:bn-b} and \eqref{E:sigman-sigma}, we find
\begin{align*}
&\E\Big[\sup_{t \leq s \leq T}|X_s^{n,t,x} - X_s^{t,x}|^p\Big] \\
&\leq \ c_p\int_t^T \E\Big[\sup_{t \leq s \leq r}|X_s^{n,t,x} - X_s^{t,x}|^p\Big] dr + c_p \int_t^T \E\big[\big|b_n(r,X_r^{t,x}) - b(r,X_r^{t,x})\big|^p\big] dr \\
&\quad \ + c_p \int_t^T \E\Big[\big|\text{tr}\big((\sigma_n(r,X_r^{t,x})-\sigma(r,X_r^{t,x}))(\sigma_n\trans(r,X_r^{t,x})-\sigma\trans(r,X_r^{t,x}))\big)\big|^{\frac{p}{2}}\Big] dr.
\end{align*}
Then, applying Gronwall's lemma to the map $r\mapsto\E[\sup_{t \leq s \leq r}|X_s^{n,t,x}-X_s^{t,x}|^p]$, we get
\begin{align*}
&\E\Big[\sup_{t \leq s \leq T}|X_s^{n,t,x} - X_s^{t,x}|^p\Big] \ \leq \ c_p \int_t^T \E\big[\big|b_n(r,X_r^{t,x}) - b(r,X_r^{t,x})\big|^p\big] dr \\
&+ c_p \int_t^T \E\Big[\big|\text{tr}\big((\sigma_n(r,X_r^{t,x})-\sigma(r,X_r^{t,x}))(\sigma_n\trans(r,X_r^{t,x})-\sigma\trans(r,X_r^{t,x}))\big)\big|^{\frac{p}{2}}\Big] dr.
\end{align*}
In conclusion, \eqref{E:X^n-->X} follows from Lebesgue's dominated convergence theorem.
\ep

\setcounter{equation}{0} \setcounter{Assumption}{0}
\setcounter{Theorem}{0} \setcounter{Proposition}{0}
\setcounter{Corollary}{0} \setcounter{Lemma}{0}
\setcounter{Definition}{0} \setcounter{Remark}{0}

\subsection{Limit theorem for BSDEs}

We prove a limit theorem for BSDEs designed for our purposes, which is inspired by the monotonic limit theorem of Peng \cite{peng00}, even if it is formulated under a different set of hypotheses. In particular, the monotonicity of the sequence $(Y^n)_n$ is not assumed. On the other hand, we impose a uniform boundedness for the sequence $(Y^n)_n$ in $\S^p(t,T)$ for some $p>2$, instead of $p=2$ as in \cite{peng00}. Furthermore, unlike \cite{peng00}, the terminal condition and the generator function of the BSDE solved by $Y^n$ are allowed to vary with $n$.

\begin{Proposition}
\label{P:LimitThmBSDE}
Let $(F_n)_n$ be a sequence of generator functions satisfying assumption \textup{(Aa)-(Ac)}, with the same constants $C_F$ and $M_F$ for all $n$. For any $n$, let $(Y^n,Z^n,K^n)\in\S^2(t,T)\times\H^2(t,T)^d\times\A^{+,2}(t,T)$, with $Y^n$ and $K^n$ continuous, satisfying, $\P$-a.s.,
\[
Y_s^n \ = \ Y_T^n + \int_s^T F_n(r,Y_r^n,Z_r^n) dr + K_T^n - K_s^n - \int_s^T \langle Z_r^n, dW_r\rangle, \qquad t \leq s \leq T
\]
and
\[
\|Y^n\|_{\S^2(t,T)}^2 + \|Z^n\|_{\H^2(t,T)^d} + \|K^n\|_{\S^2(t,T)} \ \leq \ C, \qquad \forall\,n\in\N,
\]
for some positive constant $C$, independent of $n$. Suppose that there exist a generator function $F$ satisfying conditions \textup{(Aa)-(Ac)} and a continuous process $Y\in\S^2(t,T)$, in addition $\sup_n\|Y^n\|_{\S^p(t,T)}<\infty$ for some $p>2$, and, for some null measurable sets $N_F\subset[t,T]\times\Omega$ and $N_Y\subset\Omega$,
\begin{align*}
F_n(s,\omega,y,z) \ &\overset{n\rightarrow\infty}{\longrightarrow} \ F(s,\omega,y,z), \qquad \forall\,(s,\omega,y,z)\in(([t,T]\times\Omega)\backslash N_F)\times\R\times\R^d, \\
Y_s^n(\omega) \ &\overset{n\rightarrow\infty}{\longrightarrow} \ Y_s(\omega), \qquad\qquad\;\; \forall\,(s,\omega)\in[t,T]\times(\Omega\backslash N_Y).
\end{align*}
Then, there exists a unique pair $(Z,K)\in\H^2(t,T)^d\times\A^{+,2}(t,T)$ such that, $\P$-a.s.,
\begin{equation}
\label{E:BSDELimit}
Y_s \ = \ Y_T + \int_s^T F(r,Y_r,Z_r) dr + K_T - K_s - \int_s^T \langle Z_r, dW_r\rangle, \qquad t \leq s \leq T.
\end{equation}
In addition, $Z^n$ converges strongly $($resp. weakly$)$ to $Z$ in $\L^q(t,T;\R^d)$ $($resp. $\H^2(t,T)^d$$)$, for any $q\in[1,2[$, and $K_\tau^n$ converges weakly to $K_\tau$ in $L^2(\Omega,\Fc_\tau,\P)$, for any stopping time $\tau$ valued in $[t,T]$.
\end{Proposition}
\begin{Remark}
\label{R:Y_Sp<infty}
{\rm
Notice that, under the hypotheses of Proposition \ref{P:LimitThmBSDE} (more precisely, given that $Y$ is continuous, $\sup_n\|Y^n\|_{\S^p(t,T)}<\infty$ for some $p>2$, $Y_s^n(\omega) \rightarrow Y_s(\omega)$ as $n$ tends to infinity for all $(s,\omega)\in[t,T]\times(\Omega\backslash N_Y)$), it follows that $\|Y\|_{\S^p(t,T)}<\infty$. Indeed, from Fatou's lemma we have
\begin{equation}
\label{E:Y_Sp<infty1}
\E\Big[\liminf_{n\rightarrow\infty}\sup_{t\leq s\leq T}|Y_s^n|^p\Big] \ \leq \ \liminf_{n\rightarrow\infty}\|Y^n\|_{\S^p(t,T)}^p \ < \ \infty.
\end{equation}
Moreover, since $Y$ is continuous, there exists a null measurable set $N_Y'\subset\Omega$ such that $s\mapsto Y_s(\omega)$ is continuous on $[t,T]$ for every $\omega\in\Omega\backslash N_Y'$. Then, for any $\omega\in\Omega\backslash(N_Y\cup N_Y')$, there exists $\tau(\omega)\in[t,T]$ such that
\begin{equation}
\label{E:Y_Sp<infty2}
\sup_{t \leq s \leq T}|Y_s(\omega)|^p \ = \ |Y_{\tau(\omega)}(\omega)|^p \ = \ \lim_{n\rightarrow\infty} |Y_{\tau(\omega)}^n(\omega)|^p \ \leq \ \liminf_{n\rightarrow\infty} \sup_{t\leq s\leq T}|Y_s^n(\omega)|^p.
\end{equation}
Therefore, combining \eqref{E:Y_Sp<infty1} with \eqref{E:Y_Sp<infty2}, we end up with $\|Y\|_{\S^p(t,T)}<\infty$.
\ep
}
\end{Remark}
\textbf{Proof.}
We begin proving the uniqueness of $(Z,K)$. Let $(Z,K),(Z',K')\in\H^2(t,T)^d\times\A^{+,2}(t,T)$ be two pairs satisfying \eqref{E:BSDELimit}. Taking the difference and rearranging the terms, we obtain
\[
\int_s^T \langle Z_r-Z_r',dW_r\rangle \ = \ \int_s^T \big(F(r,Y_r,Z_r)-F(r,Y_r,Z_r')\big) dr + K_T-K_s - (K_T'-K_s').
\]
Now, the right-hand side has finite variation, while the left-hand side has not finite variation, unless $Z=Z'$. This implies $Z=Z'$, from which we deduce $K=K'$.

The rest of the proof is devoted to the existence of $(Z,K)$ and it is divided in different steps.\\
\emph{Step 1. Limit BSDE.} From the hypotheses, we see that there exists a positive constant $c$, independent of $n$, such that
\[
\E\int_t^T |F_n(r,Y_r^n,Z_r^n)|^2 dr \ \leq \ c, \qquad \forall\,n\in\N.
\]
It follows that the sequence $(Z_\cdot^n,F_n(\cdot,Y_\cdot^n,Z_\cdot^n))_n$ is bounded in the Hilbert space $\H^2(t,T)^d\times\L^2(t,T;\R)$. Therefore, there exists a subsequence $(Z_\cdot^{n_k},F_{n_k}(\cdot,Y_\cdot^{n_k},Z_\cdot^{n_k}))_k$ which converges weakly to some $(Z,G)\in\H^2(t,T)^d\times\L^2(t,T;\R)$. This implies that, for any stopping time $\tau\in[t,T]$, the following weak convergences hold in $L^2(\Omega,\Fc_\tau,\P)$ as $k\rightarrow\infty$:
\[
\int_t^\tau F_{n_k}(r,Y_r^{n_k},Z_r^{n_k}) dr \ \rightharpoonup \ \int_t^\tau G(r) dr, \qquad\qquad
\int_t^\tau \langle Z_r^{n_k},dW_r\rangle \ \rightharpoonup \ \int_t^\tau \langle Z_r,dW_r\rangle.
\]
Since
\[
K_\tau^n \ = \ Y_t^n - Y_\tau^n - \int_t^\tau F_n(r,Y_r^n,Z_r^n) dr + \int_t^\tau \langle Z_r^n,dW_r\rangle
\]
and, by hypothesis, $Y_\tau^n \rightarrow Y_\tau$ strongly in $L^2(\Omega,\Fc_\tau,\P)$, we also have the weak convergence, as $k\rightarrow\infty$,
\begin{equation}
\label{E:K_tau=Y}
K_\tau^{n_k} \ \rightharpoonup \ K_\tau \ = \ \tilde K_{t,\tau},
\end{equation}
where
\[
\tilde K_{t,s} \ := \ Y_t - Y_s - \int_t^s G(r) dr + \int_t^s \langle Z_r,dW_r\rangle, \qquad t\leq s\leq T.
\]
Notice that $(\tilde K_{t,s})_{t\leq s\leq T}$ is adapted and continuous, so that it is a predictable process. We have that $\E[|K_T|^2] < \infty$. Moreover, $K^{n_k}$ converges weakly to $K$ in the Hilbert space $\L^2(t,T;\R)$. Indeed, let $\xi\in\L^2(t,T;\R)$; then, by Fubini's theorem,
\[
\E\bigg[\int_t^T \xi_s(K_s^{n_k} - K_s) ds\bigg] \ = \ \int_t^T \E\big[\xi_s(K_s^{n_k} - K_s)\big] ds.
\]
Since $\xi_s\in L^2(\Omega,\Fc_s,\P)$, for a.e. $s\in[t,T]$, we conclude, from Lebesgue's dominated convergence theorem,
\[
\int_t^T \E\big[\xi_s(K_s^{n_k} - K_s)\big] ds \ \overset{k\rightarrow\infty}{\longrightarrow} \ 0.
\]
This implies that $K$ is a predictable process. Since the process on the right-hand side of \eqref{E:K_tau=Y} is also predictable and they are equal for all stopping times valued in $[t,T]$ (it would be enough to consider the predictable ones), it follows from the predictable section theorem (see, e.g., Theorem 86, Chapter IV, in \cite{dellmeyer75}) that they are indistinguishable. In particular, $K$ is a continuous process.

Let us prove that $K$ is a nondecreasing process. For any pair $r,s$ with $t \leq r \leq s \leq T$, we have $K_r \leq K_s$, $\P$-almost surely. Indeed, let $\xi\in L^2(\Omega,\Fc_s,\P)$ be nonnegative, then, from the martingale representation theorem, we see that there exist a random variable $\zeta\in L^2(\Omega,\Fc_r,\P)$ and an $\F$-predictable square integrable process $\eta$ such that
\[
\xi \ = \ \zeta + \int_r^s \eta_u dW_u.
\]
Therefore
\begin{align*}
0 \ &\leq \ \E[\xi(K_s^n-K_r^n)] \ = \ \E[\xi K_s^n] - \E[\zeta K_r^n] - \E\bigg[\E\bigg[K_r^n\int_r^s \eta_u dW_u \bigg|\Fc_r\bigg]\bigg] \\
&= \ \E[\xi K_s^n] - \E[\zeta K_r^n] \ \overset{n\rightarrow\infty}{\longrightarrow} \ \E[\xi K_s] - \E[\zeta K_r] \ = \ \E[\xi(K_s-K_r)],
\end{align*}
which shows that $K_r \leq K_s$, $\P$-almost surely. As a consequence, there exists a null measurable set $N\subset\Omega$ such that $K_r(\omega) \leq K_s(\omega)$, for all $\omega\in\Omega\backslash N$, with $r,s\in\Q\cap[0,T]$, $r<s$. Then, from the continuity of $K$ it follows that it is a nondecreasing process, so that $K\in\A^{+,2}(t,T)$.

Finally, we notice that the process $Z$ in expression \eqref{E:K_tau=Y} is uniquely determined, as it can be seen identifying the Brownian parts and the finite variation parts in \eqref{E:K_tau=Y}. Thus, not only the subsequence $(Z^{n_k})_k$, but all the sequence $(Z^n)_n$ converges weakly to $Z$ in $\H^2(t,T)^d$. It remains to show that $G(r)$ in \eqref{E:K_tau=Y} is actually $F(r,Y_r,Z_r)$.\\
\emph{Step 2. Strong convergence of $(Z^n)_n$.} Let $\alpha\in(0,1)$ and consider the function $h_\alpha(y) = |\min(y-\alpha,0)|^2$, $y\in\R$. By applying Meyer-It\^o's formula combined with the occupation times formula (see, e.g., Theorem 70 and Corollary 1, Chapter IV, in \cite{protter05}) to $h_\alpha(Y_s^n-Y_s)$ between $t$ and $T$, observing that the second derivative of $h_\alpha$ in the sense of distributions is a $\sigma$-finite Borel measure on $\R$ absolutely continuous to the Lebesgue measure with density $1_{]-\infty,\alpha[}(\cdot)$, we obtain
\begin{align*}
&\E\big[|\min(Y_t^n-Y_t-\alpha,0)|^2\big] + \E\int_t^T 1_{\{Y_s^n-Y_s<\alpha\}} |Z_s^n-Z_s|^2 ds \\
&= \ \E\big[|\min(Y_T^n-Y_T-\alpha,0)|^2\big] + 2\E\int_t^T \min(Y_s^n-Y_s-\alpha,0) \big(F_n(s,Y_s^n,Z_s^n) - G(s)\big)ds \\
&\quad \ + 2\E\int_t^T \min(Y_s^n - Y_s-\alpha,0) dK_s^n - 2\E\int_t^T \min(Y_s^n - Y_s-\alpha,0) dK_s.
\end{align*}
Since $\min(Y_s^n - Y_s-\alpha,0) dK_s^n \leq 0$, we get
\begin{align}
\label{E:min<K}
&\E\int_t^T 1_{\{Y_s^n-Y_s<\alpha\}} |Z_s^n-Z_s|^2 ds \ \leq \ \E\big[|\min(Y_T^n-Y_T-\alpha,0)|^2\big] \\
&+ 2\E\int_t^T \min(Y_s^n-Y_s-\alpha,0) \big(F_n(s,Y_s^n,Z_s^n) - G(s)\big)ds - 2\E\int_t^T \min(Y_s^n - Y_s-\alpha,0) dK_s. \notag
\end{align}
Let us study the behavior of the right-hand side of \eqref{E:min<K} as $n$ goes to infinity. We begin noting that
\begin{equation}
\label{E:ProofConvergence1}
\E\big[|\min(Y_T^n-Y_T-\alpha,0)|^2\big] \ \overset{n\rightarrow\infty}{\longrightarrow} \ \alpha^2.
\end{equation}
Regarding the second-term on the right-hand side of \eqref{E:min<K}, since the sequence $(F_n(\cdot,Y_\cdot^n,Z_\cdot^n) - G(\cdot))_n$ is bounded in $\L^2(t,T;\R)$, we have
\[
\sup_{n\in\N}\bigg(\E\bigg[\int_t^T|F_n(s,Y_s^n,Z_s^n) - G(s)|^2ds\bigg]\bigg)^{\frac{1}{2}} \ =: \ \bar c \ < \ \infty.
\]
Therefore, by Cauchy-Schwarz inequality we find
\begin{align}
\label{E:ProofConvergence2}
&\E\int_t^T |\min(Y_s^n-Y_s-\alpha,0)| |F_n(s,Y_s^n,Z_s^n) - G(s)|ds \notag \\
&\leq \bar c \bigg(\E\bigg[\int_t^T|\min(Y_s^n-Y_s-\alpha,0)|^2ds\bigg]\bigg)^{\frac{1}{2}} \  \overset{n\rightarrow\infty}{\longrightarrow} \ \bar c\sqrt{T-t}\,\alpha.
\end{align}
Concerning the last term on the right-hand side of \eqref{E:min<K}, we notice that, by hypothesis and Remark \ref{R:Y_Sp<infty}, there exists some $p>2$ such that, from Cauchy-Schwarz inequality,
\begin{align*}
&\sup_{n\in\N}\E\bigg[\int_t^T |\min(Y_s^n - Y_s-\alpha,0)|^{\frac{p}{2}} dK_s\bigg] \\
&\leq \ \sup_{n\in\N}\bigg(\E\Big[\sup_{t \leq s \leq T} |\min(Y_s^n - Y_s-\alpha,0)|^p \Big]\bigg)^{\frac{1}{2}} \big(\E\big[|K_T|^2\big]\big)^{\frac{1}{2}} \ < \infty.
\end{align*}
It follows that $(\min(Y_\cdot^n - Y_\cdot-\alpha,0))_n$ is a uniformly integrable sequence on $([t,T]\times\Omega,\Bc([t,T])\otimes\Fc,dK_s\otimes d\P)$. Moreover, by assumption, there exists a null measurable set $N_Y\subset\Omega$ such that $Y_s^n(\omega)$ converges to $Y_s(\omega)$, for any $(s,\omega)\notin[t,T]\times N_Y$. Notice that $dK_s\otimes d\P([t,T]\times N_Y)=0$, therefore $Y^n$ converges to $Y$ pointwisely a.e. with respect to $dK_s\otimes d\P$. This implies that
\begin{equation}
\label{E:ProofConvergence3}
\E\bigg[\int_t^T |\min(Y_s^n - Y_s-\alpha,0)| dK_s\bigg] \ \overset{n\rightarrow\infty}{\longrightarrow} \ \alpha\E[K_T].
\end{equation}
From the convergence results \eqref{E:ProofConvergence1}, \eqref{E:ProofConvergence2}, and \eqref{E:ProofConvergence3}, we end up with
\begin{equation}
\label{E:min<K2}
\limsup_{n\rightarrow\infty}\E\int_t^T 1_{\{Y_s^n-Y_s<\alpha\}} |Z_s^n-Z_s|^2 ds \ \leq \ \alpha^2 + 2\bar c\sqrt{T-t}\,\alpha + 2\alpha\E[K_T].
\end{equation}
From Egoroff's theorem, for any $\delta>0$ there exists a measurable set $A\subset[t,T]\times\Omega$, with $ds\otimes d\P(A)<\delta$, such that $(Y^n)_n$ converges uniformly to $Y$ on $([t,T]\times\Omega)\backslash A$. In particular, for any $\alpha\in]0,1[$ we have $|Y_s^n(\omega)-Y_s(\omega)| < \alpha$, for all $(s,\omega)\in([t,T]\times\Omega)\backslash A$, whenever $n$ is large enough. Therefore, from \eqref{E:min<K2} we get
\begin{align*}
&\limsup_{n\rightarrow\infty}\E\int_t^T 1_{([t,T]\times\Omega)\backslash A} |Z_s^n-Z_s|^2 ds \ = \ \limsup_{n\rightarrow\infty}\E\int_t^T 1_{([t,T]\times\Omega)\backslash A} 1_{\{Y_s^n-Y_s<\alpha\}} |Z_s^n-Z_s|^2 ds \notag \\
&\leq \ \limsup_{n\rightarrow\infty}\E\int_t^T 1_{\{Y_s^n-Y_s<\alpha\}} |Z_s^n-Z_s|^2 ds \ \leq \ \alpha^2 + 2\bar c\sqrt{T-t}\,\alpha + 2\alpha\E[K_T].
\end{align*}
Sending $\alpha\rightarrow0^+$, we obtain
\begin{equation}
\label{E:alpha-->0+}
\lim_{n\rightarrow\infty}\E\int_t^T 1_{([t,T]\times\Omega)\backslash A} |Z_s^n-Z_s|^2 ds \ = \ 0.
\end{equation}
Now, let $q\in[1,2[$; by H\"older's inequality,
\begin{align*}
&\E\int_t^T |Z_s^n-Z_s|^q ds \ = \ \E\int_t^T 1_{([t,T]\times\Omega)\backslash A} |Z_s^n-Z_s|^q ds + \E\int_t^T 1_A |Z_s^n-Z_s|^q ds \\
&\leq \ \bigg(\E\int_t^T 1_{([t,T]\times\Omega)\backslash A} |Z_s^n-Z_s|^2 ds\bigg)^{\frac{q}{2}}(T-t)^{\frac{2-q}{2}} + \bigg(\E\int_t^T |Z_s^n-Z_s|^2 ds\bigg)^{\frac{q}{2}}\delta^{\frac{2-q}{2}}.
\end{align*}
Since the sequence $(Z^n)_n$ is bounded in $\H^2(t,T)^d$, we have
\[
\sup_{n\in\N}\E\int_t^T |Z_s^n-Z_s|^2 ds \ =: \ \hat c < \infty.
\]
Therefore
\[
\E\int_t^T |Z_s^n-Z_s|^q ds \ \leq \ \bigg(\E\int_t^T 1_{([t,T]\times\Omega)\backslash A} |Z_s^n-Z_s|^2 ds\bigg)^{\frac{q}{2}}(T-t)^{\frac{2-q}{2}} + \hat c^{\frac{q}{2}}\delta^{\frac{2-q}{2}},
\]
which implies, by \eqref{E:alpha-->0+},
\[
\limsup_{n\rightarrow\infty}\E\int_t^T |Z_s^n-Z_s|^q ds \ \leq \ \hat c^{\frac{q}{2}}\delta^{\frac{2-q}{2}}.
\]
Sending $\delta\rightarrow0^+$ we deduce the strong convergence of $Z^n$ towards $Z$ in $\L^q(t,T;\R^d)$, for any $q\in[1,2[$.

Notice that, for any $q\in[1,2[$, we have (recalling the standard inequality $(x+y)^q \leq 2^{q-1}(x^q+y^q)$, for any $x,y\in\R_+$)
\begin{align*}
\E\bigg[\int_t^T |F_n(s,Y_s^n,Z_s^n) - F(s,Y_s,Z_s)|^q ds\bigg] \ &\leq \ 2^{q-1}\E\bigg[\int_t^T |F_n(s,Y_s^n,Z_s^n) - F_n(s,Y_s,Z_s)|^q ds\bigg] \\
&\;\;\; + 2^{q-1}\E\bigg[\int_t^T |F_n(s,Y_s,Z_s) - F(s,Y_s,Z_s)|^q ds\bigg].
\end{align*}
Therefore, by the uniform Lipschitz condition on $F_n$ with respect to $(y,z)$, and the convergence of $F_n$ towards $F$, we deduce the strong convergence of $(F_n(\cdot,Y_\cdot^n,Z_\cdot^n))_n$ to $F(\cdot,Y_\cdot,Z_\cdot)$ in $\L^q(t,T;\R)$, $q\in[1,2[$. Since $G(\cdot)$ is the weak limit of $(F_n(\cdot,Y_\cdot^n,Z_\cdot^n))_n$ in $\L^2(t,T;\R)$, we deduce that $G(\cdot)$ $=$ $F(\cdot,Y_\cdot,Z_\cdot)$. In conclusion, the triple $(Y,Z,K)$ solves the backward stochastic differential equation \eqref{E:BSDELimit}.
\ep

\setcounter{equation}{0} \setcounter{Assumption}{0}
\setcounter{Theorem}{0} \setcounter{Proposition}{0}
\setcounter{Corollary}{0} \setcounter{Lemma}{0}
\setcounter{Definition}{0} \setcounter{Remark}{0}

\subsection{An additional result in real analysis}

\begin{Lemma}
\label{L:StabilityApp}
Let $(f_{n,k})_{n,k\in\N}$, $(f_n)_{n\in\N}$, and $f$ be $\R^q$-valued functions on $[0,T]\times X$, where $(X,\|\cdot\|)$ is a normed space which contains a countable dense subset $E$, and
\[
f_{n,k}(t,x) \ \overset{k\rightarrow\infty}{\longrightarrow} \ f_n(t,x), \qquad f_n(t,x) \ \overset{n\rightarrow\infty}{\longrightarrow} \ f(t,x), \qquad \forall\,(t,x)\in[0,T]\times X.
\]
Suppose that the double sequence $(f_{n,k})_{n,k\in\N}$ is  locally equicontinuous. Then, there exists a subsequence $(f_{n,k(n)})_n$ that converges pointwisely to $f$ on $[0,T]\times X$.
\end{Lemma}
\textbf{Proof.}
Let $(t_1,x_1),(t_2,x_2),(t_3,x_3),\ldots$ be an enumeration of the points of $(\Q\cap[0,T])\times E$. For any $n$ and $j$, it follows from the convergence $f_{n,k}(t_j,x_j) \rightarrow f_n(t_j,x_j)$, as $k\rightarrow\infty$, that there exists a positive integer $K_{n,j}$ such that
\[
|f_{n,k}(t_j,x_j) - f_n(t_j,x_j)| \ < \ \frac{1}{n}, \qquad \text{if }k \geq K_{n,j}.
\]
Let $k(n) = k(n-1)\vee K_{n,1}\vee\cdots\vee K_{n,n}$, $n\in\N$, with $k(-1)=0$. Then, for each $j$, $f_{n,k(n)}(t_j,x_j)\rightarrow f(t_j,x_j)$, as $n\rightarrow\infty$.

Now, take $(t,x)\in[0,T]\times X$ and $\eps>0$. Consider $R>0$ such that $\|x\|<R$. By the local equicontinuity, there exists $\delta>0$, depending only on $\eps$ and $R$, so that $|t-s|,\|x-y\|<\delta$, $\|y\|<R$, implies $|f(t,x)-f(s,y)|<\eps/3$ and $|f_{n,m}(t,x)-f_{n,m}(s,y)|<\eps/3$, for all $n,m$. Let $(t_j,x_j)\in(\Q\cap[0,T])\times E$ be such that $|t-t_j|,\|x-x_j\|<\delta$ and $\|x_j\|<R$. We know that there exists a positive integer $N$ for which $|f_{n,k(n)}(t_j,x_j) - f(t_j,x_j)|<\eps/3$, for any $n \geq N$. Therefore
\begin{align*}
&|f_{n,k(n)}(t,x) - f(t,x)| \\
&\leq \ |f_{n,k(n)}(t,x) - f_{n,k(n)}(t_j,x_j)| + |f_{n,k(n)}(t_j,x_j) - f(t_j,x_j)| + |f(t_j,x_j) - f(t,x)| \ < \ \eps,
\end{align*}
for all $n \geq N$.
\ep

\vspace{5mm}

\noindent\textbf{Acknowledgements.}
The present work was partially supported
by the ANR Project MASTERIE 2010 BLAN 0121 01.

\small
\bibliographystyle{plain}
\bibliography{biblio}

\end{document}